\documentclass{amsart}
\usepackage{graphicx,amsmath,amssymb,amsfonts,amscd,color}
\usepackage[all,cmtip,arrow,2cell]{xy}
\usepackage{mathrsfs}
\usepackage[colorlinks,
             linkcolor=black,
             citecolor=blue,
             bookmarksnumbered=true]{hyperref}
\usepackage{xrtwo}
\newtheorem{thm}{Theorem}[section]
\newtheorem{cor}[thm]{Corollary}
\newtheorem{lem}[thm]{Lemma}
\newtheorem{prop}[thm]{Proposition}

\theoremstyle{definition}
\newtheorem{defn}[thm]{Definition}
\theoremstyle{definition}
\newtheorem{eg}[thm]{Example}
\theoremstyle{remark}
\newtheorem{rem}[thm]{Remark}
\theoremstyle{remark}
\newtheorem*{notation}{Notation}
\numberwithin{equation}{section}
\newcommand{\abs}[1]{\left\vert#1\right\vert}
\newcommand{\set}[1]{\left\{#1\right\}}

\newcommand{\E}{\mathcal E}

\newcommand{\W}{\mathcal W}
\newcommand{\X}{\mathcal X}
\newcommand{\Y}{\mathcal Y}
\newcommand{\M}{\mathcal M}
\newcommand{\G}{\mathcal G}
\newcommand{\A}{\mathcal{A}}
\newcommand{\B}{\mathcal{B}}
\newcommand{\C}{\mathcal{C}}
\newcommand{\D}{\mathcal D}

\renewcommand{\P}{\mathcal{P}}
\newcommand{\R}{\mathcal{R}}
\newcommand{\Q}{\mathcal{Q}}
\newcommand{\J}{\mathcal{J}}

\renewcommand{\H}{\mathcal{H}}
\renewcommand{\S}{\mathcal{S}}

\newcommand{\RR}{\mathbb R}
\newcommand{\QQ}{\mathbb Q}
\newcommand{\KK}{\mathbb K}
\newcommand{\ZZ}{\mathbb Z}
\newcommand{\NN}{\mathbb N}
\renewcommand{\SS}{\mathbb S}
\newcommand{\CC}{\mathbb C}
\newcommand{\EE}{\mathbb E}
\newcommand{\ee}{\mathfrak e}

\newcommand{\Mm}{{\mathbb M}_\mathrm{m}}
\newcommand{\mM}{{{\mathcal M}_\mathrm{m}}}
\newcommand{\Gm}{{\mathbb G}_\mathrm{m}}
\newcommand{\mG}{{{\mathcal G}_\mathrm{m}}}
\newcommand{\Gad}{{\mathbb G}_\mathrm{ad}}
\newcommand{\adG}{{\mathcal G}_\mathrm{ad}}
\newcommand{\hMm}{{\hat{\mathbb M}_\mathrm{m}}}
\newcommand{\hmM}{{\hat{\mathcal M}_\mathrm{m}}}
\newcommand{\hGm}{\hat{\mathbb G}_\mathrm{m}}
\newcommand{\hmG}{{\hat{\mathcal G}_\mathrm{m}}}
\newcommand{\hGad}{\hat{\mathbb G}_\mathrm{ad}}
\newcommand{\hadG}{\hat{\mathcal G}_\mathrm{ad}}
\newcommand{\alg}{{\mathrm{alg}}}
\newcommand{\eps}{\varepsilon}
\newcommand{\To}{\longrightarrow}
\newcommand{\oT}{\longleftarrow}

\newcommand{\bC}{\mathbf{C}}
\newcommand{\bD}{\mathbf{D}}
\newcommand{\bE}{\mathbf{E}}

\newcommand{\bSE}{\mathbf{SE}}
\newcommand{\bT}{\mathbf{T}}
\newcommand{\bS}{\mathbf{S}}
\newcommand{\bR}{\mathbf{R}}

\newcommand{\bH}{\mathbf{H}}

\newcommand{\bcinf}{\mathbf{\cinf}}
\newcommand{\bscinf}{\mathbf{S}\mathbf{\cinf}}
\newcommand{\bcomega}{\mathbf{\comega}}
\newcommand{\bcom}{\mathbf{Com}}
\newcommand{\bscom}{\mathbf{SCom}}

\newcommand{\bDelta}{\mathbf{\Delta}}
\newcommand{\even}{{\underline{0}}}
\newcommand{\odd}{{\underline{1}}}
\newcommand{\cinf}{{C^\infty}}
\newcommand{\comega}{{C^\omega}}

\newcommand{\Set}{\mathbf{Set}}
\newcommand{\Mon}{\mathbf{Mon}}
\newcommand{\Grp}{\mathbf{Grp}}
\newcommand{\End}{\mathbf{End}}

\newcommand{\Ch}{\mathbf{Ch}}
\newcommand{\Ab}{\mathbf{Ab}}
\newcommand{\Mod}{\mathbf{Mod}}

\newcommand{\Alg}{\mathbf{Alg}}

\newcommand{\comsalg}{\mathbf{SComAlg}}

\newcommand{\Sym}{\mathrm{Sym_\KK}}
\newcommand{\id}{\mathrm{id}}
\newcommand{\ev}{\mathrm{ev}}
\newcommand{\op}{\mathrm{op}}
\newcommand{\cl}{\mathrm{cl}}
\newcommand{\fg}{\mathrm{fg}}
\newcommand{\dg}{\mathbf{dg}\mbox{-}}
\newcommand{\gr}{\mathbf{gr}\mbox{-}}
\newcommand{\dd}{\mathbf{d}\mbox{-}}
\newcommand{\triv}{\mathrm{triv}}

\newcommand{\Hom}{\operatorname{Hom}}
\newcommand{\Der}{\operatorname{Der}}
\newcommand{\Ker}{\operatorname{Ker}}
\newcommand{\Spec}{\mathbf{Spec}}
\renewcommand{\Im}{\operatorname{Im}}
\newcommand{\del}{\partial}
\newcommand{\om}{\Omega}
\newcommand{\wom}{\widehat{\Omega}}

\newcommand{\e}{\left[\epsilon\right]}
\newcommand{\ep}{\left[\eps\right]}

\def\oinft{ \mbox{\,\put(-1.5,0){$\bigcirc$}\put(-0.2,1.05){\hbox{\tiny$\infty$}}}\mspace{16mu}}
\def\oi{\mspace{3mu}\resizebox{0.28cm}{!}{$\oinft$}\mspace{3mu}}
\def\oinfty{\raisebox{.24ex}{$\oi$}}
\def\ve{\KK\mbox{-}\Mod}
\def\sv{\left(\ve\right)^{\mathbb{Z}_2}}
\newcommand{\IHom}{\mathbb{H}\!\operatorname{om}}
\newcommand{\ii}{\mathbf{i}}
\def\longlongleftarrow{\longleftarrow\!\!\!-\!\!\!-\!\!\!-\!\!\!-\!\!\!-\!\!\!-}
\def\longlongrightarrow{-\!\!\!-\!\!\!-\!\!\!-\!\!\!-\!\!\!-\!\!\!\longrightarrow}
\def\longlonglongrightarrow{-\!\!\!-\!\!\!-\!\!\!-\!\!\!-\!\!\!-\!\!\!-\!\!\!\longrightarrow}
\newcommand{\Adj}[4]{\xymatrix@1{#2 \ar@<-0.5ex>[r]_-{#4} & #3 \ar@<-0.5ex>[l]_-{#1}}}
\newcommand{\Adjadj}[7]{\xymatrix@1{#2 \ar@<-0.5ex>[r]_-{#4} & #3 \ar@<-0.5ex>[l]_-{#1} \ar@<-0.5ex>[r]_-{#7} & #6 \ar@<-0.5ex>[l]_-{#5} }}
\newcommand{\Adjlong}[4]{\xymatrix@C=2cm{#2 \ar@<-0.65ex>[r]_-{#4} & #3 \ar@<-0.65ex>[l]_-{#1}}}

\newcommand{\cop}[2]{ #1_{\diamond #2}}
\newcommand{\ucolim}{\underrightarrow{\lim}\,}
\newcommand{\ulim}{\underleftarrow{\lim}\,}
\begin{document}

\title{Homological Algebra for Superalgebras of Differentiable Functions}%
\author{David Carchedi}
\address{D. Carchedi \hspace{16pt}\mbox{Max Planck Institute for Mathematics, Bonn, Germany. }
}
\author{Dmitry Roytenberg}
\address{D. Roytenberg\hspace{5pt}
Department of Mathematics,
Utrecht University,
The Netherlands.}
\keywords{homotopy, differential graded, derived manifold}%

\date{\today}
%\dedicatory{}%
%\commby{}%
% ----------------------------------------------------------------
\begin{abstract}
This is the second in a series of papers laying the foundations for a differential graded approach to derived differential geometry (and other geometries in characteristic zero). In this paper, we extend the classical notion of a dg-algebra to define, in particular, the notion of a differential graded algebra in the world of $\bcinf$-rings. The opposite of the category of differential graded $\bcinf$-algebras contains the category of differential graded manifolds as a full subcategory. More generally, this notion of differential graded algebra makes sense for algebras over any (super) Fermat theory, and hence one also arrives at the definition of a differential graded algebra appropriate for the study of derived real and complex analytic manifolds and other variants. We go on to show that, for any super Fermat theory $\bS$ which admits integration, a concept we define and show is satisfied by all important examples, the category of differential graded $\bS$-algebras supports a Quillen model structure naturally extending the classical one on differential graded algebras, both in the bounded and unbounded case (as well as differential algebras with no grading). Finally, we show that, under the same assumptions, any of these categories of differential graded $\bS$-algebras have a simplicial enrichment, compatible in a suitable sense with the model structure.
%We construct Quillen model structures on several categories of $\cinf$-superalgebras with a dg-structure.
\end{abstract}

%Let $\bcom$ be the Lawvere theory of commutative algebras over $\QQ$, $\bE$ a Fermat extension of $\bcom$, $\bSE$
%its superization. After generalizing the notion of a differential graded structure to the setting of $\bE$-superalgebras ($\bSE$-algebras), we construct a Quillen model structure on the category of dg $\bE$-algebras generalizing the well-known one on the category of commutative dg-algebras over a ground ring $\KK$ containing $\QQ$. This leads to the ``dg'' approach to derived geometry for many classical geometries, such as $\cinf$ differential geometry, real or complex analytic geometry, and differential geometry over non-standard reals.

\maketitle
% ----------------------------------------------------------------
\section{Introduction.}
The purpose of this paper is to introduce the theory of differential graded algebras for a super Fermat theory, hence extending homological algebra to this setting. Super Fermat theories are theories of supercommutative algebras in which, in addition to evaluating polynomials on elements, one can evaluate other classes of infinitely differentiable functions. For details, we refer the reader to \cite{dg1}. Of central importance to our future work is that, in particular, we define in this paper the notion of a differential graded $\bcinf$-algebra, which will play a crucial role in our differential graded approach to derived smooth manifolds. For any super Fermat theory $\bS$, we define the category of differential graded $\bS$-algebras and, if $\bS$ admits integration, a concept which we define and which is satisfied in all important examples, we show it can be endowed with a Quillen model structure with an appropriately compatible simplicial enrichment. In a future paper, we will exploit this to give a model for derived smooth manifolds using differential graded manifolds which is directly amenable to calculations.

The key idea of our approach is to exploit the intimate connection between supercommutativity and classical differential graded algebras, and to generalize it to the setting of super Fermat theories. It is well known that a differential on a supercommutative algebra $\A$ , i.e. an odd derivation squaring to zero, corresponds to an action of the odd additive group $\Gad^\odd$ on the associated affine superscheme. Moreover, a grading on $\A$ corresponds to an action of the multiplicative group $\Gm,$ which if lifts to the multiplicative monoid $\Mm,$ is a grading by non-negative weights. Having the structure of a differential graded supercommutative algebra is the same as having an action of the semi-direct product $$\EE^\times\cong\Gm\ltimes\Gad^\odd,$$ which is the automorphism group of the odd line. (If the action lifts to the endomorphism monoid $\EE$, only non-positive weights occur in the grading (using cochain complex conventions), and dually, only non-negative weights will occur if the induced action of $\left(\EE^{\times}\right)^{\op}$ lifts to $\EE^{\op}$). Finally, a differential graded \emph{commutative} algebra is the same as a differential graded supercommutative algebra whose integer grading is compatible with its $\ZZ_2$-grading, which may be arranged by imposing that the action of $-1\in\Gm$ coincide with the parity involution (we call such $\EE^\times$-actions \emph{even}).

The previous paragraph provides an algebraic analogue of the observation that differential (non-positively or non-negatively) graded manifolds are the same as supermanifolds with a (left or right) action of the endomorphism monoid of the odd line, a statement made first by Kontsevich \cite{defquant} and later amplified by \v{S}evera \cite{SevSimpl}. The main idea behind our approach to defining differential graded algebras for a super Fermat theory $\bS,$ is to use appropriate versions of $\Gad^\odd,$ $\Gm,$ and $\Mm$ in this setting. In particular, in the case of $\bcinf$-rings, the monoid corresponding to $\EE$ is literally the endomorphism monoid of the smooth odd line $\End\left(\RR^{0|1}\right),$ so in particular, the opposite category of differential graded $\bcinf$-algebras contains the category of differential graded manifolds as a full subcategory. Since we show the former is simplicially enriched, this yields a simplicial category of differential graded manifolds.

Our approach to derived differential geometry is different than the existing approaches \cite{spivak,joyce,borisovnoel} as it based upon differential graded (dg) geometry, making it closer in spirit to the dg approach to derived \emph{algebraic} geometry pioneered by Ciocan-Fontanine and Kapranov in \cite{kapbg,quot,virtual}, and to later work of Behrend \cite{dgs1,dgs2}; however, it was already observed quite early on by Dold and Puppe \cite{dold,doldpuppe}, and expanded upon by Quillen in \cite{cohrings}, that such a dg approach misbehaves for general commutative rings, and that this defect can be remedied by using simplicial commutative rings. This lead to derived algebraic geometry over a ground ring $\KK$ being studied using simplicial commutative $\KK$-algebras as developed by To\"en-Vezzosi and Lurie \cite{toen,Hagy1,Hagy2,dga0,dagl}. In contrast to the general case, it is well-known that when $\KK\supset\QQ$, one gets an equivalent theory using differential graded commutative $\KK$-algebras \cite{rht}, i.e. the defect goes away for geometry over any field containing the rationals. Differential graded commutative algebras are easier to work with than simplicial commutative algebras, as they have a simpler structure, are often directly available, and are amenable to geometric intuition. Since the ground ring of differential geometry is $\RR$, it is natural to expect that an equivalent formalism based on ``dg'' structures should exist there as well, and would lead to a simpler theory. This paper is our first step in showing that this can indeed be accomplished. In the sequel, we will show that calculations, for example of derived intersections of smooth submanifolds, are incredibly easy in our setting, which we believe to be its main advantage over existing approaches. In fact, our approach is also directly applicable to derived \emph{super}geometry.

In \cite{borisovnoel}, it is shown that the simplicial category of derived manifolds it describes is (essentially) equivalent to  that of Spivak in \cite{spivak}, and it is shown in \cite{borisov2} that Joyce's $2$-category of $d$-manifolds \cite{joyce} is (essentially) the $2$-truncation of the simplicial categories of \cite{borisovnoel} and \cite{spivak}. The simplicial approach of Borisov-Noel in \cite{borisovnoel} should be related to ours via a Dold-Kan correspondence. We shall address this question in a future paper.

%When one passes from algebraic to differential geometry, the thing to do is to replace algebras by $\cinf$-algebras. Informally, a \emph{$\cinf$-algebra} can be viewed as a commutative $\RR$-algebra with an extra structure. To describe the extra structure, consider the following. Given a commutative algebra, say $\A$, any polynomial $p=p(x_1,\ldots,x_n)$ in $n$ variables with real coefficients gives rise to an $n$-ary operation $\A^n\to\A$: one substitutes elements $a_1,\ldots,a_n$ of $\A$ for the variables in $p$ and evaluates the resulting expression to obtain a well-defined element $p(a_1,\ldots,a_n)$ of $\A$, interpreting the additions and multiplications occurring in $p$ as those of $\A$. Furthermore, these operations are compatible with the substitution of polynomials into other polynomials in an obvious way. Now, in a $\cinf$-algebra, not only polynomials but all $\cinf$ functions $f:\RR^n\to\RR$ have such an interpretation as operations $\A^n\to\A$ (again, in a way compatible with compositions and identities). A morphism of $\cinf$-algebras is a map of sets $\phi:\A\to\A'$ respecting all these operations.

%A theory of derived manifolds based on simplicial $\cinf$-algebras was recently developed by David Spivak. Since the ground ring of differential geometry is $\RR$, it is natural to expect that an equivalent formalism based on "dg" structures should exist here as well, and would lead to a simpler theory. This paper demonstrates this to be the case.

%etc. etc. etc.

\subsection{Organization and main results}
In section \ref{sec:modules} (and additionally in Appendix \ref{app:modules}) we describe the theory of modules and derivations for algebras over a super Fermat theory $\bS$. Modules for an $\bS$-algebra turn out to be the same as modules for its underlying supercommutative algebra; however, the notion of derivation is different. In sharp contrast to the general case, we show that for a near-point determined $\bS$-algebra $\A,$ the notion of derivation (with values in nice enough modules) is the same as for its underlying supercommutative algebra; however, the module of K\"ahler differentials, which we also define in this section, is still different (for $\A=\cinf(M)$, $M$ a supermanifold, it coincides with the module of smooth $1$-forms on $M$).

Section \ref{sec:monoids} (and also Appendix \ref{sec:coactions}) introduces the general theory of affine algebraic groups and monoids and their actions in the setting of algebras for a super Fermat theory $\bS.$ Of central importance, is the definition of $\hGad^\odd$ $\hGm,$ and $\hMm$- the odd additive group, multiplicative group, and multiplicative monoid in this setting- which we use in Section \ref{sec:dgstruct} to define the category of differential graded $\bS$-algebras. We also describe the theory of Lie algebras and infinitesimal generators for actions.

In Section \ref{sec:dgstruct}, we define various categories of differential graded algebras for a super Fermat theory $\bS,$ including both the bounded and unbounded case. For simplicity, we first treat separately the notion of a differential and that of a grading, and then show how to use the endomorphisms of the odd line to put them together. We consider the category $\bS\Alg^{\E^\times}$ of $\bS$-algebras acted upon by $\EE^\times$. However, we observe that, for theories such as $\bS=\bcinf$, the underlying $\hGm$-action can be quite badly behaved, far from anything resembling an integer grading. Nevertheless, every algebra in $\A\in\bS\Alg^{\E^\times}$ contains a differential graded subalgebra $\A_\alg$, the algebraic part of $\A$, corresponding to the algebraic characters of $\hGm$, and we restrict attention to those $\A$ which are completely determined by $\A_\alg$, in the sense of $\bS$-completion. We call such $\A$ \emph{essentially algebraic}. It turns out that their full subcategory is equivalent to the category of differential graded supercommutative superalgebras equipped with an extra $\bS$-algebra structure on the subalgebra of $0$-cocycles, and these are what we define to be differential graded $\bS$-algebras. We obtain an ``algebraization-completion'' adjunction
\begin{equation}\label{eqn:algcomp}
\Adj{\tau_{\E^\times}^\dagger}{\dg\bS\Alg}{\bS\Alg^{\E^\times}}{\tau^{\E^\times}_\dagger}.
\end{equation}
We also show how to naturally define cohomology as a functor from $\bS\Alg^{\E^\times}$ to $\bS\Alg^{\hmG}$ ($\bS$-algebras acted upon by $\hGm$), and from differential graded $\bS$-algebras to graded $\bS$-algebras (graded superalgebras with an additional $\bS$-algebra structure on the degree-zero subalgebra) in a way that commutes with taking the algebraic part.

Section \ref{sec:diffforms} concerns itself with defining for each algebra $\A$ for a super Fermat theory $\bS,$ a differential graded $\bS$-algebra $\wom\left(\A\right)$ of differential forms, and defining what it means for a super Fermat theory $\bS$ to admit integration. The latter is of central importance, as it implies the homotopy invariance axiom, which is needed to establish a model structure on differential graded $\bS$-algebras.

In Section \ref{sec:models}, we finally show that when a super Fermat theory $\bS$ admits integration, one can construct a cofibrantly generated Quillen model structure on various versions of differential graded $\bS$-algebras, by transfer from the classical (projective) model structure on cochain complexes. In particular, we prove the following:

\begin{thm}
There is a cofibrantly generated Quillen model structure on the category $\dg\bS\Alg$ of differential graded $\bS$-algebras, with surjective maps as the fibrations, and cohomology isomorphisms as the weak equivalences. There is also a model structure on $\bS\Alg^{\E^\times}$ making \eqref{eqn:algcomp} a Quillen adjunction.
\end{thm}

From this one can easily induce model structures on bounded (positively or negatively graded) dg $\bS$-algebras via the ``inclusion-truncation'' adjunction (or obtain them directly). Alternatively, following \cite{Hagy2}, one can consider a kind of ``$t$-structure'' on the category $\dg\bS\Alg$ (or $\bS\Alg^{\E^\times}$), consisting of two subcategories (of algebras with cohomologies concentrated in non-positive, respectively non-negative degrees) whose intersection is equivalent to the category $\bS\Alg$ of discrete objects. We shall not pursue this approach here.

We conclude the paper with Section \ref{sec:simplicial}, where we define the notion of a right (and dually left) almost simplicial model category, which is a simplicial enrichment, compatible in a suitable sense with a given model structure, but is a weaker condition than being a simplicial model category in the sense of \cite{ha}; however, it does satisfy the important property that for cofibrant $X$ and fibrant $Y$, the simplicial set $$\IHom\left(X,Y\right)$$ is a Kan complex. (Right almost simplicial model structures on categories of algebras over linear operads over ground rings containing $\QQ$ were constructed by Hinich in \cite{haha}; we have merely axiomatized such simplicial structures.) We then go on to show that each of the model structures established in Section \ref{sec:models} are right almost simplicial, by using differential forms on simplices to induce a simplicial enrichment, in the spirit of \cite{pldr} and \cite{infitop}.

\vspace{0.2in}\noindent{\bf Acknowledgment:}
We are grateful to Christian Blohmann, Dennis Borisov, Vlad-imir Hinich, Dominic Joyce, Anders Kock, Ieke Moerdijk, Justin Noel, David Spivak, and Peter Teichner for useful conversations. The first author would like to additionally thank the many participants in the ``Higher Differential Geometry'' seminar (formerly known as the ``Derived Differential Geometry'' seminar) at the Max Planck Institute for Mathematics. The second author was supported by the Dutch Science Foundation ``Free Competition'' grant. He would also like to thank the Radboud University of Nijmegen, where part of this work was carried out, for hospitality.

\newpage

\subsection{Notation and conventions} We shall freely use the notation and results from \cite{dg1}. The absolute ground ring for the length of this paper will be $\QQ$; thus, $\bcom$ will mean $\bcom_\QQ$ and $\bscom$ -- $\bscom_\QQ$, and the word ``ring'' will mean ``$\QQ$-algebra''.

For a super Fermat theory $\bS$ with ground ring $\KK,$ denote by
\[
\tau_\KK:\bscom_\KK\To\bS,
\]
the canonical map of theories, to distinguish it from the structure map
\[
\tau_\bS:\bscom\To\bS.
\]
In the induced adjunction
$$\Adj{\tau_\KK^*}{\bscom_\KK\Alg}{\bS\Alg}{\tau^\KK_!}$$
the right adjoint $\tau_\KK^*$ will be referred to as the underlying $\KK$-algebra functor and will be denoted by $(\quad)_\sharp$, while the left adjoint $\tau^\KK_!$ will be denoted by $\widehat{(\quad)}$ and referred to as the $\bS$-algebra completion functor when there is no confusion about the ground ring (given a $\B\in\bS\Alg$, the completion of a $\B_\sharp$-algebra as a $\B$-algebra is different from its completion as a $\KK$-algebra).

\section{Modules and derivations.}\label{sec:modules}
In this section, we present the theory of modules in the context of a super Fermat theory $\bS.$ Many results are a natural generalization of those of \cite{1forms}, in the more general setting of a super Fermat theory. In order to introduce the correct categorical notions, we recall the general sense of module in an arbitrary category. The ideas trace back to the thesis of Jon Beck \cite{beck}.

\begin{defn}\label{dfn:beckmodule}
For $X$ an object of a category $\bC,$ a \textbf{module} over $X$ is an abelian group object $M$ in the slice category $\bC/X.$ The category of such abelian group objects, $\Ab\left(\bC/X\right)$ shall be denoted by $\Mod\left(X\right).$
\end{defn}

\begin{rem}
This concept of module is sometimes called a \emph{Beck} module.% In the case where $\bC$ is the category $\comalg$ of commutative algebras, this definition agrees with the classical notion of a module over and algebra.
\end{rem}

\begin{eg}
If $\A$ is a supercommutative $\KK$-algebra, with $\KK$ a supercommutative ring, the category $$\Mod\left(\A\right)=\Ab\left(\bscom_\KK/\A\right)$$ is canonically equivalent to the category of super $\A$-modules in the usual sense. The corresponding result in the case of commutative algebras is standard. For the reader's convenience, we provide a proof of the super case in the appendix (Proposition \ref{prop:modsame}.)
\end{eg}

\begin{prop}
Let $F:\bC \to \bD$ be a pullback-preserving functor. Then, for all $X$ in $\bC,$ there is an induced additive functor $$F_X:\Mod\left(X\right) \to \Mod\left(FX\right).$$
\end{prop}
\begin{proof}
The condition that $F$ preserves pullbacks implies that for all $X,$ the induced functor $$F_X:\bC/X \to \bD/X$$ preserves finite products. Hence, $F_X$ preserves abelian group objects and restricts to an additive functor $$F_X:\Mod\left(X\right) \to \Mod\left(FX\right).$$
\end{proof}

For our purposes, the ambient category $\bC$ shall always be assumed to be the category of algebras for an algebraic theory. Such a category $\bC$ is always Barr exact. It follows that each category of modules $$\Mod\left(X\right)=\Ab\left(\bC/X\right)$$ is abelian.

\begin{cor}
Suppose that $f:\bT \to \bT'$ is a morphism of algebraic theories. Then, for every $\bT$-algebra $\A,$ it induces an additive functor $$f^*_\A:\Mod\left(\A\right) \to \Mod\left(f^*\A\right).$$
\end{cor}

\begin{thm}\label{thm:moduelssame}
Let $\bS$ be a super Fermat theory, and let $$\tau_\KK:\bscom_{\KK} \to \bS$$ be its structure map, where $\KK=\bS\left(0|0\right)$ is the ground ring. Then for every $\bS$-algebra $\A,$ the induced functor $$\left(\tau^*_\KK\right)_\A:\Mod\left(\A\right) \to \Mod\left(\A_\sharp\right),$$ is an equivalence of categories.
\end{thm}
\begin{proof}
Let $\pi:\B \to \A_\sharp$ be an abelian group object in $\bscom_{\KK}/\A_{\sharp}$ with unit map $$z:\A_\sharp \to \B.$$ Then, $\pi$ is a square zero extension of $\A,$ so in particular a split nilpotent extension with $z$ as a splitting (See Appendix \ref{sec:sq0}). By \cite{dg1} Proposition \ref{prop:nilpextsup}, it follows that $\B$ has the canonical structure of an $\bS$-algebra such that $\pi$ becomes a morphism of $\bS$-algebras. This is an abelian group object in $\bS\Alg/\A$. If $$\pi':\B' \to \A_\sharp$$ is another abelian group object in $\bscom_{\KK}/\A_{\sharp}$, and $$\Phi:\B' \to \B$$ is a morphism of group objects (in particular this implies that $\Phi$ is a morphism over $\A_\sharp$ which respects the unit maps), then by \cite{dg1} Proposition \ref{prop:nilpextsup}, $\Phi$ is in fact a morphism of $\bS$-algebras. This defines a full and faithful functor $$\left(\tau^!_{\KK}\right)_\A:\Mod\left(\A_\sharp\right) \to \Mod\left(\A\right).$$ Moreover, it is clear that $$\left(\tau^*_\KK\right)_\A\circ\left(\tau^!_{\KK}\right)_\A=id.$$ Now suppose that $$\varphi:\C \to \A$$ is an abelian group object in $\bS\Alg/\A.$ Then $\varphi_\sharp:\C_\sharp \to \A_\sharp$ exhibits $\C$ as an split nilpotent extension of $\A$. It follows that $$\left(\tau^!_{\KK}\right)_\A\circ\left(\tau^*_\KK\right)_\A=id.$$
\end{proof}

\begin{rem}
A special case of this theorem, in the case of the Fermat theory $\bcinf,$ is proven in the Master's thesis of Herman Stel \cite{stel}.
\end{rem}

We conclude that for an $\bS$-algebra $\A,$ modules over $\A$ are the same as super modules over its underlying $\KK$-algebra $\A_\sharp$.

\subsection{Derivations}

\begin{defn}\label{dfn:derivation}
Let $\A$ be an $\bS$-algebra and $M$ an $\A$-module. A \emph{derivation of $\A$ with values in $M$} is a section $\sigma$ of $$\pi_M:M\e \to \A,$$ the square-zero extension associated to $M.$ We will also refer to this as an \emph{even derivation} of $\A$ with values in $M$; an \emph{odd derivation} of $\A$ with values in $M$ will mean a derivation of $\A$ with values in $\Pi M.$ Moreover, an even derivation of $\A$ will mean a derivation in $\A$ with values in $\A,$ where $\A$ is viewed as a module over itself in the canonical way; an odd derivation of $\A$ will mean a derivation with values in $\Pi\A.$
\end{defn}

On the surface, it may seem that this definition agrees with the notion of derivation of $\A_\sharp$ with values in $M$ in the usual sense; however, in general not every $\bscom_\KK$-algebra map $\sigma$ splitting $\pi$ will be an $\bS$-algebra map. So, even though the categories of modules $\Mod\left(\A\right)$ and $\Mod\left(A_\sharp\right)$ are naturally equivalent, the theory of derivations of $\A$ is different from that of $\A_\sharp$. We will unwind Definition \ref{dfn:derivation} to highlight the difference, but first we will start with an example.

Recall the notion of a jet algebra from \cite{dg1} (Definition \ref{dfn:jet}). Notice that, viewing $\A$ as a module over itself, one has
\begin{equation}\label{eq:aaejeo}
\A\e \cong \A \odot \J_{1|0}^1=\A \otimes \J_{1|0}^1
\end{equation}
(the latter equality holding since $\J_{1|0}^1$ is a Weil algebra), with $$\pi_\A:\A\e \to \A$$ corresponding to the map induced by pairing the identity of $\A,$ with the map $$\J_{1|0}^1=\KK\{t\}/(t^2) \to \KK \to \A,$$ where $$\KK\{t\}/\left(t^2\right) \to \KK$$ is the canonical map defining $\J_{1|0}^1$ as a Weil algebra. Moreover, notice that $$\A\e \cong \A\odot\KK\{t\}/\left(t^2\right)=\A\set{t}/\left(t^2\right),$$ and under this identification (with $\epsilon=t$ modulo $t^2$), the map $\pi_\A$ is just is determined by setting $t=0,$ and hence we denote it as $\ev_0,$ ``evaluation at $t=0.$''
Similarly, one has that $$\A\ep \cong \A \odot \J_{0|1}^1=\A \otimes \J_{0|1}^1.$$ Notice that $\J_{0|1}^1$ is nothing but $\Lambda^1$- the free $\bS$-algebra on one odd generator $\eps$, and hence we have that $$\A\ep\cong \A\set{\eps},$$ -the free $\A$-algebra on one odd generator. We also have a map (which we denote by the same name) $$ev_0:\A \odot \J_{0|1}^1=\A\set{\eps} \to \A.$$
Let $f\in \KK\set{x_1,\ldots,x_m;\xi_1,\ldots,\xi_n}.$ By the super Fermat property (\cite{dg1}, Definition \ref{dfn:superfermat}), for all $i$, there exists unique $$g_i \in \KK\set{z,w,x_1,\ldots,\widehat{x_i},\ldots,x_m;\xi_1,\ldots,\xi_n}$$ and $$h_i \in \KK\set{x_1,\ldots,x_m;\theta,\eta,\xi_1,\ldots\widehat{\xi_i},\ldots,\xi_n}$$ (where hat represents omission), such that
\begin{equation}
 \resizebox{4.5in}{!}{$f\left(x_1,\ldots,z,\ldots,x_m, \mathbf{\xi}\right)-f\left(x_1,\ldots,w,\ldots,x_m,\mathbf{\xi}\right)=\left(z-w\right)\cdot g_i\left(z,w,x_1,\ldots,\widehat{x_i},\ldots,x_n,\mathbf{\xi}\right)$}
\end{equation}
and
\begin{equation}
 \resizebox{4.5in}{!}{$f\left(\mathbf{x},\xi_1,\ldots,\theta,\ldots,\xi_n\right)-f\left(\mathbf{x},\xi_1,\ldots,\eta,\ldots,\xi_n\right)=\left(\theta-\eta\right)\cdot h_i\left(\mathbf{x},\theta,\eta,\xi_1,\ldots,\widehat{\xi_i},\ldots,\xi_n\right)$.}
\end{equation}
\begin{defn}
With $f$ as above, we define $$\frac{\partial f}{\partial x_i} \in \KK\set{x_1,\ldots,x_m;\xi_1,\ldots,\xi_n}$$ by $$\frac{\partial f}{\partial x_i}\left(\mathbf{x},\mathbf{\xi}\right)=g_i\left(x_i,x_i,x_1,\ldots,\widehat{x_i},\ldots,x_n,\mathbf{\xi}\right),$$ and $$\frac{\partial f}{\partial \xi_i} \in \KK\set{x_1,\ldots,x_m;\xi_1,\ldots,\xi_n}$$ by $$\frac{\partial f}{\partial \xi_i}\left(\mathbf{x},\mathbf{\xi}\right)=h_i\left(\mathbf{x},\xi_i,\xi_i,\xi_1,\ldots,\widehat{\xi_i},\ldots,\xi_n\right).$$
Let $\A=\KK\set{x_1,\ldots,x_m;\xi_1,\ldots,\xi_n}.$ Then we have maps of $\KK$-modules
\begin{eqnarray*}
\frac{\partial}{\partial x_i}:\A &\to& \A\\
f &\mapsto& \frac{\partial f}{\partial x_i}
\end{eqnarray*}
and
\begin{eqnarray*}
\frac{\partial}{\partial \xi_i}:\A &\to& \Pi\A\\
f &\mapsto& \frac{\partial f}{\partial \xi_i}.
\end{eqnarray*}
We define $\frac{\partial f}{\partial x_i}$ and $\frac{\partial f}{\partial \xi_i}$ to be the \emph{(left) partial derivatives} of $f$ with respect to $x_i$ and $\xi_i$ respectively.
Notice that, as $\A$-modules, one has $$\A\e=\A \oplus \A$$ and $$\A\ep=\A \oplus \Pi\A.$$
Define $\KK$-module maps:
\begin{eqnarray*}
\partial_{x_i}:\A &\to& \A\e\\
f &\mapsto& f+\epsilon\frac{\partial f}{\partial x_i}
\end{eqnarray*}
and
\begin{eqnarray*}
\partial_{\xi_i}:\A &\to& \A\ep\\
f &\mapsto& f+\eps\frac{\partial f}{\partial \xi_i}.
\end{eqnarray*}
It is easy to see that these are both in fact maps of supercommutative $\KK$-algebras, and moreover, are sections of the canonical projections down to $\A.$ By \cite{dg1}, Corollary \ref{cor:npdralgs}, it follows that they are in fact both maps of $\bS$-algebras (since $\A,$ $\A\e,$ and $\A\ep$ are all near-point determined). It follows that for all $i,$ $\partial_{x_i}$ is an even derivation of $\KK\set{x_1,\ldots,x_m;\xi_1,\ldots,\xi_n}$ and $\partial_{\xi_i}$ is an odd derivation.
\end{defn}

\begin{rem}\label{rmk:sqzfm}
Let $M$ be an $\A$-module. Then $$\pi_M:M\e \to \A$$ is a square-zero extension, and inherits an $\bS$-algebra structure. Using partial derivatives, one can describe this structure explicitly as follows:
Suppose $$f \in \KK\set{x_1,\ldots,x_m;\xi_1,\ldots,\xi_n}.$$ An arbitrary element of $M\e$ has the form $a + \epsilon u,$ with $a \in \A$ and $u \in M.$ Let $$\left(\mathbf{a},\mathbf{b}\right)=\left(\left(a^1,\ldots,a^m\right),\left(b^1,\ldots,b^n\right)\right) \in \A^m_\even \times \A^n_\odd,$$ and $$\left(\mathbf{u},\mathbf{v}\right)=\left(\left(u^1,\ldots,u^m\right),\left(v^1,\ldots,v^n\right)\right) \in M^m_\even \times M^n_\odd,$$ The $\bS$-algebra structure is determined by the formula
\begin{equation}
\resizebox{4.5in}{!}{$\M\e\left(f\right)\left(\mathbf{a}+\epsilon \mathbf{u},\mathbf{b}+\epsilon\mathbf{v}\right)=\A\left(f\right)\left(\mathbf{a},\mathbf{b}\right) + \epsilon \cdot\left(\sum \limits_{i=1}^{m} u^i \A\left(\frac{\partial f}{\partial x_i}\right)\left(\mathbf{a},\mathbf{b}\right)+\sum \limits_{j=1}^{n} v^i \A\left(\frac{\partial f}{\partial \xi_j}\right)\left(\mathbf{a},\mathbf{b}\right)\right).$}
\end{equation}
\end{rem}

Recall that the usual definition of a derivation of a supercommutative ring $\R$ with values in an super $\R$-module $M$ is an $\R$-module map $$D:\R \to M$$ such that
\begin{equation}\label{eq:usder}
D\left(a \cdot b\right)=a \cdot D\left(b\right) + D\left(a\right)\cdot b.
\end{equation}
Let us unwind Definition \ref{dfn:derivation} to express it in a similar way:

Suppose that $$\pi_M:M\e \to \A$$ is a square-zero extension with a section $\sigma.$ Then since, as a $\KK$-module, $$M\e \cong \A \oplus M,$$ we have that $\sigma$ is of the form $$\sigma\left(a\right)=a + \epsilon D\left(a\right)$$ for some $\KK$-module map $$D:\A \to M.$$ Let us now derive what properties the map $D$ must have. Suppose $$f \in \KK\set{x_1,\ldots,x_m;\xi_1,\ldots,\xi_n},$$ and $$\left(\mathbf{a},\mathbf{b}\right) \in \A^m_\even \times \A^n_\odd.$$ Suppose that $D:\A \to M$ is an arbitrary map of $\KK$-modules, and consider the module map
\begin{eqnarray*}
T:\A &\to& M\e\\
a &\mapsto& a+\epsilon D\left(a\right).
\end{eqnarray*}
In order for $T$ to be an $\bS$-algebra map, we would need the following equation to hold:
$$T\left(\A\left(f\right)\left(\mathbf{a},\mathbf{b}\right)\right)=M\e\left(f\right)\left(T\left(\mathbf{a},\mathbf{b}\right)\right).$$
The left-hand side is equal to $$\A\left(f\right)\left(\mathbf{a},\mathbf{b}\right)+\epsilon \cdot D\left(\A\left(f\right)\left(\mathbf{a},\mathbf{b}\right)\right),$$ whereas the right-hand side, by virtue of Remark \ref{rmk:sqzfm} is  $$\A\left(f\right)\left(\mathbf{a},\mathbf{b}\right) + \epsilon \cdot\left(\sum \limits_{i=1}^{m} D\left(a^i\right) \A\left(\frac{\partial f}{\partial x_i}\right)\left(\mathbf{a},\mathbf{b}\right)+\sum \limits_{j=1}^{n} D\left(b^j\right) \A\left(\frac{\partial f}{\partial \xi_j}\right)\left(\mathbf{a},\mathbf{b}\right)\right).$$ We conclude the following:

\begin{prop}\label{prop:deriv}
There is a natural bijection between $\bS$-algebra derivations of $A$ with values in $M$ and $\KK$-module maps $$D:A \to M$$ such that for all $$f \in \KK\set{x_1,\ldots,x_m;\xi_1,\ldots,\xi_n},$$ and $$\left(\mathbf{a},\mathbf{b}\right) \in \A^m_\even \times \A^n_\odd,$$
\begin{equation}\label{eq:derivationform}
D\left(\A\left(f\right)\left(\mathbf{a},\mathbf{b}\right)\right)=\sum \limits_{i=1}^{m} D\left(a^i\right) \A\left(\frac{\partial f}{\partial x_i}\right)\left(\mathbf{a},\mathbf{b}\right)+\sum \limits_{j=1}^{n} D\left(b^j\right) \A\left(\frac{\partial f}{\partial \xi_j}\right)\left(\mathbf{a},\mathbf{b}\right).
\end{equation}
\end{prop}

\begin{rem}
If $\bS=\bscom_\KK,$ it is easy to check that (\ref{eq:derivationform}) holds for all polynomials if and only if (\ref{eq:usder}) holds, hence one recovers the usual definition of derivation.
\end{rem}

\begin{rem}
There are ``sign rules'' hidden in (\ref{eq:derivationform}), although the formula itself contains no signs. For example, suppose that $D$ is an odd derivation, and that $\alpha^1$ and $\alpha^2$ are odd elements of $\A.$ Consider the multiplication function
$$f\left(\xi^1,\xi^2\right)=\xi^1\xi^2 \in \bS\left(0|2\right)_0.$$ The partial derivatives are
$$\frac{\partial f}{\partial \xi^1}=\xi^2$$ and $$\frac{\partial f}{\partial \xi^2}=-\xi^1$$ (the sign occurs since we are taking \emph{left} partial derivatives.) Note that
$$A\left(f\right)\left(\alpha^1,\alpha^2\right)=\alpha^1 \cdot \alpha^2.$$
Applying (\ref{eq:derivationform}), one sees that
$$D\left(\alpha^1\cdot\alpha^2\right)=D\left(\alpha^1\right)\cdot \alpha^2-\alpha^1 \cdot D\left(\alpha^2\right),$$
thus recovering the usual Koszul sign in the Leibniz rule.
\end{rem}

\begin{defn}
If $\A$ is a $\bS$-algebra and $M$ an $\A$-module, define the $\A$-module $\Der\left(\A,M\right),$ to be the super $\KK$-module such that $\Der\left(\A,M\right)_\even$ is the $\KK$-module of even derivations of $\A$ with values in $M,$ and $\Der\left(\A,M\right)_\odd$ is the $\KK$-module of odd derivations of $\A$ with values in $M,$ with the obvious $\A$-module structure.
\end{defn}

\begin{rem}
From Proposition \ref{prop:deriv}, it follows that the underlying super $\KK$-module of $\Der\left(\A,M\right)$ is a submodule of $\IHom_{\sv}\left(\A,M\right).$
\end{rem}

\begin{defn}
Given an $\bS$-algebra $\A,$ equip $\Der\left(\A,\A\right)$ with a bracket by defining
$$\left[\X,\Y\right]\left(f\right)=\X\left(f\right)\Y\left(f\right)-\left(-1\right)^{\deg\left(\X\right)\deg\left(\Y\right)} \Y\left(f\right) \X\left(f\right),$$ for all elements $f$ of $\A.$ This gives $\Der\left(\A,\A\right)$ the canonical structure of a Lie superalgebra, called the Lie superalgebra of derivations of $\A.$
\end{defn}

\begin{prop}\label{prop:npdfreeder}
Let $\A$ be a near-point determined $\bS$-algebra (\cite{dg1}, Definition \ref{dfn:npd}) and $M$ a free super $\A$-module. Then every $\KK$-algebra derivation of $\A_\sharp$ with values in $M$ is an $\bS$-algebra derivation of $\A$ with values in $M.$
\end{prop}

\begin{proof}
The isomorphism (\ref{eq:aaejeo}) easily generalizes to
\begin{equation}\label{eq:annn}
\A^{m|n}\e\cong \A \odot \J_{n|m}^1,
\end{equation}
where $\A^{n|m}$ is the free $\A$-module on $m$ even and $n$ odd generators, i.e,
$$\A^{n|m}=\A^m \oplus \left(\Pi \A\right)^{m}.$$
Notice that $\J_{n|m}^1$ is a Weil algebra. More generally, if $S$ and $T$ are (not necessarily finite) sets, one can analogously define the $1^{st}$ jet algebra on the set $S$ of even generators, and the set $T$ of odd generators, $\J_{S|T}^1$ which in general is a formal Weil algebra, and an analogous equation to (\ref{eq:annn}) holds. By \cite{dg1}, Proposition \ref{prop:tensweilnpd}, it follows that if $M$ is a free super $\A$-module, then $M\e$ is near-point determined. By \cite{dg1}, Corollary \ref{cor:npdralgs}, the result now follows.
\end{proof}

\begin{cor}
If $\A$ is a near-point determined $\bS$-algebra, and $N$ a submodule of a free super $\A$-module, then every $\KK$-algebra derivation of $\A_\sharp$ with values in $N$ is an $\bS$-algebra derivation of $\A$ with values in $N.$ In particular, this holds for projective $\A$-modules.
\end{cor}

\begin{proof}
By Theorem \ref{thm:moduelssame} and Proposition \ref{prop:modsame}, it follows that $N\e$ is a subalgebra of $M\e$ (or this can be checked directly). By \cite{dg1}, Remark \ref{rem:subpt}, it follows that $N \e$ is near-point determined. The rest is the same as the proof of Proposition \ref{prop:npdfreeder}.
\end{proof}

\begin{eg}
Let $M$ be a (paracompact Hausdorff $2^{nd}$ countable) supermanifold. Then even/odd $\bS\bcinf$-derivations of $\bcinf\left(M\right)$ are in agreement with derivations of the underlying supercommutative algebra, and hence correspond to even/odd vector fields on $M.$
\end{eg}

\begin{rem}
We warn the reader that even for a near-point determined $\bS$-algebra $\A$, the module of K\"ahler differentials on $\A$ as an $\bS$-algebra is different (a quotient of) the algebraic K\"ahler differentials of $\A_\sharp,$ as we will see in the next subsection.
\end{rem}

\subsection{K\"ahler differentials}

Let $\A$ be an $\bS$-algebra. Consider the codiagonal $$\nabla:\A \odot \A \to \A.$$ Denote its kernel by $I$. It comes equipped with a super $\KK$-module map
\begin{eqnarray*}
\bar d:\A &\to& I\\
a &\mapsto& a\otimes 1 -1 \otimes a,
\end{eqnarray*}
where we have abused notation by using the canonical embedding of supercommutative $\KK$-algebras $$\A \otimes \A \to \A \odot \A.$$ By composition with the projection, there is a $\KK$-module map $$d:\A \to I/I^2.$$

\begin{lem}
The map $d$ is an $\bS$-algebra derivation of $\A,$ with values in the $\A$-module $I/I^2.$
\end{lem}

\begin{proof}
Denote by
\begin{eqnarray*}
i_1:\A &\to& \A \odot \A\\
a &\mapsto& a \otimes 1,
\end{eqnarray*}
and similarly for $i_2.$ Both $i_1$ and $i_2$ are $\bS$-algebra maps. It follows that for $$f \in \KK\set{x_1,\ldots,x_m;\xi_1,\ldots,\xi_n},$$
$$1\otimes\A\left(f\left(a_1\ldots,a_m;b_1,\ldots,b_n\right)\right)=\left(\A\odot\A\right)\left(f\right)\left(1 \otimes a_1,\ldots,1\otimes a_m;1\otimes b_1,\ldots,1 \otimes b_n\right)$$ and
$$\A\left(f\left(a_1\ldots,a_m;b_1,\ldots,b_n\right)\right)\otimes 1=\left(\A\odot\A\right)\left(f\right)\left(a_1\otimes 1,\ldots,a_n \otimes 1;b_1 \otimes 1,\ldots,b_n \otimes 1\right).$$
From the first order Taylor expansion \cite{dg1}, Corollary \ref{cor:taylor}, it follows that the difference between these two terms (which is $\bar d\left(\A\left(f\left(a_1\ldots,a_m;b_1,\ldots,b_n\right)\right)\right)$) is equal to
$$\resizebox{5in}{!}{$\sum\limits_{i=1}^{m}
\left(a_i \otimes 1 - 1\otimes a_i\right)\A\left(\frac{\partial f}{\partial
x_i}\right)\left(a_1,\ldots,a_m;b_1,\ldots,
b_n\right)+ \sum\limits_{j=1}^{n}
\left(b_i \otimes 1 - 1\otimes b_i\right)\A\left(\frac{\partial f}{\partial
\xi_j}\right)\left(a_1,\ldots,a_m;b_1,\ldots,b_n\right)$}$$ modulo the ideal $I^2.$ This is precisely the condition that $d$ is a derivation, by Proposition \ref{prop:deriv}.
\end{proof}

\begin{defn}\label{dfn:kahlerforms}
The super $\A$-module $I/I^2,$ is the \emph{module of $\bS$-K\"ahler differentals} on $\A,$ and will be denoted by $\Omega^1\left(\A\right).$
\end{defn}

\begin{thm}
If $X:\A \to M$ is a $\bS$-derivation of $\A$ with values in $M,$ then there is a unique map of $\A$-modules $$f_X:\Omega^1\left(\A\right)  \to M,$$ such that $X=f_X \circ d.$
\end{thm}

\begin{proof}
Denote by $$\phi_X:\A \odot \A \to M\e$$ the map induced by the zero section paired with $X.$ Note that the following diagram commutes
$$\xymatrix{\A \odot \A \ar[r]^{\phi_X} \ar[rd]_-{\nabla} & \M\e \ar[d]^-{\pi_M}\\
& \A,}$$
where $\nabla$ denotes the codiagonal. It follows that $$\phi_X\left(I\right) \subset M=\Ker\left(\pi_M\right).$$ Consider the induced map of super $\A$-modules $$\tilde \phi_X:I \to M.$$ Since $M^2=0,$ it follows that $$I^2 \subset \Ker\left(\tilde \phi_X\right).$$ Hence, there is an induced map of $\A$-modules $$f_X:I/I^2 \to M.$$  By construction, it follows that $$f_X \circ d=X.$$ It suffices to show that $f_X,$ as constructed, is the unique map with this property. The rest of this proof follows \cite{1forms} almost verbatim, however we include it here for convenience. Suppose that $g_X:I/I^2 \to M$ is an $\A$-module map such that $$g_X \circ d=X.$$ We wish to show that $f_X=g_X.$ Since $$\left(f_X - g_X\right) \circ d= 0,$$ it suffices to show that if $v:I/I^2 \to M$ is an $\A$-module map such that $v \circ d=0,$ then $v=0.$
An arbitrary element of $\A \odot \A,$ is of the form
$$s=\varphi\left(a_1 \otimes 1,\ldots,a_k \otimes 1,1\otimes a_{k+1},\ldots,1 \otimes a_m; b_1 \otimes 1,\ldots,b_l \otimes 1,1 \otimes b_{l+1},\ldots,1\otimes b_n\right),$$ for some $\varphi \in \KK\set{x_1,\ldots,x_m;\xi_1,\ldots,\xi_n}.$
Let $s \in I \subset \A \odot \A.$ Then
\begin{equation}\label{eq:swt}
\nabla\left(s\right)=\varphi\left(a_1,\ldots,a_m,b_1\ldots,b_n\right)=0.
\end{equation}
Notice also that we can rewrite $s$ as
$$\resizebox{5in}{!}{$s=\varphi\left(1\otimes a_1-da_1,\ldots,1\otimes a_k -da_k,1\otimes a_{k+1},\ldots,1\otimes a_m,1\otimes b_1-db_1,\ldots,b_l -db_l,1\otimes b_{l+1},\ldots,1\otimes b_n\right).$}$$
Again by the Taylor formula, it follows that
%\begin{eqnarray*}
%s &\equiv& \varphi\left(1\otimes a_1,\ldots,1\otimes a_m,1\otimes b_1,\ldots,1\otimes b_n\right)-\\
%\sum\limits_{i=1}^{m} A\odot \A\left(\frac{\partial \varphi}{\partial
%x_i}\right)\left(1\otimes a_1,\ldots,1\otimes a_m,1\otimes b_1,\ldots,
%1 \otimes b_n\right)\left(da_i\right)\\
%- \sum\limits_{j=1}^{n}
%A\odot \A\left(\frac{\partial \varphi}{\partial
%\xi_j}\right)\left(1\otimes a_1,\ldots,1\otimes a_m,1\otimes b_1,\ldots,
%1 \otimes b_n\right)\left(db_i\right) \mbox{ mod }I^2
%\end{eqnarray*}
$$\resizebox{4.9in}{!}{$$\xymatrix@C=0.1cm{s\equiv &\varphi\left(1\otimes a_1,\ldots,1\otimes a_m,1\otimes b_1,\ldots,1\otimes b_n\right)-  &\sum\limits_{i=1}^{m} A\odot \A\left(\frac{\partial \varphi}{\partial
x_i}\right)\left(1\otimes a_1,\ldots,1\otimes a_m,1\otimes b_1,\ldots,
1 \otimes b_n\right)\left(da_i\right) &\\
 & &  -\sum\limits_{j=1}^{n}
A\odot \A\left(\frac{\partial \varphi}{\partial
\xi_j}\right)\left(1\otimes a_1,\ldots,1\otimes a_m,1\otimes b_1,\ldots,
1 \otimes b_n\right)\left(db_j\right)& \mbox{mod } I^2}$$}$$
Notice by (\ref{eq:swt}) we have
$$\varphi\left(1\otimes a_1,\ldots,1\otimes a_m,1\otimes,1\otimes b_1,\ldots,1\otimes b_n\right)=1\otimes \varphi\left(a_1,\ldots,a_m,b_1\ldots,b_n\right)=0.$$
Hence, we have that
$$\resizebox{4.9in}{!}{$$\xymatrix@C=0.5cm{s &\equiv & -\sum\limits_{i=1}^{m} A\odot \A\left(\frac{\partial \varphi}{\partial
x_i}\right)\left(1\otimes a_1,\ldots,1\otimes a_m,1\otimes b_1,\ldots,
1 \otimes b_n\right)\left(da_i\right)&\\
&-&\sum\limits_{j=1}^{n}
A\odot \A\left(\frac{\partial \varphi}{\partial
\xi_j}\right)\left(1\otimes a_1,\ldots,1\otimes a_m,1\otimes b_1,\ldots,
1 \otimes b_n\right)\left(db_i\right) &\mbox{mod } I^2\\
& =&-\sum\limits_{i=1}^{m} 1\otimes \A\left(\frac{\partial \varphi}{\partial
x_i}\right)\left(a_1,\ldots,a_m,1b_1,\ldots,
b_n\right)\left(da_i\right)&\\
&&-\sum\limits_{j=1}^{n}
1\otimes \A\left(\frac{\partial \varphi}{\partial
\xi_j}\right)\left(a_1,\ldots,a_m,b_1,\ldots,b_n\right)\left(db_j\right).&}$$}$$
Since $v$ is $\A$-linear, it follows that
$$\resizebox{4.9in}{!}{$v\left(s\right)=-\sum\limits_{i=1}^{m} \A\left(\frac{\partial \varphi}{\partial
x_i}\right)\left(a_1,\ldots,a_m,b_1,\ldots,
b_n\right)\left(vdb_i\right)-\sum\limits_{j=1}^{n}
\A\left(\frac{\partial \varphi}{\partial
\xi_j}\right)\left(a_1,\ldots,a_m,b_1,\ldots,b_n\right)\left(vdb_j\right)=0.$}$$
Since $s$ was arbitrary, this means $v=0.$
\end{proof} 
\section{Affine $\bS$-algebraic monoids and groups, and their actions.}\label{sec:monoids}
%In this section: define the above notions in $\cinf$ context. Remark that both limits and colimits are computed by forgetting the action. Adjunction between algebraic and smooth actions. Infinitesimal generators.\\
\subsection{Comonoids and coactions.}
The notion of an \emph{affine algebraic group} is well known in algebraic geometry (cf. \cite{demgab}). When thinking algebraically, the name affine algebraic group is slightly misleading, as these are not group objects in commutative rings, but rather group objects in the opposite category, i.e. group objects in the category of affine schemes (technically affine schemes of finite type). Similarly for monoids. To describe such objects algebraically, one uses cogroup (resp. comonoid) objects in the category $\bcom_\KK\Alg_\fg$ of finitely generated commutative algebras over some ground ring $\KK$, with coproduct $\otimes_\KK$ and initial object $\KK.$ Given such a comonoid object $\H,$ for any finitely generated algebra $\A,$ the set $\Hom\left(\H,\A\right)$ inherits a canonical monoid structure, that is to say, there is a canonical lift
$$\xymatrix@C=2.5cm{& \mathbf{Mon} \ar[d]\\
\bcom_\KK\Alg_\fg \ar@{-->}[ur] \ar[r]_-{Y_{\bS\Alg_\fg^\op}(\H)}& \Set,}$$
where $Y_{\bS\Alg_\fg^\op}(\H)$ is the functor corepresented by the underlying algebra of $\H.$ Put yet another way, these are monoid objects in the functor category $$\mathbf{Fun}\left(\bcom_\KK\Alg_\fg,\Set\right)$$ which happen to be corepresentable. This notion extends immediately to $\bS$-algebras for any (super) Fermat theory $\bS$.

\begin{defn} Let $\bS$ be a (super) Fermat theory with ground ring $\KK$. A monoid- (resp. group-) valued functor
$H:\bS\Alg_\fg\to\mathbf{Mon}$ (resp. $H:\bS\Alg_\fg\to\mathbf{Grp}$) is called an \emph{affine $\bS$-algebraic monoid} (resp. \emph{group}) if it is isomorphic to $Y_{\bS\Alg_\fg^\op}(\H)$ for a comonoid (resp. cogroup) object $\H$ in the symmetric monoidal category $(\bS\Alg_\fg,\odot,\KK),$ where $Y$ denotes the Yoneda embedding.

For a Fermat theory $\bE$, an \emph{affine $\bE$-algebraic supermonoid} (resp. \emph{supergroup}) is an affine $\bSE$-algebraic monoid (resp. group), where $\bSE$ is the superization of a $\bE$ (cf. \cite{dg1}).
\end{defn}

The category $\bS\Alg\mathbf{Mon}$ of affine $\bS$-algebraic monoids is the full subcategory of $\mathbf{Mon}^{\bS\Alg_\fg}$ spanned by such functors $H$, and is also equivalent to a (non full) subcategory of $\bS\Alg_\fg^\op,$ whose objects are the corresponding comonoids $\H,$ but whose morphisms are only those preserving the comonoid structures. We shall be concerned with actions of a given affine $\bS$-algebraic monoid $H$ on functors of the form $$\bS\Alg(\A,i\circ(\quad)):\bS\Alg_\fg\to\Set$$ (where $i:\bS\Alg_\fg\to\bS\Alg$ is the inclusion) arising from coactions of $\H$ on some $\bS$-algebra $\A.$ The category of algebras in $\bS\Alg$ equipped with an $\H$-coaction will be denoted by $\bS\Alg^\H$ (see Appendix \ref{sec:coactions} for generalities on monoid objects and their coactions).

%\begin{rem}
%Take $\bS=\bcinf$. Thinking of $\bcinf$-algebras informally as ``functions on smooth spaces'', the codiagonal $\nabla:\H\oinfty\H\to\H$ corresponds to restricting to the diagonal. For instance, in case $\H=\cinf(\RR^n)$, we have
%\[
%\nabla:\cinf(\RR^{2n})\To\cinf(\RR^n),\quad h(x,y)\mapsto h(x,x).
%\]
%The restriction of $\nabla$ to $\H\otimes\H\subset\H\oinfty\H$ then corresponds to the usual pointwise multiplication of functions:
%\[
%f\otimes g\mapsto f(x)g(y)\mapsto f(x)g(x)
%\]
%(in fact, the restriction of $\nabla$ to $\H\otimes\H\subset\H\oinfty\H$ is always the multiplication in the underlying algebra $\H_\sharp$).

%It is natural to consider replacing $\nabla$ by a possibly non-commutative operation
%\[
%\star:\H\oinfty\H\To\H
%\]
%to obtain a $\bcinf$ version of a \emph{quantum group}. As far as we know, this approach to non-commutative geometry has not been tried yet; however, this line of inquiry lies beyond the scope of the present paper, so we shall not pursue it any further.
%\end{rem}

\begin{eg}
Let $\bS=\bS\bcinf$. Every Lie supergroup $G$ gives rise to an affine $\bcinf$-algebraic supergroup $\H=\cinf(G)$ with
\[
\Delta:\cinf(G)\oinfty\cinf(G)\simeq\cinf(G\times G)\To\cinf(G)\quad\mathrm{and}\quad\epsilon:\cinf(G)\To\RR
\]
given by the pullback of the corresponding group operations. Here, the codiagonal $\nabla:\H\oinfty\H\to\H$ is the pullback by the diagonal inclusion
\[
\delta:G\To G\times G.
\]
The restriction of $\nabla$ to $\H\otimes\H\subset\H\oinfty\H$ then corresponds to the usual pointwise multiplication of functions:
\[
f\otimes g\mapsto f(x)g(y)\mapsto f(x)g(x)
\]
(in fact, the restriction of $\nabla$ to $\H\otimes\H\subset\H\oinfty\H$ is always the multiplication in the underlying $\RR$-algebra $\H_\sharp$ for any $\H$).

Similarly, a smooth action of $G$ on a supermanifold $M$ is equivalent to a coaction of $\cinf(G)$ on $\cinf(M)$.
\end{eg}

\subsection{Completion and algebraization.}
Recall that the structure map $$\tau=\tau_\bS:\bscom_\KK\to\bS$$ gives rise to an algebraic morphism
\[
\Adj{\tau^*}{\bscom_\KK\Alg}{\bS\Alg}{\tau_!}
\]
where $\tau^*=(\quad)_\sharp$ is the underlying $\KK$-algebra functor and $\tau_!=\widehat{(\quad)}$ is the $\bS$-algebra completion functor. Since the $\bS$-completion preserves coproducts and initial objects, and also takes finitely generated $\bscom_\KK$-algebras to finitely generated $\bS$-algebras, so the completion $\hat{\H}$ of a  comonoid (cogroup) $\H$ in $\bscom_\KK\Alg_\fg$ is a comonoid (cogroup) in $\bS\Alg_\fg$. Hence, if $H=Y_{\bscom_\KK\Alg_\fg^\op}(\H)$ is an affine algebraic monoid (group), $\hat{H}=Y_{\bS\Alg^\op}(\hat\H)=H\circ(\quad)_\sharp$ is an affine $\bS$-algebraic monoid (group).

Specializing the general Proposition \ref{prop:indadj} to the present setting, we obtain
\begin{prop}\label{prop:adjalg}
There is an induced adjunction
\[
\Adj{\tau^*_\H}{\bscom_\KK\Alg^\H}{\bS\Alg^{\hat\H}}{\tau_!^\H},
\]
with
\[
\tau_!^\H(\A,\Phi:\A\to\H\otimes\A)=(\hat\A,\hat\Phi:\hat\A\to\hat\H\odot\hat\A)
\]
and
\[
\tau^*_\H(\B,\Psi:\B\to\hat\H\odot\B)=(\B_\alg,\Psi_\alg:\B_\alg\to\H\otimes\B_\alg),
\]
where
\[
\B_\alg=\B_\sharp\times_{(\hat\H\odot\B)_\sharp}(\H\otimes\B_\sharp)
\]
and $\Psi_\alg$ is the induced coaction, with $\tau_!^\H \dashv \tau^*_\H.$
\end{prop}

\begin{rem}
If the canonical map $\H\otimes\B_\sharp\to(\hat\H\odot\B)_\sharp$ is injective, then $$\B_\mathrm{alg}=\Psi_\sharp^{-1}(\H\otimes\B_\sharp).$$
\end{rem}

Given a coaction $\Phi:\B\to\hat\H\odot\B$, one should think of the $\H$-comodule $\B_\alg$ as the ``algebraic part'' of the $\hat\H$-comodule $\B$. It is natural to isolate the class of coactions which are in some sense freely generated by their algebraic parts. At first glance, it is tempting to take this to mean those coactions in the essential image of $\tau^\H_!$, i.e. those of the form $\hat\Phi:\hat\A\to\hat\H\odot\hat\A$ for some coaction $\Phi:\A\to\H\otimes\A$. However, this would fail to include all trivial coactions, since $\widehat{(\quad)}$ is not essentially surjective. This is not quite what we would like, since we want the $\hat\H$-\emph{coaction} to be freely generated by an $\H$-coaction, rather than the $\bS$-algebra structure to be freely generated by the underlying $\KK$-algebra structure.

To remedy this, we work \emph{relative to the invariant subalgebras}. Recall (Appendix \ref{sec:coactions}) that there is an adjunction
\[
\Adjlong{(\quad)_\H}{\bscom_\KK\Alg}{\bscom_\KK\Alg^\H}{(\quad)_\triv}
\]
where the left adjoint $(\quad)_\triv$ assigns to an $\A\in\bscom_\KK\Alg$ the trivial coaction $(\A,j_2:\A\to\H\otimes\A)$, while the right adjoint $(\quad)_\H$ assigns to a coaction $(\A,\Phi)$ the invariant subalgebra $\A_\H$ which is the equalizer of $\Phi$ and $j_2$. There is an analogous adjunction for $\hat\H$-coactions. Observe that (Remark \ref{rem:triv}), for a coaction $$\Psi:\B\to\hat\H\odot\B,$$ the superalgebras $(\B_\alg)_\H$ and $(\B_{\hat\H})_\sharp$ are naturally isomorphic (in particular, if $\Psi$ is trivial, $\B_\alg$ is isomorphic to $\B_\sharp$).

Consider the category
\[
\bS\Alg^{\hat\H}_\alg=\bscom_\KK\Alg^\H\times_{\bscom_\KK\Alg}\bS\Alg,
\]
the homotopy pullback in the $(2,1)$ category of categories, functors and natural \emph{isomorphisms} of the diagram
\[
\bscom_\KK\Alg^\H\mathrel{\mathop{\longlongrightarrow}^{(\quad)_\H}}\bscom_\KK\Alg\mathrel{\mathop{\longlongleftarrow}^{(\quad)_\sharp}}\bS\Alg.
\]
Its objects are triples $$((\A,\Phi:\A\to\H\otimes\A),\A^0,\phi:\A^0_\sharp\mathrel{\mathop{\longrightarrow}^{\sim}}\A_\H)$$ and morphisms given by the appropriate commutative diagrams (in plain terms: $\H$-comodules with a specified extra $\bS$-algebra structure on the invariant subalgebras, with morphisms $\H$-equivariant maps inducing $\bS$-algebra maps on the invariant subalgebras).

Combining Remark \ref{rem:triv} and Proposition \ref{prop:indadjrel} we arrive at

\begin{prop}\label{prop:adjalg1}
There is an adjunction
\[
\Adj{\tau^\dagger_\H}{\bS\Alg^{\hat\H}_\alg}{\bS\Alg^{\hat\H}}{\tau^\H_\dagger}
\]
such that
\[
\tau^\dagger_\H((\A,\Phi))=((\A_\alg,\Phi_\alg),\A_{\hat\H},\phi_\A:(\A_{\hat\H})_\sharp\to(\A_\alg)_\H)
\]
and
\[
\tau^\H_\dagger((\B,\Psi:\B\to\H\otimes\B),\B^0,\psi)=(\tau_!^{\B^0}(\B),\tau_!^{\B^0}\Psi),
\]
where
\[
\Adj{\tau^*_{\B^0}}{\B^0_\sharp/\bscom_\KK\Alg}{\B^0/\bS\Alg}{\tau_!^{\B^0}}
\]
is the algebraic morphism induced by the morphism of theories $$\tau_{\B^0}:\bscom_{\B^0_\sharp}\to\bS_{\B^0}.$$
\end{prop}

\begin{defn}\label{def:algebraic}
A coaction $\Phi:\A\to\hat\H\odot\A$ is \emph{essentially algebraic} if it is in the essential image of $\tau^\H_\dagger$. The functor $\tau^\H_\dagger$ itself is referred to as the \emph{relative completion}, its right adjoint $\tau^\dagger_\H$ -- as the \emph{algebraization}.
\end{defn}

%\begin{rem}
%If the counit $\widehat{(\A_{\hat\H})_\sharp}\to\A_{\hat\H}$ is an isomorphism, $\Phi$ is essentially algebraic if and only if $\widehat{(\A_\alg)}\to\A$ is an isomorphism.
%\end{rem}

%Oh, and it gives the right answer for $\H=\RR[\lambda]$ the multiplicative monoid and $\A=\cinf(\RR^{\sum_{k\geq0}n_k})$ with $\lambda$ acting %on each $\RR^{n_k}$ by multiplication by $\lambda^k$. To see this, consider an $f\in\A$ such that
%\[
%\Phi(f)(\lambda,x_0,x_1,...)=f(x_0,\lambda x_1,\lambda^2 x_2,\ldots)=\sum_i p_i(\lambda)h_i(x),
%\]
%where the last sum is finite and $p_i$'s are polynomials.

Lastly,

\begin{prop}\label{prop:adjalg2}
There is an adjunction
\[
\Adj{\tau_\H^\circ}{\bscom_\KK\Alg^\H}{\bS\Alg^{\hat\H}_\alg}{\tau^\H_\circ}
\]
(cf. \eqref{eq:adjcirc}) such that the composition
\[
\Adjadj{\tau_\H^\circ}{\bscom_\KK\Alg^\H}{\bS\Alg^{\hat\H}_\alg}{\tau^\H_\circ}{\tau_\H^\dagger}{\bS\Alg^{\hat\H}}{\tau^\H_\dagger}
\]
equals
\[
\Adj{\tau_\H^*}{\bscom_\KK\Alg^\H}{\bS\Alg^{\hat\H}}{\tau^\H_!}.
\]
\end{prop}

\subsection{Linear coactions.}
Let $V$ be a $\KK$-module and $\G$ a comonoid in $\bS\Alg$. Let $\bS(V)$ be the free $\bS$-algebra on (the underlying $\ZZ_2$-graded set of) $V$. A \emph{linear coaction} of $\G$ on $V$ is, by definition, a coaction $$\Phi:\bS(V)\to\G\odot\bS(V)$$ such that $\Phi(V)\subset\G\otimes V$. To put it another way, it is a map of $\KK$-modules
\[
\phi:V\To\G\otimes V
\]
such that its unique extension to a map of $\bS$-algebras
\[
\Phi:\bS(V)\To\G\odot\bS(V)
\]
is a coaction in $\bS\Alg^\G$ (the extension is determined by the fact that $\G\odot\bS(V)$ is the free $\bS_\G$-algebra generated by the $\G$-module $\G\otimes V$).

If $\G$ is finitely generated and $H\in\Mon^{\bS\Alg_\fg}$ is the corresponding affine $\bS$-algebraic monoid, an action of $H$ on a functor represented by a $\KK$-module $V$ is then the same thing as a coaction of $\G$ on the dual module $V^\vee=\Hom(V,\KK)$.

\begin{eg}
If $G$ is a Lie group and $V$ is a finite-dimensional vector space over $\RR$, a smooth linear action of $G$ on $V$ is the same thing as a linear coaction of $\G=\cinf(G)$ on the dual space $V^*$. Indeed, the smooth action map
\[
\rho:G\times V\To V
\]
is equivalent to a morphism of $\bcinf$-algebras
\[
\Phi=\rho^*:\cinf(V)\To\cinf(G\times V)=\cinf(G)\oinfty\cinf(V)
\]
which is a coaction. Furthermore, $\rho$ is linear if and only if $\rho^*$ maps linear functions on $V$ to fiberwise linear functions on $G\times V$, i.e induces a map
\[
V^*\To\cinf(G)\otimes V^*.
\]
\end{eg}

Now suppose $\G=\hat\H$ for some comonoid $\H$ in $\bscom_\KK\Alg$. Then, for every linear $\H$-coaction
\[
\Phi:\mathrm{Sym}_\KK(V)\To\H\otimes_\KK\mathrm{Sym}_\KK(V)
\]
(where $\mathrm{Sym}_\KK(V)$ denotes the free $\bscom_\KK$ algebra on $V$), coming from a map
\[
\phi:V\To\H\otimes V,
\]
the coaction
\[
\hat\Phi:\bS(V)\To\hat\H\odot\bS(V).
\]
is also linear. Conversely, the algebraization of a linear $\hat\H$-coaction is a linear $\H$-coaction.

%Let us now also describe those coactions which are linear. Let $V$ be a super vector space, $\cinf(V)$ the free $\cinf$-superalgebra generated by $V$. The latter contains $V^*$, the linear dual of $V$, as a subspace. We say that a coaction
%\[
%\Phi:\cinf(V)\To\H\oinfty\cinf(V)
%\]
%is \emph{linear} if $\Phi(V^*)\subset\H\otimes V^*\subset\H\oinfty\cinf(V)$. In fact, one can equivalently describe a linear coaction of $\H$ on $V$ as a linear map
%\[
%\phi:V^*\To\H\otimes V^*
%\]
%such that
%\[
%(\Delta\otimes\mathrm{id})\circ\phi=\iota\circ(\mathrm{id}\otimes\phi)\circ\phi,
%\]
%where $\iota:\H\otimes\H\otimes V^*\To\H\oinfty\H\otimes V^*$ is the canonical inclusion. Denote the resulting category of $\H$-comodules by $\mathbf{SVect}_\H$, a subcategory of $\cinfsalg_\H$.\\

%\emph{\textcolor[rgb]{0.00,0.00,1.00}{In fact, comodules can also be described ``a la Quillen'' as abelian group objects in the appropriate \textbf{under}category, in this case $\H/\mathbf{\cinf Bialg}$, where $\mathbf{\cinf Bialg}$ denotes the category of smooth bialgebras, i.e. comonoids in $\cinfsalg$.}}\\

\subsection{The additive and multiplicative groups.}

The following functors are well-known in algebraic geometry: \emph{the even and odd additive group, the multiplicative group, and the multiplicative monoid}
\[
\Gad^\even, \Gad^\odd,\Gm:\bscom_\KK\Alg\To\mathbf{Grp},\quad\Mm:\bscom_\KK\Alg\To\mathbf{Mon}
\]
given, respectively, by
\[
\Gad^\even(\A)=(\A_\even,+),\quad\Gad^\odd(\A)=(\A_\odd,+),\quad\Gm(\A)=(\A_\even^\times,\cdot),\quad\Mm(\A)=(\A_\even,\cdot),
\]
and corepresented, respectively, by $\adG^\even$, $\adG^\odd$, $\mG$ and $\mM$, defined as follows:

\begin{itemize}
\item  $\adG^\even=\KK[x]$ with
\[
\Delta(f)(x,y)=f(x+y),\quad\epsilon(f)=f(0),\quad S(f)(x)=f(-x);
\]
\item  $\adG^\odd=\KK[\xi]=\Lambda_1$ with
\[
\Delta(f)(\xi,\eta)=f(\xi+\eta),\quad\epsilon(f)=f(0),\quad S(f)(\xi)=f(-\xi);
\]
\item  $\mM=\KK[x]$ with
\[
\Delta(f)(x,y)=f(xy),\quad\epsilon(f)=f(1);
\]
\item  $\mG=\KK[x,x^{-1}]$ with $\Delta$ and $\epsilon$ as for $\mM$ and
\[
S(f)(x)=f(x^{-1}).
\]
\end{itemize}
These functors extend to $\bS\Alg$ in an obvious way, by precomposition with $(\quad)_\sharp$. These extensions, denoted by $\hGad^\even$, $\hGad^\odd$, $\hGm$ and $\hMm$, are corepresented by the $\bS$-completions of the corresponding algebras, namely $\hadG^\even=\KK\set{x}$, $\hadG^\odd=\KK\set{\xi}=\Lambda_1$, $\hmG=\KK\set{x,x^{-1}}$ (since $\bS$-completion commutes with localization) and $$\hmM=\KK\set{x},$$ with the structure maps given by the same formulas.

 There is an inclusion of functors $\Gm\hookrightarrow\Mm$ (resp. $\hGm\hookrightarrow\hMm$) corresponding to the localization homomorphism $\KK[x]\to\KK[x,x^{-1}]$ (resp. $\KK\{x\}\to\KK\{x,x^{-1}\}$).

In addition to the structure maps, the multiplicative monoid has an operation of which we shall make frequent use, namely, \emph{evaluation at zero}:
\[
\zeta:\hmM\To\KK,\quad\zeta(f)=f(0).
\]
This operation satisfies a condition dual to that of a zero in a semigroup, namely, the following diagram commutes:
\begin{equation}\label{eq:multmon}
\xymatrix@C=2cm{\KK \ar[d]_{e} & \hmM \ar[l]_-{\zeta} \ar[d]_-{\Delta} \ar[r]^-{\zeta} & \KK \ar[d]^-{e}\\
\hmM & \hmM \odot \hmM \ar[l]^-{\mathrm{id}\odot\zeta} \ar[r]_-{\zeta\odot\mathrm{id}} & \hmM.}
\end{equation}
In particular, $(\zeta\odot\zeta)\circ\Delta=\zeta$ (a zero is an idempotent).

\subsection{Lie superalgebras and infinitesimal generators.}

This subsection follows \cite{demgab} quite closely.

Let $H:\bS\Alg\to\Grp$ be any functor. Define $$T^\even H:\bS\Alg\to\Grp$$ by setting $$T^\even H=H\circ(-\odot\J_{1|0}^1),$$ where $\J_{1|0}^1=\KK\{t\}/(t^2)$. Likewise, define $$T^\odd H=H\circ(-\odot\J_{0|1}^1),$$ where $J_{0|1}^1=\Lambda_\KK^1=\KK\{\tau\}$ ($\tau$ odd). The canonical homomorphisms
\[
\KK\mathrel{\mathop{\longrightarrow}^{i}}\J_{1|0}^1\mathrel{\mathop{\longrightarrow}^{\pi}}\KK
\]
and
\[
\KK\mathrel{\mathop{\longrightarrow}^{i}}\J_{0|1}^1\mathrel{\mathop{\longrightarrow}^{\pi}}\KK
\]
induce morphisms of functors
\[
H\mathrel{\mathop{\longrightarrow}^{i}}T^\even H\mathrel{\mathop{\longrightarrow}^{\pi}}H
\]
and
\[
H\mathrel{\mathop{\longrightarrow}^{i}}T^\odd H\mathrel{\mathop{\longrightarrow}^{\pi}}H
\]
such that $\pi\circ i=\id_H$.

Now let $e:*\to H$ be the group unit. Let $T_e^\even H=\Ker\pi=\pi^{-1}(e),$ namely the pullback
$$\xymatrix{\pi^{-1}\left(e\right) \ar[d] \ar[r] & T^\even H \ar[d]_-{\pi} \\ {*} \ar[r]^{e} & H,}$$
and similarly for $T_e^\odd$. It can be shown (cf. \cite{demgab}) that the functors $T_e^\even$ and $T_e^\odd$ actually take values in $\Ab,$ the category of abelian groups; moreover, the multiplicative monoid actions on $\J^1_{1|0}$ and $\J^1_{0|1}$  induce the same on $T_e^\even$ and $T_e^\odd$. It follows that the functor $T_e H=(T_e^\even H,T_e^\odd H)$ takes values in $\KK$-modules, and $T_e^\odd H=\Pi T_e^\even H$.

Suppose now that $H=Y_{\bS\Alg^\op}(\H)$. Then, since $\J^1_{1|0}$ and $\J^1_{0|1}$ are Weil algebras, hence co-exponentiable, $TH$, and hence $T_eH$, is co-representable. Denote the corresponding $\KK$-module by $\mathfrak{h}=(\mathfrak{h}_\even,\mathfrak{h}_\odd)$. The module $\mathfrak{h}_\even$ consists of all $$v:\H\to\J^1_{1|0}$$ such that $\pi\circ v=\epsilon$, where $\epsilon$ is the counit of the comonoid $\H$ as in Appendix \ref{sec:coactions}, and similarly for $\mathfrak{h}_\odd$.

In exactly the same way as for algebraic groups (cf.\cite{demgab}), one constructs a bracket on $\mathfrak{h}$ and shows that it defines a Lie superalgebra. In fact, if $H=Y_{\bscom_\KK\Alg^\op}(\H)$ and $\hat{H}=H\circ(\quad)_\sharp=Y_{\bS\Alg^\op}(\hat\H)$, then the Lie superalgebras of $H$ and $\hat{H}$ are isomorphic.

Now let $\Phi:\A\to\H\odot\A$ be a coaction, and let $v\in\mathfrak{h}_\even$. Define the infinitesimal generator of the coaction to be
\[
\Phi_v=(\id_\A\odot v)\circ\Phi:\A\To\A \odot \J^1_{1|0}\cong\A\e.
\]
It is easy to check that $\Phi_v$ is a section of the canonical projection, hence corresponds to an even derivation $\phi_v$ of $\A$. Similarly, every $\xi\in\mathfrak{h}_\odd$ gives rise to an odd derivation $\phi_\xi$. The map
\[
\mathfrak{h}\To\Der(\A,\A),\quad(v,\xi)\mapsto(\phi_v,\phi_\xi)
\]
is a homomorphism of Lie superalgebras.

\section{Differential graded structures on $\bS$-algebras.}\label{sec:dgstruct}
%In this section: Graded (or homogeneity) structures as actions of the multiplicative group or monoid. Differentials as actions of the odd line. DG structures as %smooth actions of the endomorphisms or automorphisms of the odd line.\\

%\smallskip

%Our goal here is to give the notions of ``graded'' and ``differential'' meaning in the smooth context.
%\bigskip

\subsection{Differentials.}

Consider an action of the odd additive group $\hGad^\odd$ on a corepresentable functor $Y_{\bS\Alg^\op}(\A)$ corresponding to a coaction
\[
\Phi:\A\To\hadG^\odd\odot\A=\Lambda^1\odot\A\cong\A\oplus\tau\A.
\]
Since this is a coaction, it is a section of the augmentation map $\Lambda^1\odot\A\to\A$, so we can write
\[
\Phi=\id_\A\oplus\tau d,
\]
where $d:\A\to\Pi\A$ is an odd derivation. In this case, the coaction is uniquely determined by its infinitesimal generator. Since the infinitesimal action must be a homomorphism of Lie superalgebras, it follows that
\[
d^2=\frac{1}{2}[d,d]=0.
\]

\begin{defn}
A \emph{differential $\bS$-algebra} is an $\bS$-algebra endowed with an odd $\bS$-derivation $d:\A\to\Pi\A$ (referred to as \emph{the differential}) such that $d^2=0$. A morphism of differential $\bS$-algebras is a morphism of $\bS$-algebras respecting the differential. Denote the corresponding category by $\dd\bS\Alg.$

A \emph{differential $\KK$-superalgebra} is a differential $\bscom_\KK$-algebra; for a Fermat theory $\bE$, a \emph{differential $\bE$-superalgebra} is a differential $\bSE$-algebra.
\end{defn}

The preceding discussion can be summarized by saying that we have the following equivalences of categories:
\[
\bscom_\KK\Alg^{\adG^\odd}\simeq\dd\bscom_\KK\Alg
\]
and
\[
\bS\Alg^{\hadG^\odd}\simeq\dd\bS\Alg.
\]

Applying Propositions \ref{prop:adjalg}, \ref{prop:adjalg1} and \ref{prop:adjalg2}, we get adjunctions
\[
\Adjadj{\tau_{\adG^\odd}^\circ}{\dd\bscom_\KK\Alg}{\bS\Alg^{\hat\adG^\odd}_\alg}{\tau^{\adG^\odd}_\circ}{\tau_{\adG^\odd}^\dagger}{\dd\bS\Alg}{\tau^{\adG^\odd}_\dagger}
\]
composing to
\[
\Adj{\tau_{\adG^\odd}^*}{\dd\bscom_\KK\Alg}{\dd\bS\Alg}{\tau^{\adG^\odd}_!}.
\]
The objects of the category $\bS\Alg^{\hat\adG^\odd}_\alg$ are differential supercommutative $\KK$-algebras $\A$ with an additional $\bS$-algebra structure on the $\Gad^\odd$-invariant subalgebra $\A^\cl$ consisting of the elements annihilated by $d$ (cocycles).

%Forgetting the $\bS$-algebra structure and only remembering the $\KK$-module structure gives rise to the forgetful functor to the category of differential %$\KK$-modules, i.e. $\KK$-modules equipped with an odd $\KK$-linear square-zero operator. When $\KK$ is a ring (as opposed to a superring), these are the same as %$\ZZ_2$-graded chain complexes of $\KK$-modules. As in the general case, this functor has a left adjoint.

\subsection{Grading.}

It is well-known (and easy to see) that an action of $\Gm$ on a superalgebra $\A\in\bscom_\KK\Alg$ amounts to an integer grading on $\A$. Indeed, a coaction
\[
\Phi:\A\To\A\otimes\mG=\A[u,u^{-1}]
\]
gives rise to a decomposition
\[
\A\cong\bigoplus_{n\in\ZZ}\A^n
\]
of $\A$ as a $\KK$-module, where
\[
\A^n=\{f\in\A|\Phi(f)=u^nf\}
\]
is the submodule corresponding to the character $n$ of $\Gm$, such that $1\in\A^0$ and $\A^i\A^j\subset\A^{i+j}$ for each $i,j\in\ZZ.$ (In particular, $\A^0$ coincides with the invariant subalgebra $\A_\mG.$)

If the $\Gm$-action extends to $\Mm$, only the non-negative weights occur in the above decomposition, so that the $\ZZ$-grading is actually an $\NN$-grading.

Notice that the grading coming from the $\Gm$-action is generally independent of the $\ZZ_2$-grading given by the Grassmann parity. A compatibility between the two means that the latter coincides with the modulo-$2$ reduction of the former; in plain terms, components of even (resp. odd) degree are purely even (resp. odd). In particular, it implies that the ground ring $\KK$ is purely even. There is a way of expressing this compatibility purely in terms of the coaction $\Phi$. Introduce the parity involution
\[
P_\Phi=\ev_{-1}\circ\Phi:\A\To\A,
\]
where $\ev_{-1}$ is the evaluation at $x=-1$. It acts on each $\A^n$ by multiplication by $(-1)^n$. We also have the Grassmann parity involution $P_{\mathrm{Grass}}$ which acts as identity on $\A_\even$ and by multiplication by $-1$ on $\A_\odd$. The action grading is compatible with parity if and only if
\[
P_\Phi=P_{\mathrm{Grass}}.
\]
Denote the full subcategory $\bscom_\KK\Alg^\mG$ (resp. $\bscom_\KK\Alg^\mM$) consisting of such coactions by $\bscom_\KK\Alg^\mG_\even$ (resp. $\bscom_\KK\Alg^\mM_\even$). Denote by $\gr\bcom_\KK\Alg$ (resp. $\gr\bscom_\KK\Alg$) the category of $\ZZ$-graded (super) commutative $\KK$-algebras, and denote by $\gr\bcom_\KK\Alg_{\geq0}$ (resp. $\gr\bscom_\KK\Alg_{\geq0}$) the category of $\NN$-graded (super) commutative $\KK$-algebras.

To summarize: we have equivalences of categories
\[
\bscom_\KK\Alg^\mG\simeq\gr\bscom_\KK\Alg,\quad\bscom_\KK\Alg^\mM\simeq\gr\bscom_\KK\Alg_{\geq0}
\]
and
\[
\bscom_\KK\Alg^\mG_\even\simeq\gr\bcom_\KK\Alg,\quad\bscom_\KK\Alg^\mM_\even\simeq\gr\bcom_\KK\Alg_{\geq0}.
\]
In the latter case, the sign rule is determined by the integer weight:
\[
ab=(-1)^{ij}ba
\]
if $a\in\A^i$, $b\in\A^j$, whereas in the former case the weight has no bearing on signs. The invariant subalgebra $\A_\mG=\A^0$ is purely even if the weight is compatible with parity.

\subsection{What happens upon completion?}

The $\bS$-completion of $\mG$ and $\mM$ and the behavior of coactions depends very much on the theory $\bS$ and its ground ring. For instance, for $\bS=\bcomega$, we can make the following easy observation:
\begin{lem}\label{lem:multgroup}
\[
\hGm\cong\hGad\times\mu_2,
\]
where $\mu_2$ is the group of square roots of unity: as a functor, it associates to every algebra $\A$ the (multiplicative) subgroup of elements of order $2$.
\end{lem}
\begin{proof}
This is just saying that we have an isomorphism of Lie groups
\[
(\RR-\{0\},\times)\stackrel{\simeq}{\To}(\RR,+)\times\{\pm 1\}
\]
mapping a $\lambda\in\RR-\{0\}$ to $(\log\abs{\lambda},\frac{\lambda}{\abs{\lambda}})$; the inverse sends $(t,\sigma)$ to $\sigma\exp{t}$.
\end{proof}

\begin{rem}
It follows that the same is true for every extension of $\bcomega$ with ground ring $\RR$, such as $\bcinf$, $\bS\bcomega$ and $\bS\bcinf$.
\end{rem}

In view of this, an action of $\hGm$ consists of a flow and an involution commuting with the flow; in particular, every $\hGad$-action induces a $\hGm$-action, with $\mu_2$ acting trivially. This is in stark contrast with the algebraic situation: the above isomorphism is transcendental, and Lemma \ref{lem:multgroup} above is, in fact, utterly false for $\Gm$ over any ground ring.

Nevertheless, for any $\bS$, any $\hGm$-algebra $\A$ contains a $\Gm$-subalgebra $\A_\mathrm{alg},$ i.e. a $\ZZ$-graded subalgebra: \newline \\
By definition, if
\[
\Phi:\A\To\A\set{u}= \hGm\odot\A
\]
is the coaction map, then $\A_\mathrm{alg}=\Phi^{-1}(\A_\sharp[u])$, hence it splits as a direct sum
\[
\A_\mathrm{alg}=\bigoplus_{n\in\ZZ}\A^n,
\]
where $\A^n$ is the submodule of $\A$ corresponding to the character $n$ of $\Gm$, i.e.
\[
\A^n=\{f\in\A|\Phi(f)=u^nf\}.
\]
In other words, we are restricting attention to the algebraic characters of $\hGm$ only, i.e. those which come from characters of $\Gm$.

In particular, $\A^0=(\A_\mathrm{alg})_\mG=(\A_\hmG)_\sharp$ (cf. Remark \ref{rem:triv}) so $\A^0$ has the extra structure of an $\bS$-algebra. Thus, the category $\bS\Alg^\hmG_\alg$ has as objects $\ZZ$-graded algebras $\A$ with a specified $\bS$-algebra structure on $\A^0$. Morphisms in $\bS\Alg^\hmG_\alg$ preserve the degree as well as the $\bS$-algebra structure on the degree-zero part.

\begin{defn}
A \emph{$\ZZ$-graded $\bS$-algebra} is an object of $\bS\Alg^\hmG_\alg$; a $\ZZ_{\geq0}$- (or $\NN$-) graded $\bS$-algebra is an object of $\bS\Alg^\hmM_\alg$. The category $\bS\Alg^\hmG_\alg$ (resp. $\bS\Alg^\hmM_\alg$) will henceforth be denoted by $\gr\bS\Alg$ (resp. $\gr\bS\Alg_{\geq0}$).
\end{defn}

%We shall denote this category by $\gr\bS\Alg$ from now on and call its objects \emph{$\ZZ$-graded $\bS$-algebras}.

\begin{rem}\label{rem:badact}
For badly behaved $\hGm$-actions, the algebraic part can be very small, even trivial. For example, consider the torus $T^2$ equipped with the Kronecker (irrational rotation) flow. In view of Lemma \ref{lem:multgroup}, this gives a $\hmG$-coaction on $\A=\cinf(T^2)$. For this coaction, $\A^0=\RR$ while $\A^n=0$ for $n>0$. Thus $\widehat{\A_\mathrm{alg}}=\RR$, which is much smaller than $\A$.

For this reason, it would be misleading to call arbitrary objects of $\bS\Alg^\hmG$ ``graded $\bS$-algebras'' (except for $\bS=\bscom_\KK$, where these two notions agree).
\end{rem}

\begin{rem}\label{rem:Mmact}
$\hMm$-actions are better behaved than $\hGm$-actions, even for $\bS=\bcinf$, thanks to having a zero around. Given an $\hmM$-coaction
\[
\Phi:\A\To\A\{u\}=\KK\{u\}\odot\A,
\]
the composition $z=(\zeta\odot\id_\A)\circ\Phi:\A\to\A$ with evaluation at $u=0,$ where $\zeta$ is as in (\ref{eq:multmon}), is an idempotent splitting the inclusion $\A^0\to\A$ (i.e. $\zeta\circ\Phi(\A)\simeq\A^0$). More generally, taking the $k^{th}$ jet of the coaction at $u=0$ -- i.e., composing $\Phi$ with the projection
\[
\A\{u\}\To\A\{u\}/(u^{k+1})
\]
and co-exponentiating -- we get a map
\[
z_k:\A_{\J^k_1}\To\A
\]
and a subalgebra $\A_k=\Im(z_k)$. These subalgebras form an exhausting (and for finitely generated $\A$ -- even stabilizing) filtration
\[
\A^0=\A_0\subset\A_1\subset\A_2\subset\cdots,\quad\A=\bigcup_k\A_k
\]
of $\A$. Moreover, it was observed in \cite{Sev} (see also \cite{Grab}) that, for $\A=\cinf(M)$, $M$ a supermanifold, \emph{every} $\hmM$-coaction is essentially algebraic, so that $M$ has a ``higher vector bundle'' structure corresponding to the above filtration. Furthermore, $$\A^0=\cinf(M_0),$$ where $M_0$ is the supermanifold of fixed points of this action (which is the image of the map induced by the zero element) and $\A$ is the completion (relative to $\A^0$) of an $\NN$-graded algebra which is locally free (over $M_0$).
\end{rem}

\begin{rem}
When the ground ring is $\CC$, $\hGm$-actions are also better behaved, but for a different reason: due to the connectivity of the complex multiplicative group $\CC^\times$, in contrast with its real counterpart. For instance, a complex manifold with a holomorphic $\CC^\times$-action can be covered by $\CC^\times$-invariant homogeneous affine charts.
\end{rem}

%In general, the counit map $\widehat{\A_\mathrm{alg}}\to\A$ is not an isomorphism (for badly behaved $\hGm$-actions, $\widehat{\A_\mathrm{alg}}$ is much smaller %than $\A$), but for the algebras of interest to us -- namely, essentially algebraic coactions with $\bLambda$-point determined invariant subalgebra -- it will %be. In fact, in the Quillen model structure we construct below, this counit map will be a weak equivalence.

%\smallskip

%Somewhat surprisingly, $\hMm$-actions are better behaved than $\hGm$-actions. We \emph{\textcolor[rgb]{0.00,0.00,1.00}{should be able to prove}} the following
%\begin{thm}
%In the adjunction $R\vdash L$ (see Subsection \ref{sec:monoids}):
%\[
%R:\cinfsalg_{\hmM}\leftrightarrows\comsalg_\mM\times_{\comsalg}\cinfsalg:L,
%\]
%$L$ is fully faithful while $R$ is faithful and essentially surjective.
%\end{thm}
%\begin{proof}
%We must prove that the unit of this adjunction is an isomorphism (easy), and that the counit is surjective.
%\end{proof}

\begin{eg} In spite of the above remarks, there exist $\hmG$-coactions for any super Fermat theory $\bS$ having rings of formal power series as algebras, which are not essentially algebraic in the sense we defined. Consider the $\bS$-algebra of formal power series $\A=\KK[[x]]$ equipped with the standard dilation action
\[
\Phi:\KK[[x]]\To\KK\{u\}[[x]],\quad x\mapsto ux
\]
Then $\A^0=\KK$, $\A_\alg=\KK[x]$, but $\widehat{(\A_\alg)}=\KK\{x\}$ which maps into $\A$ but not isomorphically. Nevertheless, we can still say that $\A$ is \emph{some} -- in this case, formal -- completion of $\A_\alg$.
\end{eg}

%When $\A^0$ is $\bLambda$-point determined, we \emph{\textcolor[rgb]{0.00,0.00,1.00}{should be able to}} say more:

%\begin{thm}
%The $((\quad)_\mathrm{alg},\widehat{(\quad)})$ adjunction restricts to a fully faithful reflective embedding
%\[
%(\quad)_\mathrm{alg}:\cinfsalg_{\mathbf{\Lambda}\mathrm{det},\hmM}\overleftarrow{\hookrightarrow}\comsalg_\mM:\widehat{(\quad)}.
%\]
%In particular, for a $\hMm$-algebra $\A$ that is $\bLambda$-point determined, the counit $\widehat{\A_\mathrm{alg}}\to\A$ is an isomorphism.
%\end{thm}

%(Part of this is showing that if $\A^0$ is $\bLambda$-point determined, then so is all of $\A$.)

%The invariant part of a $\hMm$-algebra $\A$ is easy to calculate, thanks to having a zero around. We have a map
%\[
%\pi=(\zeta\oinfty\id)\circ\Phi:\A\To\A
%\]
%which is an idempotent whose image turns out to be precisely $\A^0$. So $\pi$ splits as
%\[
%\A^0\leftrightarrows\A.
%\]

We can describe binary coproducts in $\gr\bS\Alg$ explicitly (see Appendix \ref{sec:coactions} for a general discussion). Let $\A,\B$ be graded $\bS$-algebras. We first observe that we have an inclusion of (super)commutative algebras
\[
\A^0\otimes\B^0\To(\A\otimes\B)^0=\bigoplus_{i+j=0}\A^i\otimes\B^j.
\]
Form an algebra $(\A\otimes\B)^{0'}$ as the pushout of
\[
\A^0\odot\B^0\oT\A^0\otimes\B^0\To(\A\otimes\B)^0.
\]
Then $(\A\otimes\B)^{0'}$ is naturally an $\A^0\odot\B^0$-algebra. Define $(\A\circledast\B)^0$ to be its $\bS$-algebra completion relative to $\A^0\odot\B^0$, and finally, let the coproduct $\A\circledast\B$ be the pushout of
\[
(\A\circledast\B)^0\oT(\A\otimes\B)^0\To\A\otimes\B.
\]
When both $\A$ and $\B$ are non-negatively (or non-positively) graded, the construction simplifies drastically, for then we have
\[
\A^0\otimes\B^0\cong(\A\otimes\B)^0,
\]
and hence $\A\circledast\B$ is just the pushout of
\[
\A^0\odot\B^0\oT\A^0\otimes\B^0\To\A\otimes\B,
\]
with $(\A\circledast\B)^0=\A^0\odot\B^0$. More general pushouts of graded $\bS$-algebras are defined analogously.

\begin{eg}
Let $M$ and $N$ be $\NN$-graded manifolds as defined in \cite{roygr}, with base manifolds $M_0$ and $N_0$, respectively, and let $\C(M)$ and $\C(N)$ be their algebras of polynomial functions. They are non-negatively graded $\bcinf$-algebras in the sense we have defined. It is easy to see that $\C(M)\circledast\C(N)=\C(M\times N)$: in the construction, one completes $$\C(M)^0\otimes\C(N)^0=\cinf(M_0)\otimes\cinf(N_0)$$ to $\cinf(M_0)\oinfty\cinf(N_0)=\cinf(M_0\times N_0)$ and then applies the change of base.

Now, if we apply to $\C(M)$ the smooth completion relative to $\C(M)^0=\cinf(M_0)$, we get the $\bscinf$-algebra $\cinf(\hat{M})$ of smooth functions on the supermanifold $\hat{M}$ equipped with a smooth action of the multiplicative monoid $\hMm$, which is an $\NN$-graded manifold in the sense defined in \cite{Sev}. Clearly we have $$\cinf(\hat{M})\oinfty\cinf(\hat{N})=\cinf(\hat{M}\times\hat{N})$$ and, since completions preserve coproducts, also $$\widehat{M\times N}=\hat{M}\times\hat{N}.$$ Furthermore, in view of Remark \ref{rem:Mmact} we can conclude that the categories of $\NN$-graded manifolds as defined in \cite{roygr}  and in \cite{Sev} are equivalent. On the level of algebras, this equivalence is induced by the adjunction
\[
\Adj{\tau^\dagger_\mM}{\gr\bcinf\Alg_{\geq0}=\bscinf\Alg^{\hmM}_{\alg,\even}}{\bscinf\Alg^{\hmM}_\even.}{\tau^\mM_\dagger}
\]
\end{eg}

\subsection{The Euler derivation.}

The Lie algebra of $\Gm$ (and hence of $\hGm$) is one-dimensional abelian. The infinitesimal generator of a $\hmG$-coaction is commonly known as the \emph{Euler derivation} and denoted by $\epsilon:\A\to\A$. The coaction is essentially algebraic if and only if the operator $\epsilon$ is diagonalizable with integer eigenvalues (which are the weights) and complete (integrates to a flow).  ``Diagonalizable'' here means that $\A$ is generated as an $\bS$-algebra (over $\A^0$) by the graded subalgebra formed by the eigenspaces $\A^n$, $n\in\ZZ$. Indeed, it is only in this case that the name ``Euler derivation'' is appropriate.

Let us now give an example of calculating $\A_\alg$ for an $\bS$-algebra $\A$, which shows that we get the expected answer; in fact, this is the motivating example for introducing the algebraization functor.

\begin{eg}
Let
\[
\A=\KK\{x_0,x_1,x_2\ldots\}
\]
(the set of generators can be finite or infinite; it suffices to consider only even generators). Define the $\Mm$-action via the Euler derivation
\[
\epsilon=\sum_{k\geq0}kx_k\frac{\del}{\del x_k}
\]
In other words,
\[
\Phi(f)(x_0,x_1,x_2,\ldots)=f(x_0,ux_1,u^2x_2,\ldots)
\]
We claim that
\[
\A_\alg=\KK\{x_0\}[x_1,x_2,\ldots].
\]
Indeed, $f\in\A_\alg$ if and only if
\[
f(x_0,ux_1,u^2x_2,\ldots)=\sum_{k=0}^Na_ku^k
\]
for some $N\in\NN$ and $a_k\in\A$, $k=0,\ldots,N$. Taking the $n^{th}$ derivative with respect to $x_i$ for each $i>0$ on both sides, we obtain
\[
u^{in}f^{(n)}_i(x_0,ux_1,u^2x_2,\ldots)=\sum_{k=0}^Na^{(n)}_{k,i}u^k,
\]
where for a $g\in\A$, $g^{(n)}_i$ denotes the $n^{th}$ partial derivative of $g$ with respect to $x_i$. Comparing the right and left-hand sides, we conclude that $f^{(n)}_i$ must vanish identically for $in>N$. Since $f$ depends only on finitely many of the generators (by definition of the free $\bS$-algebra), we conclude that \emph{all} partial derivatives of $f$ of sufficiently high order with respect to the generators of positive degree must vanish. Therefore, by Taylor's formula (cf. \cite{1forms} or \cite{dg1}, Corollary \ref{cor:taylor}), we conclude that $f$ is polynomial in those variables, and we are done.

Notice that the invariant subalgebra $\A^0$ of $\A$ is $\KK\{x_0\}$, and the $\bS$-completion of $\A_\alg$ relative to $\A^0$ is isomorphic to $\A$.

Observe finally that the same argument applies more generally, for $\A$ the free $\bS$-algebra on a generating set
\[
S=\coprod_{k\geq0}S_k
\]
where $S_k$ are arbitrary sets and the $\Mm$-action is defined by assigning each generator corresponding to an element of $S_k$ degree $k$.
\end{eg}

\subsection{The odd line and its endomorphisms.}

By the \emph{odd (affine) line} we mean the functor $\KK^{0|1}:\bscom_\KK\Alg\to\Set$ associating to each $\A\in\bscom_\KK\Alg$ its odd part $\A_\odd$. Notice that this functor actually takes values in $\KK$-modules; the underlying abelian group-valued functor is the odd additive group $\Gad^\odd$ considered earlier. We shall, however, distinguish the $\Set$-valued functor $\KK^{0|1}$ from $\Gad^\odd$; the former is a torsor over the latter, justifying the name ``odd affine line''. The functor $\KK^{0|1}$ is co-represented by the superalgebra $\Lambda^1=\KK[\theta]$ (with $\theta$ odd), which is to be viewed as the ``algebra of functions'' on $\KK^{0|1}$.

As any Weil superalgebra, $\Lambda^1$ is co-exponentiable (\cite{dg1}, Proposition \ref{prop:coexp}), the functor $F^{\KK^{0|1}}$ is co-representable for every co-representable $$F:\bscom_\KK\Alg\to\Set.$$ In particular, $\EE={\KK^{0|1}}^{\KK^{0|1}}$ is co-representable. To find the co-representing algebra $\E$, just observe that, for any superalgebra $\A$, any homomorphism $$\phi:\KK[\theta]\to\A[\tau]=\A\otimes\KK[\tau]$$ is of the form
\begin{equation}\label{eq:endo}
\phi(\theta)=\alpha+a\tau,
\end{equation}
with $\alpha\in\A_\odd$, $a\in\A_\even$, and conversely, every such pair $(\alpha,a)$ defines a map, since $\KK[\theta]$ is free on one odd generator. It follows that $\E=\KK[u,\nu]$ (with $u$ even, $\nu$ odd).

Since $\EE$ is obviously a monoid in $\Set^{\bscom_\KK\Alg}$, $\E$ is a comonoid in $\bscom_\KK\Alg$. One easily computes the comultiplication to be
\[
\Delta:\E=\KK[u,\nu]\To\KK[u',u'',\nu',\nu'']=\E\otimes\E,\quad u\mapsto u'u'',\quad\nu\mapsto\nu'+u'\nu''.
\]
and counit
\[
c:\KK[u,\nu]\To\KK,\quad u\mapsto 1,\quad\nu\mapsto 0.
\]
Thus, we have the semidirect product decomposition
\[
\EE\cong\Mm\ltimes\Gad^\odd.
\]
In other words, all endomorphisms of the odd line are affine.

Observe now that the endomorphism $\phi$ in \eqref{eq:endo} is invertible if and only if $a$ is a unit in $\A$. Hence, the maximal subgroup $\EE^\times\subset\EE$ of automorphisms of $\KK^{0|1}$ is
\[
\EE^\times\cong\Gm\ltimes\Gad^\odd.
\]
It is corepresented by $\E^\times=\KK[u,u^{-1},\nu]$ with the same formulas for the comultiplication and counit, and the antipode given by
\[
S:\KK[u,u^{-1},\nu]\To\KK[u,u^{-1},\nu],\quad u\mapsto u^{-1},\quad\nu\mapsto-\nu.
\]

The Lie algebra $\ee$ of $\EE$ is the Lie superalgebra $\Der(\KK[\theta])$ of $\KK$-linear derivations of $\KK[\theta]$. As a $\KK$-module, it is isomorphic to $\KK^{1|1}$, with the even generator $\eps=-\theta\frac{d}{d\theta}$ of dilations and odd generator $\delta=\frac{d}{d\theta}$ of translations, subject to bracket relations
\begin{equation}\label{eq:brackets}
[\eps,\delta]=\delta,\quad [\eps,\eps]=[\delta,\delta]=0.
\end{equation}
(The minus sign in the definition of $\eps$ is a matter of convention; it corresponds to the inversion in the $\Gm$-component of $\EE.$)

\subsection{Differential graded superalgebras.}

Consider an $\EE^\times$-action ($\E^\times$-coaction) on a superalgebra $\A\in\bscom_\KK\Alg$. The action of the subgroup $\Gm\subset\EE^\times$ induces a $\ZZ$-grading on $\A$, while the action of $\Gad^\odd\subset\EE^\times$ gives an odd-time flow along an odd derivation $d$. The corresponding infinitesimal action is a homomorphism of Lie superalgebras
\[
\ee\To\Der(\A),\quad\eps\mapsto\epsilon,\quad\delta\mapsto d,
\]
where $\epsilon$ is the Euler derivation inducing the weight grading. The relations \eqref{eq:brackets} imply that the operator  $d$ is a differential (since $d^2=\frac{1}{2}[d,d]=0$), and that it increases the weight by $1$ (since $[\epsilon, d]=d$).

\begin{defn}
Such an action is said to be \emph{even} if the action of $\mu_2\subset\Gm\subset\EE^\times$ (corresponding to the $\ev_{(-1,0)}:\E^\times\to\KK$) coincides with the Grassmann parity involution. (Recall that this simply means that the modulo-$2$ reduction of the weight grading coincides with the Grassman parity.)
\end{defn}

In the case of an even action, the fact that the differential $d$ is of degree $+1$ already implies that it is Grassmann odd (with respect to sign rules). We thus recover the known fact that an even $\EE^\times$-action on a superalgebra $\A$ is the same thing as a differential graded commutative algebra structure on $\A$. In other words, we have an equivalence of categories
\begin{equation}\label{eq:dgeven}
\bscom_\KK\Alg^{\E^\times}_\even\simeq\dg\bcom_\KK\Alg,
\end{equation}
between supercommutative $\KK$-algebras with an even coaction of $\E^{\times},$ and the category of differential graded $\KK$-algebras.
For general (non-even) actions, the weight grading is independent of the Grassmann parity, the differential $d$ is both odd and of degree $+1$, and we have
\begin{equation}\label{eq:dg}
\bscom_\KK\Alg^{\E^\times}\simeq\dg\bscom_\KK\Alg.
\end{equation}

\begin{rem}
The minus sign in our choice of the generator $\eps$ of dilations ensures, via the relations \eqref{eq:brackets}, that the differential increases the degree. Without the sign, we would have $[\epsilon,d]=-d$, so $d$ would \emph{decrease} the degree instead, leading to \emph{chain} instead of cochain complex grading conventions.

The same minus sign also means that the restricted $\mG$-coaction maps an $a\in\A^n$ to $u^{-n}a$.
\end{rem}

If the $\EE^\times$ action extends to an $\EE$-action, then, according to our convention, only \emph{non-positive} weights occur. This leads to equivalences of categories
\begin{equation}\label{eq:dgneg}
\bscom_\KK\Alg^{\E}\simeq\dg\bscom_\KK\Alg_{\leq0}
\end{equation}
and
\begin{equation}\label{eq:dgnegeven}
\bscom_\KK\Alg^{\E}_\even\simeq\dg\bcom_\KK\Alg_{\leq0}.
\end{equation}
Dually, if the corresponding $(\EE^\times)^\op$-action extends to an $\EE^\op$-action, only \emph{non-negative} weights occur, and we have
\begin{equation}\label{eq:dgpos}
\bscom_\KK\Alg^{\E^\op}\simeq\dg\bscom_\KK\Alg_{\geq0}
\end{equation}
and
\begin{equation}\label{eq:dgposeven}
\bscom_\KK\Alg^{\E^\op}_\even\simeq\dg\bcom_\KK\Alg_{\geq0}.
\end{equation}

\begin{rem}
Since the monoid $\EE$ is neither commutative nor a group, $\EE$ and $\EE^\op$-actions are genuinely different from each other.
\end{rem}

To conclude this subsection, observe that the $\E^\times$-invariant subalgebra $\A_{\E^\times}\subset\A$ consists of elements annihilated \emph{both} by $\epsilon$ \emph{and} by $d$. In other words,
\[
\A_{\E^\times}=\A^{0,\cl},
\]
the subalgebra of $0$-cocycles. The same holds for $\E$ and $\E^\op$-coactions, except for in $\E$-coactions, no positive weights occur, so $\A^{0,\cl}$ is \emph{all} of $\A^0$. Moreover, in this case, the differential $d$ is $\A^0$-linear (in general, it is only $\A^{0,\cl}$-linear).

\subsection{Differential graded $\bS$-algebras.}

The functors $\KK^{0|1}$, $\EE$ and $\EE^\times$ extend along the forgetful functor to functors $\hat\KK^{0|1}$, $\hat\EE$ and $\hat\EE^\times$ on $\bS\Alg$ co-represented, respectively, by $\KK\{\theta\}=\Lambda^1$, $\KK\{u,\nu\}$ and $\KK\{u,u^{-1},\nu\}$. We have semidirect product decompositions
\[
\hat\EE\cong\hMm\ltimes\hGad^\odd
\]
and
\[
\hat\EE^\times\cong\hGm\ltimes\hGad^\odd
\]
with the formulas defining the comonoid and cogroup structures exactly the same as in the algebraic case. The Lie superalgebra $\ee$ given by \eqref{eq:brackets} is also the Lie algebra of $\hat\EE^\times$ of derivations of $\KK\{\theta\}=\Lambda^1$, viewed as an $\bS$-algebra.

$\hat\E^\times$-actions differ from $\E^\times$-actions just insofar as $\hGm$-actions differ from $\Gm$-actions. Given an $\hat\E^\times$-coaction on an $\A\in\bS\Alg$, we still have the differential $d$ and an Euler derivation $\epsilon$ satisfying the Lie superalgebra relations \eqref{eq:brackets} of $\ee$, only now $\epsilon$ can be quite badly behaved and not determined by any graded subalgebra of $\A$. It is still true, however, that
\[
\A_{\hat\E^\times}=\A^{0,\cl}
\]
and that, for $\hat\E$-coactions, $\A^{0,\cl}$ coincides with $\A^0$. As in Proposition \ref{prop:adjalg}, we have the adjunctions:
\[
\Adj{\tau_{\E^\times}^*}{\dg\bscom_\KK\Alg}{\bS\Alg^{\hat\E^\times}}{\tau^{\E^\times}_!},
\]
\[
\Adj{\tau_{\E}^*}{\dg\bscom_\KK\Alg_{\leq0}}{\bS\Alg^{\hat\E}}{\tau^{\E}_!}
\]
and
\[
\Adj{\tau_{\E^\op}^*}{\dg\bscom_\KK\Alg_{\geq0}}{\bS\Alg^{\hat\E^\op}}{\tau^{\E^\op}_!}
\]
(taking into account the equivalences \eqref{eq:dg}, \eqref{eq:dgneg} and \eqref{eq:dgpos}), which decompose as
\begin{equation}\label{eq:adjdg}
\Adjadj{\tau_{\E^\times}^\circ}{\dg\bscom_\KK\Alg}{\bS\Alg^{\hat\E^\times}_\alg}{\tau^{\E^\times}_\circ}{\tau_{\E^\times}^\dagger}{\bS\Alg^{\hat\E^\times}}{\tau^{\E^\times}_\dagger},
\end{equation}
\[
\Adjadj{\tau_{\E}^\circ}{\dg\bscom_\KK\Alg_\leq0}{\bS\Alg^{\hat\E}_\alg}{\tau^{\E}_\circ}{\tau_{\E}^\dagger}{\bS\Alg^{\hat\E}}{\tau^{\E}_\dagger}
\]
and
\[
\Adjadj{\tau_{\E^\op}^\circ}{\dg\bscom_\KK\Alg_\geq0}{\bS\Alg^{\hat\E^\op}_\alg}{\tau^{\E^\op}_\circ}{\tau_{\E^\op}^\dagger}{\bS\Alg^{\hat\E^\op}}{\tau^{\E^\op}_\dagger},
\]
respectively (cf. Propositions \ref{prop:adjalg1} and \ref{prop:adjalg2}).

When $\bS=\bSE$, the superization of a Fermat theory $\bE$ (so the ground ring $\KK$ is purely even), it also makes sense to consider even coactions, for which we have the adjunctions
\[
\Adjadj{\tau_{\E^\times}^\circ}{\dg\bcom_\KK\Alg}{\bSE\Alg^{\hat\E^\times}_{\even,\alg}}{\tau^{\E^\times}_\circ}{\tau_{\E^\times}^\dagger}{\bSE\Alg^{\hat\E^\times}_\even}{\tau^{\E^\times}_\dagger},
\]
\[
\Adjadj{\tau_{\E}^\circ}{\dg\bcom_\KK\Alg_{\leq0}}{\bSE\Alg^{\hat\E}_{\even,\alg}}{\tau^{\E}_\circ}{\tau_{\E}^\dagger}{\bSE\Alg^{\hat\E}}{\tau^{\E}_\dagger}
\]
and
\[
\Adjadj{\tau_{\E^\op}^\circ}{\dg\bcom_\KK\Alg_{\geq0}}{\bSE\Alg^{\hat\E^\op}_{\even,\alg}}{\tau^{\E^\op}_\circ}{\tau_{\E^\op}^\dagger}{\bSE\Alg^{\hat\E^\op}}{\tau^{\E^\op}_\dagger}.
\]

The objects of the category $\bS\Alg^{\hat\E^\times}_\alg$ are differential $\ZZ$-graded $\KK$-algebras $\A$ equipped with an extra $\bS$-algebra structure on $\A^{0,\cl}$ (and similarly for $\hat\E$ and $\hat\E^\op$). The functor $\tau^{\E^\times}_\dagger$ takes such an algebra to its $\bS_{\A^{0,\cl}}$-completion, while the right adjoint $\tau_{\E^\times}^\dagger$ takes an $\hat\E^\times$-algebra $\A$ to the differential graded $\KK$-algebra $\A_\alg$ assembled from the components of $\A$ corresponding to the algebraic characters of $\hGm$ (i.e. the integers), while remembering the $\bS$-algebra structure on $\A^{0,\cl}$.

We are finally ready to define what it means to have a ``differential graded structure'' on an $\bS$-algebra.

\begin{defn}
A \emph{differential} $\ZZ$ (resp. $\ZZ_{\leq0}$, $\ZZ_{\geq0}$) \emph{-graded $\bS$-algebra} is an object of $\bS\Alg^{\E^\times}_\alg$ (resp. $\bS\Alg^{\E}_\alg$, $\bS\Alg^{\E^\op}_\alg$). In other words, it is a differential graded superalgebra $\A$ whose subalgebra $\A^{0,\cl}$ of zero-cocycles is equipped with an additional $\bS$-algebra structure. A morphism of differential graded $\bS$-algebras $$\A \to \B,$$ is a morphism of underlying differential graded superalgebras such that the induced morphism $$\A^{0,\cl} \to \B^{0,\cl}$$ is a morphism of $\bS$-algebras. When $\bS=\bSE$, we speak of a \emph{differential graded $\bE$-superalgebra}, or, when the coaction is even, a \emph{differential graded $\bE$-algebra}.

A \emph{differential} $\ZZ$ (resp.$\ZZ_{\leq0}$, $\ZZ_{\geq0}$) -\emph{graded structure} on an $\bS$-algebra $\B$ is an essentially algebraic $\E^\times$ (resp. $\E$, $\E^\op$) -coaction on $\B$. In other words, a coaction and an isomorphism from $\B$ to $\tau^{\E^\times}_\dagger(\A)$ (resp. $\tau^{\E}_\dagger(\A)$, $\tau^{\E^\op}_\dagger(\A)$) for some differential graded $\bS$-algebra $\A$.
\end{defn}

\begin{notation}
We shall henceforth denote the category $\bS\Alg^{\E^\times}_\alg$ (resp. $\bS\Alg^{\E}_\alg$, $\bS\Alg^{\E^\op}_\alg$) by $\dg\bS\Alg$ (resp. $\dg\bS\Alg_{\leq0}$, $\dg\bS\Alg_\geq0$).

When $\bS=\bSE$, we shall use $\dg\bE\Alg$ for $\bSE\Alg^{\E^\times}_{\alg,\even}$, and similarly for the bounded degree cases.

%Denote the category $\bS\Alg^{\mG}_\alg$ (resp. $\bS\Alg^{\mM}_\alg$, $\bS\Alg^{\mM^\op}_\alg$) by $\gr\bS\Alg$ (resp. $\gr\bS\Alg_{\leq0}$, $\gr\bS\Alg_\geq0$; %of course, the latter two categories are equivalent).
\end{notation}

\begin{rem}
Since the functors $\tau^{\E^\times}_\dagger(\A)$, $\tau^{\E}_\dagger(\A)$ and $\tau^{\E^\op}_\dagger(\A)$ are fully faithful, the categories of differential graded $\bS$-algebras, and of $\bS$-algebras with a differential graded structure, are equivalent.
\end{rem}

As in the graded case, we can describe coproducts (and more general pushouts) in $\dg\bS\Alg$ (and the related categories) quite explicitly. Given differential graded $\bS$-algebras $\A$ and $\B$, we first form the pushout $(\A\otimes\B)^{0'}$ of
\[
\A^{0,\cl}\odot\B^{0,\cl}\oT\A^{0,\cl}\otimes\B^{0,\cl}\To(\A\otimes\B).^{0,\cl}
\]
Then, let $(\A\circledast\B)^{0,\cl}$ be the completion of $(\A\otimes\B)^{0'}$ relative to $\A^{0,\cl}\odot\B^{0,\cl}$, and finally define $\A\circledast\B$ to be the pushout of
\[
(\A\circledast\B)^{0,\cl}\oT(\A\otimes\B)^{0,\cl}\To\A\otimes\B.
\]
When both $\A$ and $\B$ are non-\emph{positively} graded, we have $\A^{0,\cl}=\A^0$ and $$\A^0\otimes\B^0=(\A\otimes\B)^0,$$ hence $$(\A\circledast\B)^0=\A^0\odot\B^0$$ and $\A\circledast\B$ is the pushout of
\[
\A^0\odot\B^0\oT\A^0\otimes\B^0\To\A\otimes\B,
\]
as in the graded case (without differentials).

\subsection{Cohomology.}

Given a differential $\bS$-algebra $(\A,d)\in\dd\bS\Alg$, we can define its cohomology in the usual way:
\[
H(\A,d)=\Ker(d)/\Im(d).
\]
Observe that $\Ker(d)$ is an $\bS$-subalgebra of $\A$: in fact, $\Ker(d)$ is the equalizer of the section $\id\oplus\tau d$ of the canonical projection
\[
\A\oplus\tau\A\To\A
\]
and the zero section. Furthermore, since $d$ is an $\bS$-algebra derivation, it is in particular a derivation of $\A$ as a $\KK$-algebra. Therefore, $\Im(d)$ is an ideal in $\Ker(d)$, so it defines an $\bS$-congruence, and hence the quotient has a canonical $\bS$-algebra structure. Moreover, every morphism of differential $\bS$-algebras clearly induces an $\bS$-algebra morphism on cohomology. We conclude that cohomology is a functor
\[
H:\dd\bS\Alg\To\bS\Alg.
\]
Now suppose $\A\in\bS\Alg^{\hat\E^\times}$. We can take its cohomology with respect to the differential $d$. Since $\hGm$ is normal in $\hat\EE^\times$, $H(\A,d)$ inherits a $\hGm$-action. Thus we get a functor
\[
H:\bS\Alg^{\hat\E^\times}\To\bS\Alg^{\hmG}.
\]
Similarly, we conclude that cohomology defines a functor
\[
H^\bullet_\alg:\dg\bS\Alg\To\gr\bS\Alg.
\]
Here, we only need to observe that, since $\A^{0,\cl}$ is an $\bS$-algebra and $d(\A^{-1})\subset\A^{0,\cl}$ is an ideal, $H^{0}(\A,d)$ inherits an induced $\bS$-algebra structure.

Finally, we have the usual cohomology functor as a special case:
\[
H^\bullet:\bscom_\KK\Alg^{\E^\times}\simeq\dg\bscom_\KK\Alg\To\gr\bscom_\KK\Alg\simeq\bscom_\KK\Alg^{\mG}.
\]

There are also versions of these functors for $\EE$ and $\EE^\op$ actions, and for $\bS=\bSE$, and also even versions. We leave these to the reader.

\begin{prop}\label{prop:coh}
The diagram of functors
$$\xymatrix@C=2cm{\dg\bscom_\KK\Alg \ar[d]_-{H^\bullet} &  \dg\bS\Alg \ar[l]_-{\tau_{\E^\times}^\circ} \ar[d]_-{H^\bullet_\alg} & \bS\Alg^{\hat\E^\times} \ar[l]_-{\tau_{\E^\times}^{\dagger}} \ar[d]^-{H}\\
\gr\bscom_\KK\Alg & \gr\bS\Alg \ar[l]^-{\tau_\mG^\circ} & \bS\Alg^\hmG \ar[l]^-{\tau_\mG^\dagger}}$$
commutes up to natural isomorphisms. Likewise for $\E$, $\E^\op$ and (for $\bS=\bSE$) even versions.
\end{prop}
\begin{proof}
The diagrams commute because the right adjoint functors involved either forget the extra $\bS$-algebra structure or restrict to a subset of characters of $\hGm$ (namely, the integers). The forgetful functor preserves equalizers and quotients by ideals (i.e. regular epimorphisms, a particular case of sifted colimits), hence cohomology, whereas the subset $\ZZ$ of characters is preserved and reflected by $d$ (``reflected'' means that if $a$ is not in $\A^n$ for any $n\in\ZZ$, then neither is $da$).
\end{proof}

\begin{rem}
The left adjoints (completions) generally fail to commute with cohomology. This is because taking cohomology involves computing equalizers, and those are generally destroyed by completion. The following counterexample is based on an example we learned from Robert Bryant. Consider the polynomial vector field
\[
v=(x^2+1)\frac{\del}{\del x}+2xy\frac{\del}{\del y}
\]
on $\RR^2$. It defines a differential $\ZZ_{\geq0}$-graded algebra structure on
\[
\A=\RR[x,y,\xi]
\]
with
\[
\epsilon=\xi\frac{\del}{\del\xi},\quad d=\xi v.
\]
The $0^{th}$ cohomology $H^0(\A,d)$ is the subalgebra of $\RR[x,y]$ consisting of the polynomial integrals of motion of $v$, but it is easy to see that there are no non-constant ones, so $H^0(\A,d)=\RR$. Now consider the completion $\hat\A=\RR\{x,y,\xi\}$ of $\A$ as an $\bS\bcomega$-algebra, with the same grading and differential. The function
\[
u(x,y)=\frac{y}{x^2+1}\in\RR\{x,y\}
\]
is easily seen to be an integral of motion for $v$. Therefore, $H^0(\hat\A,d)$ contains at least $\RR\{u\}$, whereas $\widehat{H^0(\A,d)}=\hat\RR=\RR$.
\end{rem}

\section{Differential forms.}\label{sec:diffforms}
\subsection{The algebra of differential forms.}
Given a superalgebra $\A\in\bscom_\KK\Alg$, its superalgebra of differential forms $\om(\A)$ can be defined as the universal differential superalgebra generated by $\A$. This universal property can be expressed in the more general setting of $\bS$-algebras by interpreting a differential as a $\adG^\odd$-coaction and taking advantage of the co-exponentiability of $\adG^\odd=\Lambda^1$ in $\bS\Alg$.

\begin{defn}
Given an $\bS$-algebra $\A,$ we define its algebra of differential forms to be $$\wom(\A):=\cop{\A}{\adG^\odd},$$ via co-exponentiation.
\end{defn}

\begin{rem}
Co-exponentiation corresponds to exponentiation in the opposite category, so one may think of the algebra of differential forms $\wom(\A)$ as the algebra of functions on the affine $\bS$-scheme $\Spec_{\bS}\left(\A\right)^{\KK^{0|1}}.$
\end{rem}

\begin{eg}
In particular, if $\bS=\bS\bcinf,$ and $\A= \bcinf\left(M\right)$ is smooth functions on a (super) manifold, then $M^{\RR^{0|1}}=\Pi TM,$ the odd tangent bundle, and $$\wom(\bcinf\left(M\right))=\bcinf\left(\Pi TM\right)$$ which is by definition, due to Bernstein-Leites, the algebra of pseudo-differential forms on $M$ (see also the discussion in \cite{pfield} Section 3.2). Moreover, it readily follows that $M^{\RR^{0|1}}=\Pi TM$ admits an action of the smooth endomorphism monoid $\End\left(\RR^{0|1}\right),$ hence $\wom(\bcinf\left(M\right))$ has an algebraic coaction of the comonoid
$$\hat\EE =\bcinf\left(\End\left(\RR^{0|1}\right)\right),$$
making it a differential graded $\bS\bcinf$-algebra. The algebraic part of this $\hat\EE$-comodule, $$\wom(\bcinf\left(M\right))_{\alg},$$ corresponds to the subalgebra of $\bcinf\left(\Pi TM\right)$ consisting of those smooth functions polynomial in the fiber coordinates, and inherits the structure of a differential graded (super) $\RR$-algebra, which is the classical dg-algebra of differential forms. As we shall see, this is a general phenomena.
\end{eg}

Indeed, by the general Remark \ref{rem:freecoact}, $\wom(\A)=\cop{\A}{\adG^\odd}$ carries a canonical $\adG^\odd$-coaction, and the functor
\[
\wom:\bS\Alg\To\bS\Alg^{\adG^\odd}
\]
is left adjoint to the forgetful functor to $\bS\Alg$ (and similarly for $\om$ in the algebraic case). The notation $\wom$ is justified by observing that co-exponentiation is a left adjoint, hence commutes with (relative) completion. For the same reason $\om$ and $\wom$ take free algebras to free algebras. Explicitly, we have
\[
\om(\KK[x^1,\ldots,x^m\xi^1,\ldots,\xi^n])=\KK[x^1,\ldots,x^m,d\xi^1,\ldots,d\xi^n;\xi^1,\ldots,\xi^n,dx^1,\ldots,dx^m]
\]
and
\[
\wom(\KK\{x^1,\ldots,x^m\xi^1,\ldots,\xi^n\})=\KK\{x^1,\ldots,x^m,d\xi^1,\ldots,d\xi^n;\xi^1,\ldots,\xi^n,dx^1,\ldots,dx^m\}.
\]
Denote these algebras by $\om_{m|n}$ and $\wom_{m|n}$, respectively (and for $n=0$, $\om_m$ and $\wom_m$). The differential $d$ is odd and takes $x^\mu$ to $dx^\mu$, and $dx^\mu$ to $0$, and similarly for the $\xi^i$'s. By Clairaut's theorem (which holds for every (super) Fermat theory, cf. \cite{dg1}, Proposition \ref{prop:clairaut}), $d^2=0$.

The description of $\wom(\A)$ as $\cop{\A}{\adG^\odd}$ reveals that $\wom(\A)$ carries not only a $\adG^\odd$-coaction, but an $\hat\E^\op$-coaction, and $\wom(\A)$ is also universal with this property. Since $\wom(\A)$ is the completion of $\om(\A)$, this coaction is algebraic; the Euler derivation is defined by setting
\[
\epsilon(a)=0,\quad\epsilon(da)=da\quad\forall a\in\A.
\]
For $\wom_{m|n}$ we have
\[
\epsilon=dx^\mu\frac{\del}{\del dx^\mu}+d\xi^i\frac{\del}{\del d\xi^i},\quad d=dx^\mu\frac{\del}{\del x^\mu}+d\xi^i\frac{\del}{\del\xi^i}.
\]
The degree-$1$ submodule $\wom^1(\A)$ coincides with the module $\Pi\Omega^1\left(\A\right),$ the recipient of the universal odd derivation, as in Definition \ref{dfn:kahlerforms}.

\begin{rem}
Some authors refer to elements the subalgebra $\wom_\alg(\A)$ as differential forms proper (since they are polynomial in the $d\alpha$'s for $\alpha\in\A_\odd$), while referring to general elements of $\wom(\A)$ as \emph{pseudo-differential forms}, as these include non-polynomial functions like $e^{d\alpha}$ (for $\bS=\bS\bcinf$ or $\bS\bcomega$).
\end{rem}

The co-exponential description of $\wom(\A)$ also produces, for every derivation $D$ of $\A$ (even or odd), the derivations $\iota_D$ (contraction, of parity opposite to that of $D$) and $L_D$ (Lie derivative, of the same parity as $D$) of $\wom(\A)$, satisfying the usual commutation relations:
\[
[\epsilon,\iota_D]=-\iota_D,\quad[\epsilon,L_D]=0,\quad[d,\iota_D]=L_D,
\]
\[
[\iota_D,\iota_{D'}]=0,\quad[L_D,\iota_{D'}]=\iota_{[D,D']},\quad[L_D,L_{D'}]=L_{[D,D']}.
\]

\subsection{Homotopy invariance.}
\begin{defn}\label{def:hoinv}
We say that a super Fermat theory $\bS$ satisfies the \emph{homotopy invariance axiom} if for every $\A\in\bS\Alg^{\adG^\odd}$, the canonical inclusion
\[
\A\To\A\odot\wom_1=\A\{t,dt\}
\]
induces an isomorphism on cohomology. We say that $\bS$ satisfies the \emph{graded homotopy invariance} axiom if for all $\A\in\bS\Alg^{\E^\times}_\alg$, the canonical inclusion
\[
\A\To\A\circledast\wom_{1,\alg}
\]
is a homology isomorphism (where $\circledast$ denotes the coproduct in $\bS\Alg^{\E^\times}_\alg$).
\end{defn}

\begin{rem}\label{rem:poinc}
The homotopy invariance axiom (resp. graded homotopy invariance axiom) implies the Poincar\'e lemma for $\wom_n$ (resp. $\wom_{n,\alg}$) for all $n$.
\end{rem}

\begin{rem}
By using $\wom_{0|1}$ instead of $\wom_1$, we get the odd version of the homotopy invariance axiom. Both the even and odd versions together imply the Poincar\'e Lemma for $\wom_{m|n}$ for all $m$ and $n$. We shall not be concerned with the odd version in this paper.
\end{rem}
The following is standard:
\begin{prop}\label{prop:hialg}
The theory $\bscom_\KK$ satisfies the homotopy invariance axiom.

\end{prop}
\begin{proof}
Let $\A\in\dg\bscom_\KK\Alg$. We must show that the canonical inclusion
\[
\A\To\A\otimes_\KK\Omega_1=\A[t,dt]
\]
is a homology isomorphism. To this end, consider the Euler derivation $e=t\frac{\del}{\del t}$ on $\A[t]$ corresponding to the standard grading of polynomials. We claim that the corresponding contraction
\[
\iota_e=t\frac{\del}{\del dt}
\]
of differential forms provides a retraction of $\A[t,dt]$ onto $\A$. Indeed, let $D=d+dt\frac{\del}{\del t}$ be the differential on $\A[t,dt]$ (where $d$ is the differential on $\A$). Then $d$ obviously commutes with $\iota_e$, so we have
\[
L_e=[D,\iota_e]=\left[dt\frac{\del}{\del t},t\frac{\del}{\del dt}\right]=t\frac{\del}{\del t}+dt\frac{\del}{\del dt}.
\]
The module $\A[t,dt]$ decomposes into the direct sum of eigenmodules of $L_e$, corresponding to the eigenvalues $n\in\NN$. These eigenmodules are $D$-invariant and, for $n\neq0$, contractible (since $\KK\supset\QQ$); the zero-eigenmodule is $\A$, so we are done.
\end{proof}

\begin{cor}\emph{(Poincar\'e Lemma, algebraic version.)}
The initial object inclusion
\[
\KK\To\Omega_n=\KK[t^1,\ldots,t^n,dt^1,\ldots, dt^n]
\]
is a homology isomorphism.
\end{cor}

\begin{rem}
The odd version of the homotopy invariance axiom also holds for $\bscom_\KK$; in fact, the algebraic Poincar\'e Lemma holds for infinitely generated free (super)algebras. To see this, let $S=(S_\even,S_\odd)$ be a generating set and consider the Euler vector field
\[
e=\sum_{s\in S_\even}t^s\frac{\del}{\del t^s}+\sum_{\sigma\in S_\odd}\tau^\sigma\frac{\del}{\del\tau^\sigma}
\]
generating the canonical grading on $\bcom_\KK\Alg(S)$. Then use the same argument as in the proof of Proposition \ref{prop:hialg}.
\end{rem}

\subsection{Integration.}

For theories other than $\bscom_\KK$, the homotopy invariance can be deduced by standard arguments if differential forms can be integrated, and the integration operation has the expected properties (in fact, the above argument involving the Euler vector field is integration in disguise). It turns out that, to develop an integration theory sufficient for our purposes -- namely, integrals over polyhedral chains satisfying Stokes' Theorem -- it suffices to impose just one axiom: every function has an anti-derivative, and any two such differ by a constant (cf \cite{sgm}).

\begin{defn}\label{def:int}
We say that a super Fermat theory $\bS$ \emph{admits integration} if for every $f\in\bS(m+1,n)=\bS(m|n)\{t\}$ there exists an $F\in\bS(m|n)\{t\}$ such that
\[
\frac{\del F}{\del t}=f,
\]
and moreover, if $F$ and $G$ satisfy this property, then $F-G\in\bS(m|n)$.
%An \emph{integral} on a super Fermat theory $\bS$ is a an operation associating to each pair $(a,b)\in\KK\times\KK$ and each superdimension $(m|n)$ a map of %$\bS(m|n)$-modules
%\[
%\int_a^b:\bS(m|n)\odot\wom_1=\bS(m|n)\{t,dt\}\To\Pi\bS(m|n)
%\]
%satisfying the following properties
%\begin{enumerate}
%\item $\int_a^b$ has degree $-1$ with respect to the form grading,
%\item $\int_a^b=-\int_b^a$,
%\item $\int_a^b$ is natural with respect to maps $\bS(n|m)\to\bS(p|q)$, and
%\item for any $f\in\bS(m|n)\odot\wom_0=\bS(m|n)\{t\}$,
%\[
%\int_a^b\!\!df=f(b)-f(a),
%\]
%where $df=dt\frac{\del f}{\del t}$.
%\end{enumerate}
\end{defn}

If $\bS$ admits integration, given $a,b\in\bS(m|n)_\even$ and $f\in\bS(m|n)\{t\}$, we can \emph{define} the integral by forcing the fundamental theorem of calculus:
\begin{equation}\label{eqn:int}
\int_a^b\!\!dtf=F(b)-F(a)=\int_a^b\!\!dF,
\end{equation}
where we have used the canonical evaluation maps given by the $\bS$-algebra structure.

\smallskip

A number of remarks are in order.

\begin{rem}
The integral defines a (parity-preserving!) $\bS(m|n)$-linear map
\[
\int_a^b\!\!dt:\bS(m|n)\{t\}\To\bS(m|n),
\]
which is moreover natural with respect to change of base morphisms $$\bS(m|n)\to\bS(p|q).$$ Further properties -- additivity, antisymmetry, integration by parts and the change of variables formula -- are all easily deduced from \eqref{eqn:int} (exactly as it is usually done in calculus, after paying lip service to Riemann sums). The latter formula implies that the integral $\int_a^b$ is a well-defined (odd) operation on $1$-forms. With a little more work, one can define integrals of forms of higher degree over polyhedral chains and prove Stokes' theorem. We leave the details to the reader (or see \cite{sgm}).
\end{rem}

\begin{rem}
We assumed from the outset that $\KK\supset\QQ$. However, this property is \emph{implied} by the existence of an integral, which also implies that the theory $\bS$ must be reduced in the sense of \cite{dg1}.
%However, if $\KK$ is a \emph{field} and $\bS$ admits integration, it follows that $\KK$ must be of characteristic $0$. Indeed, if the characteristic is $p\neq0$, %then $dt^p=pt^{p-1}=0$, contradicting property $(4)$ of the integral.
\end{rem}

\begin{rem}
One can also speak of integration with respect to odd variables -- the so-called \emph{Berezin integral}. This can be defined for \emph{every} theory $\bS$ simply by setting
\[
\int\!\!d\theta(1)=0,\quad\int\!\!d\theta(\theta)=1
\]
and extending by linearity (notice that there are no limits of integration as they do not make sense for an odd variable). However, we shall have no use for this concept in this paper.
\end{rem}

We now observe that integrals admit base extensions:

\begin{lem}\label{lem:basextint}
Let $\bS$ be a super Fermat theory with integration, and $\B \in \bS\Alg$ any $\bS$-algebra. Then the integral extends to a map of $\B$-modules
\[
\int_a^b\!\!dt:\B\{t\}\To\B,
\]
natural with respects to morphisms of $\bS$-algebras.
\end{lem}
\begin{proof}
We first start with the finitely generated case. Suppose that $\B=\bS(m|n)/I$ for a homogeneous ideal $I$. By naturality and linearity of the integral, if $g\in I$, and $j(g)\in j_*I$ is its image under the canonical inclusion $j:\bS(m|n)\to\bS(m|n)\{t\}$, then
$$\int_a^b\!\!dt (j(g))=g\int_a^b\!\!dt(1)=(b-a)g$$
by definition of the integral. Hence $\int_a^b\!\!dt(j_*I)\subseteq I$, and so the integral descends to a $\B$-module map
\[
\int_a^b\!\!dt:\B\{t\}\To\B.
\]
To conclude the argument, we invoke the fact that arbitrary $\B$-algebras are filtered colimits of finitely-generated ones, and such colimits are computed pointwise. The naturality is clear.
\end{proof}

\begin{rem}\label{rem:egthwithint}
Examples of theories with integration are $\bcom$, $\bcinf$, $\bcomega$, $\bH$ and their superizations; the integrals are the standard ones. Moreover, by Lemma \ref{lem:basextint} it follows that if $\bS$ admits integration, so does $\bS_\B$ for any $\bS$-algebra $\B$.
\end{rem}

\begin{eg}\label{eg:noint}
An example of a Fermat theory not admitting integration is the theory $\bR^\RR$ of global rational functions with real coefficients. In fact, the Poincar\'{e} Lemma fails for forms with rational function coefficients: e.g. the closed $1$-form
\[
\alpha=\frac{dt}{1+t^2}
\]
has no rational primitive.
\end{eg}

\underline{For the remainder of this paper, let $\bS$ be a super Fermat theory with integration.} We now proceed to prove the homotopy invariance of cohomology for algebras over such a theory. The argument is completely standard -- essentially the same one used in proving the homotopy invariance of de Rham cohomology. We include it here for the reader's convenience. We start first with a definition.

\begin{defn}
Let $\A_1,\A_2\in\bS\Alg^{\adG^\odd}$ be differential $\bS$-algebras. Two maps $\phi_0,\phi_1:\A_1\to\A_2$ are said to be \emph{homotopic} if there exists a map
\[
\Phi:\A_1\To\A_2\odot\wom_1=\A_2\{t,dt\}
\]
such that $p_i\circ\Phi=\phi_i$, $i=0,1$, where $p_i$ is the evaluation at $t=i$, $dt=0$.
\end{defn}

\begin{lem}\label{lem:hoinv}
If $\phi_0$ is homotopic to $\phi_1$, $H(\phi_0)=H(\phi_1)$.
\end{lem}
\begin{proof}
We claim that the $\KK$-linear map $h=\int_0^1\circ\mspace{3mu}\Phi$ provides the requisite chain homotopy. Let $d_i$ be the differential on $\A_i$, $i=1,2$. Then the differential on $\A_2\{t,dt\}$ is $d=d_2+dt\frac{\del}{\del t}$. Let $\omega\in\A_1$. Then $\Phi\omega=\alpha+dt\beta$ for some $\alpha,\beta\in\A_2\{t\}$, and so \[
h\omega=\int_0^1\!\!dt\beta.
\]
Furthermore,
\[
\Phi d_1\omega=d\Phi\omega=\left(d_2+dt\frac{\del}{\del t}\right)\left(\alpha+dt\beta\right)=d_2\alpha+dt\left(\frac{\del\alpha}{\del t}-d_2\beta\right).
\]
Therefore,
\begin{eqnarray*}
h d_1\omega=\int_0^1\!\! dt\left(\frac{\del\alpha}{\del t}-d_2\beta\right)&=&\alpha(1)-\alpha(0)-\int_0^1\!\!dt(d_2\beta)\\
&=& \left(p_1\circ\Phi-p_0\circ\Phi\right)\omega-d_2\int_0^1\!\!dt\beta\\
&=& \phi_1(\omega)-\phi_0(\omega)-d_2\int_0^1\!\!dt\beta\\
&=& \phi_1(\omega)-\phi_0(\omega)-d_2h\omega
\end{eqnarray*}
(the interchange of $d_2$ and the integral is valid by the naturality of the integral), and we are done.
\end{proof}

\begin{cor}\label{cor:hoinv}
The canonical map
\[
j:\A\To\A\odot\wom_1
\]
is a homology isomorphism (in other words, $\bS$ satisfies the homotopy invariance axiom).
\end{cor}
\begin{proof}
We claim that the evaluation at $0$,
\[
p_0:\A\odot\wom_1\To\A
\]
exhibits $\A$ as a deformation retract of $\A\odot\wom_1$. Clearly, $p_0\circ j=\id_\A$. To produce the requisite homotopy from $j\circ p_0$ to $\id_{\A\odot\wom_1}$, consider the comultiplication map
\[
\delta:\KK\{u\}\To\KK\{s,t\},\quad u\mapsto st.
\]
Then
\[
\id_\A\odot\wom(\delta)
\]
is the desired homotopy. We conclude by invoking Lemma \ref{lem:hoinv}.
\end{proof}

\begin{cor}\label{cor:poincare} (\emph{The Poincar\'e Lemma}). For any differential $\bS$-algebra $\A$, the canonical map
\[
j:\A\To\A\odot\wom_n=:\wom_n(\A)
\]
is a homology isomorphism for all $n\geq0$, with homotopy inverse given by evaluation at $0$:
\[
p_0:\wom_n(\A)\To\A.
\]
\end{cor}

\begin{rem}\label{rem:grhoinv}
The \emph{graded} homotopy invariance and Poincar\'e lemma for differential graded $\bS$-algebras can be deduced by an analogous argument, working with $\circledast$ instead of $\odot$ and $\om$ instead of $\wom$.
\end{rem} 
\section{The model structures.}\label{sec:models}
\subsection{Cochain complexes.}
Let $\dg\ve$ denote the category of ($\ZZ$-graded) cochain complexes of $\KK$-modules. Since $\KK$ is in general a superalgebra, modules over it are also $\ZZ_2$-graded by default, so objects of $\dg\ve$ are $\ZZ\times\ZZ_2$-graded, and we require that the differential both increase the degree by $1$ and be Grassmann-odd. Let us also introduce the self-explanatory notation $\dg\ve_{\leq0}$, $\dg\ve_{\geq0}$ and, for $\KK$-modules equipped with an odd differential \emph{without} any integer grading, $\dd\ve$. When $\KK$ is purely even, we also consider the subcategory $\dg\ve_\even$ of those complexes for which the $\ZZ$-grading is compatible with parity (and also its positively and negatively graded versions $\dg\ve_{\geq0,\even}$ and $\dg\ve_{\leq0,\even}$); these subcategories are equivalent to the usual categories $\Ch$, $\Ch_{\geq0}$ and $\Ch_{\leq0}$ of (co)chain complexes of $\KK$-modules.

For all of the above categories, we have the standard cohomology functors
\[
\dd\ve\To\ve,\quad\dg\ve\To\gr\ve,
\]
and so on.

The category $\dg\ve$ is equivalent to the category of \emph{linear} $\EE^\times$-actions ($\E^\times$-coactions), via the left adjoint $\Sym$ (the free $\bscom_\KK$-algebra functor) of the adjunction
\begin{equation}\label{eq:usym}
\Adj{U}{\dg\ve}{\dg\bscom_\KK\Alg\simeq\bscom_\KK\Alg^{\E^\times}}{\Sym},
\end{equation}
(where $U$ denotes the underlying $\KK$-module functor), and similarly for the other versions of categories of chain complexes. The forgetful functor $U$ obviously commutes with cohomology; by the K\"{u}nneth formula, so does $\Sym$.

Each of the above categories of complexes admits a \emph{projective} Quillen model structure, with surjective maps as fibrations, and cohomology isomorphisms (quasi-isomorphisms) as the weak equivalences (with the exception of the $\ZZ_{\leq0}$-graded cases, for which the fibrations are required to be surjective only in strictly negative degrees). For $\Ch$, $\Ch_{\geq0}$ and $\Ch_{\leq0}$ these results are classical; for the other categories, exactly the same proofs apply. All these model categories are cofibrantly generated.

\begin{eg}\label{eg:supercofs}
The generating sets of cofibrations and of acyclic cofibrations for $\dd\ve$ are very small: there are two of each. Define the cells
\[
D_\even=\KK\oplus\Pi\KK=t\KK\oplus\theta\KK,\quad d=\theta\frac{\del}{\del t},
\]
and
\[
D_\odd=\Pi\KK\oplus\KK=\theta\KK\oplus t\KK,\quad d=t\frac{\del}{\del\theta}.
\]
Then define the spheres
\[
S_\even=\KK=t\KK,\quad S_\odd=\Pi\KK=\theta\KK,
\]
with trivial differentials. The set of generating acyclic cofibrations is $$J=\{j_\even,j_\odd\},$$ where
\[
j_\even:0\To D_\even,\quad j_\odd:0\To D_\odd
\]
are the initial object inclusions; the set of generating cofibrations is $$I=\{i_\even,i_\odd\},$$ where
\[
i_\even:S_\odd=\theta\KK\to t\KK\oplus\theta\KK= D_\even,\quad i_\odd: S_\even=t\KK\to\theta\KK\oplus t\KK=D_\odd
\]
are the boundary inclusions.
\end{eg}

\begin{eg}\label{eg:supergrcofs}
The generating (acyclic) cofibrations for $\dg\ve$ are a combination of those for $\Ch$ (cf. eg. \cite{haha}) with those of the previous example. For each $(n,\eps)\in\ZZ\times\ZZ_2,$ we introduce generators $t^n_\eps$ of degree $n$ and parity $\eps$, and define the cells
\[
D^n_\eps=t^n_\eps\KK\oplus t^{n+1}_{\eps+1}\KK,\quad d=t^{n+1}_{\eps+1}\frac{\del}{\del t^n_\eps}
\]
and the spheres
\[
S^n_\eps=t^n_\eps\KK, \quad d=0.
\]
Then define the generating acyclic cofibrations to be the initial object inclusions
\[
j^n_\eps:0\To D^n_\eps,
\]
and the generating cofibrations -- the boundary inclusions
\[
i^n_\eps:S^{n+1}_{\eps+1}\To D^n_\eps.
\]
\end{eg}

\subsection{The general transfer theorem.}
Our model structures on the various categories of differential graded $\bS$-algebras will be induced from those on the corresponding categories of cochain complexes by transfer along the appropriate ``free-forgetful'' adjunctions. The original argument is due to Quillen:

\begin{thm}\label{thm:transfer}\cite{ha,Crans}
Let
\[
\Adj{R}{\bC}{\bD}{L}
\]
be an adjunction, and suppose that $\bC$ is a cofibrantly generated model category, with sets $I$ and $J$ of generating cofibrations and acyclic cofibrations, respectively. Define a morphism $f$ in $\bD$ to be a fibration (resp. weak equivalence) if $Rf$ is one in $\bC$; define $f$ to be a cofibration if it has the left lifting property with respect to all acyclic fibrations. Suppose further than
\begin{enumerate}
\item the right adjoint $R$ preserves sequential colimits; and
\item every map in $\bD$ with the left lifting property with respect to all fibrations is a weak equivalence.
\end{enumerate}
Then $\bD$ becomes a cofibrantly generated model category, with the sets $L(I)$ and $L(J)$ forming generating sets of cofibrations and acyclic cofibrations, respectively. The adjunction $(L\dashv R)$ becomes a Quillen adjunction.
\end{thm}

A useful criterion for checking hypothesis $(2)$ of Theorem \ref{thm:transfer} is

\begin{thm}\textbf{Quillen's path object argument} \label{thm:path}\cite{ha,rezk,ssmod}
Let the category $\bD$ and the classes of fibrations, cofibrations and weak equivalences be as above. Suppose the following conditions hold:
\begin{enumerate}
\item $\bD$ admits a fibrant replacement endofunctor; and
\item every object has a natural path object. In other words, for every $\D\in\bD$ there exists a $P(\D)$ such that the diagonal map $\D\to\D\times\D$ factors as
    \[
    \D\mathrel{\mathop{\To}^i}P(\D)\mathrel{\mathop{\To}^q}\D\times\D,
    \]
    where $i$ is a weak equivalence and $q$ is a fibration.
\end{enumerate}
Then condition $(2)$ of Theorem \ref{thm:transfer} holds.
\end{thm}

We shall presently apply Theorems \ref{thm:transfer} and \ref{thm:path} to obtain model structures on all our categories of differential (graded) algebras by transferring the projective model structures on the corresponding versions of cochain complexes along the appropriate adjunctions. The proofs in all cases are nearly identical. The condition $(1)$ of Theorem \ref{thm:transfer} holds in each case because the right adjoints of algebraic morphisms preserve sifted (hence, in particular, sequential) colimits. All that remains is to pick the appropriate co-interval object in each case by invoking the respective version of homotopy invariance of cohomology.

Notice in particular that each transfer goes through the appropriate category of differential (graded) supercommutative algebras but the general proof does not depend on this intermediate case. In fact, we get the model structure on differential graded supercommutative $\KK$-algebras as a special case of $\bS$-algebras, for $\bS=\bscom_\KK$. This is a very slight generalization of the classical case of differential graded \emph{commutative} algebras. We highlight it here for convenience.

\subsection{The model structure for differential graded commutative superalgebras.}
From the $(\Sym\dashv U)$ adjunction \eqref{eq:usym}, we obtain

\begin{thm}\label{thm:dgKmod}
There is a cofibrantly generated model structure on each of the categories $\dd\bscom_\KK\Alg$, $\dg\bscom_\KK\Alg$, $\dg\bscom_\KK\Alg_{\geq0}$ and $\dg\bscom_\KK\Alg_{\leq0}$ obtained from the projective model structure on, respectively, $\dd\ve$, $\dg\ve$, $\dg\ve_{\geq0}$ and $\dg\ve_{\leq0}$ by transfer along the $(\Sym\dashv U)$ adjunction \eqref{eq:usym}.
\end{thm}

\begin{rem}
For $\KK$ purely even, the model structures on $\dg\bcom_\KK\Alg_{\leq0}$, \newline  $\dg\bcom_\KK\Alg_{\geq0}$ and $\dg\bcom_\KK\Alg$ are obtained by transfer from $\Ch_{\leq0}$ $\Ch_{\geq0}$ and $\Ch$, respectively. These model structures were derived in, respectively and in historical sequence, \cite{rht}, \cite{pldr} and \cite{haha}. The same proofs apply in the corresponding super versions (with the non-graded version $\dd\bscom_\KK\Alg$ analogous to the unbounded $\ZZ$-graded version of \cite{haha}).
\end{rem}

Let us sketch our proof; the argument is different from the more explicit ones in the sources cited above, and will apply, with very small modifications, to algebras over other super Fermat theories (see below).

\begin{proof}
The condition $(1)$ of Theorem \ref{thm:transfer} holds because the right adjoints of algebraic morphisms preserve sifted (hence, in particular, sequential) colimits. To derive condition $(2),$ we apply Theorem \ref{thm:path}. Since fibrations are just the surjective maps, every object is fibrant, so we can take the identity endofunctor to satisfy condition $(1)$ of the theorem. To check condition $(2),$ we observe that the algebra
\[
\Omega_1=\KK[t,dt]
\]
is a cofibrant co-interval object in $\dg\bscom_\KK\Alg$. The fact that the canonical inclusion
\[
\A\To\A[t,dt]=\A\otimes_\KK\Omega_1
\]
is a weak equivalence (in fact, an acyclic cofibration) is the homotopy invariance of de Rham cohomology (algebraic version) (Proposition \ref{prop:hialg}).
\end{proof}

\begin{eg}\label{eg:supergralgcofs}
The sets of generating (acyclic) cofibrations for $\dg\bscom_\KK\Alg$ are obtained from those for $\dg\ve$ (cf. Example \ref{eg:supergrcofs}) by applying the free algebra functor. To wit, they are
\[
j^n_\eps:\KK\To\D^n_\eps,\quad i^n_\eps:\S^{n+1}_{\eps+1}\To\D^n_\eps,
\]
where
\[
\D^n_\eps=\KK[t^n_\eps,t^{n+1}_{\eps+1}],\quad d=t^{n+1}_{\eps+1}\frac{\del}{\del t^n_\eps}
\]
and
\[
\S^n_\eps=\KK[t^n_\eps],\quad d=0.
\]
\end{eg}

%\medskip

%For the remainder of the paper, we assume that the theory $\bS$ admits integration; in particular, this includes all of the examples listed in Remark %\ref{rem:egthwithint}.

\subsection{The model structure for differential $\bS$-algebras.}
The ungraded case is the most clear-cut. We consider the adjunction
\begin{equation}\label{eq:adjmod}
\Adj{U}{\dd\ve}{\dd\bS\Alg}{\bS(\,)}
\end{equation}
where $U$ is the forgetful functor and $\bS(V)$ is the free $\bS$-algebra on a $\KK$-module $V$. We endow $\dd\ve$ with the projective model structure. Let us first dispense with the showing that $U$ preserves sequential colimits. First notice that in general, arbitrary colimits of coactions are computed by forgetting the coaction (Remark \ref{rem:cofreecoact}). Furthermore, in the adjunction (\ref{eq:adjmod}), each category is an algebraic category, so sifted colimits (and in particular sequential colimits) are computed pointwise, and since $U$ is induced by an algebraic morphism between algebraic categories, it \emph{preserves} sifted colimits. With this aside, to apply the transfer theorem \ref{thm:transfer}, we need only to show the conditions of Theorem \ref{thm:path} are satisfied. Just as in the proof of Theorem \ref{thm:dgKmod}, role of a fibrant replacement functor can be played by the identity functor. Finally, notice that $\wom_1$ can be taken as the co-interval object, so the path object functor will be
\[
P(\A)=\A\odot\wom_1=\A\{t,dt\}.
\]
Indeed, we have a factorization of the diagonal
\[
\A\mathrel{\mathop{\To}^j}\A\{t,dt\}\mathrel{\mathop{\To}^q}\A\times\A,
\]
where $j$ is the canonical inclusion and $q=(p_0,p_1)$. Notice that $q$ is surjective since for all $a_0,a_1\in\A$,
\[
q(a_0(1-t)+a_1t)=(a_0,a_1),
\]
while $j$ is a homology isomorphism by Corollary \ref{cor:hoinv}. Therefore, Theorem \ref{thm:path}, and hence Theorem \ref{thm:transfer} applies, and we have

\begin{thm}\label{thm:modelsuper}
The category $\dd\bS\Alg$ of differential $\bS$-algebras admits a cofibrantly generated Quillen model structure with surjective maps as fibrations and homology isomorphisms as weak equivalences.
\end{thm}

\begin{eg}\label{eg:supersalgcofs}
The generating (acyclic) cofibrations for $\dd\bS\Alg$ are obtained from those for $\dd\ve$ (cf. Example \ref{eg:supercofs}) by applying the free $\bS$-algebra functor. To wit, they are
\[
j_\even:\KK\To\hat\D_\even,\quad j_\odd:\KK\To\hat\D_\odd
\]
and
\[
i_\even:\hat\S_\odd\To\hat\D_\even,\quad i_\odd:\hat\S_\even\To\hat\D_\odd,
\]
where
\[
\hat\D_\even=\KK\{t,\theta\},\quad d=\theta\frac{\del}{\del t},
\]
\[
\hat\D_\odd=\KK\{\theta,t\},\quad d=t\frac{\del}{\del \theta}
\]
and
\[
\hat\S_\even=\KK\{t\},\quad\hat\S_\odd=\KK\{\theta\},\quad d=0
\]
for both.
\end{eg}

\subsection{The model structure for differential graded $\bS$-algebras.}
%Our next goal is to transfer these structures further to $\bS\Alg^{\hat\E^\times}_\alg$, and finally to $\bS\Alg^{\hat\E^\times}$, along the adjunctions %$(\tau^{\E^\times}_\circ\dashv\tau_{\E^\times}^\circ)$ and $(\tau^{\E^\times}_\dagger\dashv\tau_{\E^\times}^\dagger)$ \eqref{eq:adjdg}, and likewise for the %bounded, ungraded and even versions.

%Now we are going to transfer the models structure on $\dg\bscom_\KK\Alg$ (resp. $\dg\bcom_\KK\Alg$) up to $\dg\bS\Alg$ (resp. $\bS\Alg^{\E^\times}_{\even,\alg}$) %along the adjunction
%\[
%\Adj{\tau_{\E^\times}^\circ}{\dg\bscom_\KK\Alg}{\dg\bS\Alg}{\tau^{\E^\times}_\circ}
%\]
%(resp. the corresponding even version).
Here, the difference from the ungraded case of $\dd\bS\Alg$ and the ``completely algebraic'' case of $\dg\bscom_\KK\Alg$ is that we use $\wom_1^\alg$ -- which is nothing but $\om_1$ with the canonical (initial!) $\bS$-algebra structure on $\om_1^{0,\cl}=\KK$ -- and the coproduct $\circledast$ instead of $\otimes$ or $\odot$, as in the graded version of homotopy invariance (Remark \ref{rem:grhoinv}). The same argument goes through, and we obtain

\begin{thm}\label{thm:modelsalg}
There is a cofibrantly generated Quillen model structure on the category $\dg\bS\Alg$ with surjective maps as fibrations and cohomology isomorphisms as the weak equivalences.
%Likewise for $\bS\Alg^{\E^\times}_{\even,\alg}$.
\end{thm}

\begin{eg}\label{eg:dgsalgcofs}
The generating (acyclic) cofibrations for this model structure are the same as in Example \ref{eg:supergralgcofs}, except we use
\[
\hat\D^{-1}_\odd=\KK[t^{-1}_\odd]\{t^0_\even\}
\]
and
\[
\hat\S^0_\even=\KK\{t^0_\even\}
\]
instead of $\D^{-1}_\even$ and $\S^0_\even$.

%The generating sets for $\bS\Alg^{\E^\times}_{\even,\alg}$ form subsets of these corresponding to those generators $t^n_\eps$ for which $\eps=n$ modulo $2$.
\end{eg}

\begin{rem}
For $\bS=\bSE$, the theorem specializes to yield a model structure on $\dg\bE\Alg$, with generating cofibrations a subset of those for $\dg\bSE\Alg$ corresponding to those generators $t^n_\eps$ for which $\eps=n$ modulo $2$.
\end{rem}

\begin{rem}
The model structures on $\dg\bS\Alg_{\leq0}$ and $\dg\bS\Alg_{\geq0}$ can be either obtained directly by an analogous argument, or induced from $\dg\bS\Alg$ along the inclusion/truncation adjunction.
\end{rem}

\subsection{The model structure on $\bS\Alg^{\hat\E^\times}$.}
Here we must be a little careful, as there are \emph{two} ways to induce a model structure. One is by forgetting the $\hmG$-coaction and transferring from $\ZZ_2$-graded complexes, via $\dd\bS\Alg$, in which case the weak equivalences would be isomorphisms in cohomology with respect to the differential. The other is via algebraization, i.e. restricting to the algebraic characters of the multiplicative group (the integers) -- in other words, by transferring from $\ZZ_2\times\ZZ$-graded cochain complexes, via $\dg\bS\Alg$ along the adjunction
\[
\Adj{\tau_{\E^\times}^\dagger}{\dg\bS\Alg}{\bS\Alg^{\E^\times}}{\tau^{\E^\times}_\dagger}.
\]
For this model structure, the weak equivalences are isomorphisms in the \emph{algebraic part} of the cohomology (recall that cohomology commutes with algebraization by Proposition \ref{prop:coh}). The two models structures are \emph{inequivalent} due to the existence of pathological nontrivial $\hmG$-coactions with trivial algebraic part. For instance, the algebra in Remark \ref{rem:badact} (endowed with the zero differential) is equivalent to $\RR$ in the latter model structure but not in the former. This was exactly the point of introducing algebraization: to weed out such pathological coactions by forcing them to be weakly equivalent to algebraic ones. Of course, it is the latter -- algebraic -- model structure that is relevant in derived geometry.

To deduce this model structure, we make two observations. The first one is that colimits in $\bS\Alg^{\hat\E^\times}$ are computed by forgetting the coaction, and sifted (hence sequential) colimits of $\bS$-algebras are computed on underlying modules (or even on $\ZZ_2$-graded sets). It follows immediately that these colimits preserve the character decomposition, hence are preserved by algebraization. The second observation concerns the co-interval object, $\wom_1$. The $\E^\times$-coaction on $\wom_1$ is clearly essentially algebraic, and for any $\A\in\bS\Alg^{\hat\E^\times}$, the canonical map
\[
j:\A\To\A\odot\wom_1
\]
induces an isomorphism in \emph{all} of cohomology (Corollary \ref{cor:hoinv}), hence in particular on the algebraic part. Thus, the same argument as in the proof of Theorem \ref{thm:modelsuper} goes through, and we have

\begin{thm}\label{thm:modelss}
There is a cofibrantly generated Quillen model structure on the category $\bS\Alg^{\E^\times}$ with surjective maps as fibrations and maps inducing isomorphisms in $H_\alg$ as the weak equivalences. For $\bS=\bSE$, we also have a model structure on $\bSE\Alg^{\E^\times}_\even$.
\end{thm}

\begin{eg}\label{eg:dgscofs}
The generating (acyclic) cofibrations for $\bS\Alg^{\E^\times}$ are the same as in Example \ref{eg:supergralgcofs}, except we use
\[
\hat\D^n_\eps=\KK\{t^n_\eps,t^{n+1}_{\eps+1}\},\quad\hat\S^n_\eps=\KK\{t^n_\eps\}
\]
instead of $\D^n_\eps$ and $\S^n_\eps$. For $\bSE\Alg^{\E^\times}_\even$, we only use those generators for which the degree is compatible with parity.
\end{eg}

\begin{rem}
As for $\dg\bS\Alg$, the model structures on $\bS\Alg^{\E}$ and $\bS\Alg^{\E^\op}$ can be either constructed directly using the same argument, or induced from $\bS\Alg^{\E^\times}$ via the inclusion/truncation adjunction.
\end{rem}

\begin{rem}
In all of the model structures we have constructed, every object is fibrant.
\end{rem} 
\section{The simplicial structure.}\label{sec:simplicial}
\subsection{Right and left almost simplicial model categories.}

Recall that a \emph{simplicial category} is a category $\bC$ enriched in the category $\SS=\Set^{\bDelta^\op}$ of simplicial sets. We have the simplicial mapping space functor
\[
\IHom(-,-):\bC^\op\times\bC\To\SS
\]
with an associative and unital composition, the unit being the terminal simplicial set. Moreover, $\IHom(\C,\D)_0=\bC(\C,\D)$.

\begin{defn}
We say that a simplicial category $\bC$ \emph{admits (finite) tensors} if there is a functor
\[
-\otimes-:\SS_{(f)}\times\bC\To\bC,\quad (K,\C)\mapsto K\otimes\C
\]
and a natural isomorphism of functors
\[
\SS(-,\IHom(-,-))\simeq\bC(-\otimes-,-):\SS_{(f)}^\op\times\bC^\op\times\bC\To\Set.
\]
Dually, we say that $\bC$ \emph{admits (finite) cotensors} if there is a functor
\[
(-)^{(-)}:\bC\times\SS_{(f)}^\op\To\bC,\quad (\C,K)\mapsto\C^K
\]
and a natural isomorphism of functors
\[
\SS(-,\IHom(-,-))\simeq\bC(-,(-)^{(-)}):\SS_{(f)}^\op\times\bC^\op\times\bC\To\Set.
\]
\end{defn}

\begin{rem}
We do not make the usual requirement that $(\C^K)^L\simeq\C^{K\times L}$; in fact, this is false for all the categories of differential graded algebras we shall consider.
\end{rem}

We state the following lemma for the cotensor case only; the dual tensor case can be stated and proved by going to the opposite category.

\begin{lem}\label{lem:cotensor}
Let $\bC$ be a category with coproducts and (finite) limits, $\W_\bullet$ a simplicial object in $\bC$. Then the assignment
\begin{equation}\label{eqn:enrich}
\IHom(\A,\B)_n=\bC(\A,\B\amalg\W_n)
\end{equation}
defines a simplicial enrichment of $\bC$. Furthermore, it admits (finite) cotensors given by the formula
\begin{equation}\label{eqn:cotensor}
\A^K=\mathrel{\mathop{{\ulim}}_{\Delta[n]\to K}}\!(\A\amalg\W_n)
\end{equation}
where
\[
K=\mathrel{\mathop{{\ucolim}}_{\Delta[n]\to K}}\!\Delta[n]
\]
is expressed as the canonical colimit of its (non-degenerate) simplices.
\end{lem}
\begin{proof}
Since the assignment \eqref{eqn:enrich} is obviously functorial with respect to maps in $\bDelta^\op$, $\IHom(\A,\B)_\bullet$ is a simplicial set for each $\A,\B\in\bC$. To define composition, given $f\in\IHom(\A,\B)$, $g\in\IHom(\B,\C)$, we let $g\circ f\in\IHom(\A,\C)$ be the following composite:
\[
\A\mathrel{\mathop{\To}^f}\B\amalg\W_n\mathrel{\mathop{\longlongrightarrow}^{g\amalg\id}}\C\amalg\W_n\amalg\W_n\mathrel{\mathop{\longlongrightarrow}^{\id\amalg\nabla}}\C\amalg\W_n,
\]
where $\nabla:\W_n\amalg\W_n\to\W_n$ is the codiagonal, and we have omitted the obvious associativity isomorphisms for the coproduct. It is obvious that this composition is compatible with the simplicial structure, associative and unital with the unit
\[
\id_\A:*\To\IHom(\A,\A)
\]
given by the canonical inclusion into the coproduct:
\[
i_n:\A\To\A\amalg\W_n.
\]
To see that this simplicial structure admits cotensors given by the formula \eqref{eqn:cotensor}, we pick a simplicial set $K$ and compute:
\begin{eqnarray*}
\SS(K,\IHom(\A,\B))&=&\SS\left(\mathrel{\mathop{{\ucolim}}_{\Delta[n]\to K}}\!\Delta[n],\IHom(\A,\B)\right)\\
&=&\mathrel{\mathop{{\ulim}}_{\Delta[n]\to K}}\!\SS(\Delta[n],\IHom(\A,\B))\\
&=& \mathrel{\mathop{{\ulim}}_{\Delta[n]\to K}}\!\IHom(\A,\B)_n\\
&=& \mathrel{\mathop{{\ulim}}_{\Delta[n]\to K}}\!\bC(\A,\B\amalg\W_n)\\
&=& \bC\left(\A,\mathrel{\mathop{{\ulim}}_{\Delta[n]\to K}}\!\B\amalg\W_n\right).
\end{eqnarray*}
Notice in particular that
\[
\B^{\Delta[n]}=\B\amalg\W_n,
\]
since the category of simplices of $\Delta[n]$ contains a terminal object, namely $$\id:\Delta[n]\to\Delta[n].$$
\end{proof}

\begin{defn}
Suppose $\bC$ is a Quillen model category. Consider $\SS$ with the classical (Kan-Quillen) model structure. $\bC$ is said to \emph{satisfy the corner axiom} (axiom (SM7) of \cite{ha}) if, for every cofibration $i:\A\to\B$ and a fibration $p:\Y\to\X$ in $\bC$, the natural map
\[
\IHom(\B,\Y)\mathrel{\mathop{\longlonglongrightarrow}^{(i^*,p_*)}}\IHom(\A,\Y)\times_{\IHom(\A,\X)}\IHom(\B,\X)
\]
is a Kan fibration of simplicial sets, which is furthermore acyclic if either $i$ or $p$ is.
\end{defn}

\begin{defn}
A simplicial category $\bC$ which is also a model category is said to be a \emph{left almost simplicial model category} if it satisfies the corner axiom and in addition
\begin{description}
\item[(i)] $\bC$ admits finite tensors, and
\item[(ii)] tensoring preserves the initial object $\varnothing\in\bC$, i.e. $K\otimes\varnothing\cong\varnothing$ for any finite simplicial set $K$.
\end{description}
Dually, $\bC$ is said to be a right almost simplicial model category if it satisfies the corner axiom and in addition
\begin{description}
\item[(i')] $\bC$ admits finite cotensors, and
\item[(ii')] cotensoring preserves the terminal object $\star\in\bC$, i.e. $\star^K \cong \star$ for any finite simplicial set $K$.
\end{description}
\end{defn}

\begin{rem}
$\bC$ is left almost simplicial if and only if $\bC^\op$ is right almost simplicial, and vice versa. If $\bC$ admits finite cotensors, then for any $\C\in\bC$ we have
\[
\IHom(\varnothing,\C)_n=\SS(\Delta[n],\IHom(\varnothing,\C))=\bC(\varnothing,\C^{\Delta[n]})=*
\]
for every $n$, so $\varnothing$ is ``simplicially initial''; however to conclude that $\star$ is ``simplicially terminal'', in the sense that $\IHom(\C,\star)=*$ for every $\C$, we need condition $(ii')$. Similarly for condition $(ii)$ with the roles of tensors and cotensors, and $\star$ and $\varnothing$ reversed. Of course, when $\bC$ admits \emph{both} finite tensors \emph{and} cotensors, conditions $(ii)$ and $(ii')$ are redundant.
\end{rem}

\begin{rem}
By the above remark, if $\bC$ is either right or left simplicial, $\B\in\bC$ is cofibrant, and $\Y\in\bC$ is fibrant, then $\IHom(\B,\Y)$ is a Kan complex. To see this, just use the corner axiom with $\A=\varnothing$, $\X=*$.
\end{rem}

\begin{rem}\label{rem:corner}
It is a standard observation that, if $\bC$ admits finite cotensors, the corner axiom is equivalent to the statement that, for every cofibration
\[
j:K\To L
\]
of finite simplicial sets and a fibration
\[
p:\Y\To\X
\]
in $\bC$, the map
\[
\Y^L\mathrel{\mathop{\longlonglongrightarrow}^{(\Y^j,p^L)}}\Y^K\times_{\X^K}\X^L
\]
is a fibration in $\bC$ which furthermore is acyclic if either $j$ or $p$ is. In fact, it suffices to restrict $j$ to lie in the generating set of (acyclic) cofibrations in $\SS$, i.e. horn and boundary inclusions.
\end{rem}

\subsection{Simplicial structure of differential forms.}

We shall now proceed to equip each of our model categories of differential graded $\bS$-algebras with a \emph{right} almost simplicial structure. It is a variation of the argument in (\cite{pldr}, Section $5$). We will use the simplicial object formed by the algebras of differential forms, in exactly the same way as done in \cite{pldr}. We consider the algebras $\wom_n$, $n\in\NN$, but with a slightly different presentation, namely
\[
\wom_n=\KK\{t^0,\ldots,t^n,dt^0,\ldots,dt^n\}/\left(T_n,dT_n\right),
\]
where
\[
T_n=1-\sum_{i=0}^nt^i,
\]
making the simplicial structure apparent: given an $\phi:[m]\to[n]$ in $\bDelta$, the formula
\[
\phi^*t^k=\sum_{\phi(j)=k}t^j,\quad0\leq k\leq n
\]
defines a morphism $\phi^*:\wom_n\to\wom_m$ respecting the differentials and gradings. Given an algebra $\B\in\bS\Alg^{\E^\times}$ (or $\dd\bS\Alg$), define also the simplicial object
\[
\wom_\bullet(\B)=\B\odot\wom_\bullet
\]
in the respective category (``differential forms on simplices with coefficients in $\B$''); it has two differentials (and gradings, when appropriate) -- one internal, coming from $\B$, the other the standard one on differential forms. In case $\bS=\bcom$ or $\bscom$, we write $\om_\bullet$ instead of $\wom_\bullet$. To work in $\dg\bS\Alg=\bS\Alg^{\E^\times}_\alg$, we use $\wom^\alg_\bullet$, which is just $\om_\bullet$, except we remember the $\bS$-algebra structure on its subalgebra of $0$-cocycles, namely, $\KK$.

Finally, given a simplicial set $K\in\SS$, define
\[
\B^K=\SS(K,\wom_\bullet(\B)).
\]
It is clear that this definition coincides with \eqref{eqn:cotensor} (with $\W_\bullet=\wom_\bullet$) and gives rise to a simplicial enrichment with cotensors of each of our categories of differential (graded) $\bS$-algebras via Lemma \ref{lem:cotensor}. The following two results are key:

\begin{lem}\label{lem:contract}
For each $(\epsilon,n)\in\ZZ_2\times\ZZ$, $\wom_\bullet^{(\epsilon,n)}(\B)$ is a contractible Kan complex for any $\B$.
\end{lem}
\begin{proof}
For $\bS=\bcom_\KK$ or $\bcinf$ this is classical (see eg. \cite{pldr} or \cite{infitop}). However, a brief inspection of the short proof given in (\cite{linftynilp}, Lemma 3.2) makes clear that it works for any $\bS$ and any $\B$, even without assuming integration.
\end{proof}

The second result is an explicit form of the de Rham theorem for $\B^K$. Again, for $\bS=\bcinf$ and $\B=\RR$ this is classical and is due to Dupont (cf. \cite{dupont76,dupont78}); a different proof (not using an explicit contraction) was given for $\bS=\bcom$ and $\B=\QQ$ in \cite{pldr}. However, by closely examining the argument given in \cite{linftynilp}, we again discover that it remains valid for any super Fermat theory $\bS$ with integration and any $\B$. We shall briefly sketch the argument here, referring the reader to \cite{linftynilp} or \cite{dupont78} for details.

Recall first the following construction, due in its classical form to Whitney \cite{geomint}. Given an ordered $k+1$-tuple $(i_0,\ldots,i_k)$ of elements of the set $\{0,\ldots,n\}$, let
\[
I_{i_0\ldots i_k}:\wom_n\To\KK
\]
be the integral over the $k$-simplex spanned by the vertices $(e_{i_0},\ldots,e_{i_k})$ in $\KK^{n+1}$. By Lemma \ref{lem:basextint}, this extends to a map of $\B$-modules
\[
I_{i_0\ldots i_k}:\wom_n(\B)\To\B.
\]
Now consider the elementary forms (also known as the \emph{Whitney forms})
\[
\omega_{i_0\ldots i_k}=k!\sum_{q=0}^k(-1)^qt_{i_q}dt^{i_0}\cdots\widehat{dt^{i_q}}\cdots dt^{i_k}.
\]
Notice that these span a subcomplex $C_n$ of $\wom_n$ since
\[
d\omega_{i_0\ldots i_k}=\sum_{i=0}^n\omega_{ii_0\ldots i_k}.
\]
Denote by $C_n(\B)$ the $\B$-submodule of $\wom_n(\B)$ spanned by the same forms; since the Whitney forms are \emph{linear} in the $t^i$'s, it follows that
$$C_n(\B)=C_n\otimes_\KK\B.$$
The inclusions $C_\bullet(\B)\hookrightarrow\wom_\bullet(\B)$ are clearly compatible with the simplicial structure, and with both differentials. The complex $C_\bullet(\B)$ is isomorphic to the complex of simplicial cochains on $\Delta[n]$ with coefficients in $\B$. Furthermore, there is a projection $P_\bullet:\wom_\bullet\to C_\bullet$ (due to Whitney \cite{geomint}) splitting the inclusion and given by the formula
\[
P_n\omega=\sum_{k=0}^n\sum_{i_0<\cdots<i_k}\omega_{i_0\ldots i_k}I_{i_0\ldots i_k}(\omega).
\]
This projection is also extendable to $P_\bullet:\wom_\bullet(\B)\to C_\bullet(\B)$ by $\B$-linearity and is compatible with both the simplicial structure and the two differentials. Lastly, for any simplicial set $K$ we have an induced inclusion of subcomplexes
\[
C(K,\B)=\SS(K,C_\bullet(\B))\subset\B^K,
\]
which is split by the induced projection
\[
P:B^K\To C(K,\B).
\]
Dupont \cite{dupont76,dupont78} found a contracting homotopy for $P_\bullet$ (for $\bS=\cinf$), thus giving an explicit proof of de Rham's theorem. Again, an inspection of his formulas (\cite{dupont78}, Chapter $1$, or \cite{linftynilp}, Section $3$) reveals that they remain valid for any $\bS$ with integration, and any $\B$. We have
\begin{thm}\label{thm:derham} \emph{(The de Rham theorem.)}
There is a simplicial endomorphism $s_\bullet:\wom_\bullet\to\wom_\bullet[\odd,-1]$ (i.e. $s_\bullet$ is odd and decreases the degree by $1$) such that
\[
\id-P_\bullet=[d,s_\bullet].
\]
Furthermore, $s_\bullet^2=0$. Consequently, for every simplicial set $K$ and every $\B$, there is an induced contracting homotopy $s:\B^K\to\B^K[\odd,-1]$ such that
\[
\id-P=[d,s],
\]
retracting $\B^K$ onto $C(K,\B)$.
\end{thm}

\begin{thm}
The construction of Lemma \ref{lem:cotensor} with $\W_\bullet=\wom_\bullet$ equips the model categories $\dg\bS\Alg$, $\dd\bS\Alg$ and $\bS\Alg^{\E^\times}$ with right almost simplicial structures.
\end{thm}
\begin{proof}
That $0^K=0$ is obvious, since the terminal algebra $0$ is a zero with respect to coproducts in all the categories of algebras considered. To prove the corner axiom, we follow the argument of (\cite{pldr}, Proposition 5.3) closely, with the exception that we will not take advantage of the exactness of coproducts, which generally fails to hold. The arguments for all the mentioned categories are analogous, so we will only do one -- say, $\dd\bS\Alg$. By Remark \ref{rem:corner}, it suffices to show that, given a surjective map $p:\Y\to\X$ of differential $\bS$-algebras and an injective map $j:K\to L$ of simplicial sets, the map
\[
\Y^L\mathrel{\mathop{\longlonglongrightarrow}^{(\Y^j,p^L)}}\Y^K\times_{\X^K}\X^L
\]
is surjective, and in addition, is a homology isomorphism if either $p$ is a homology isomorphism, or if $j$ is acyclic. In fact, for the first part it suffices to take $j$ to be a boundary inclusion $i_n:\del\Delta[n]\to\Delta[n]$. For surjectivity, it suffices to show that both
$$\Y^{i_n}:\Y^{\Delta[n]}\To\Y^{\del\Delta[n]}$$
and
$$p^{\Delta[n]}:\Y^{\Delta[n]}\To\X^{\Delta[n]}$$
are surjective for all $n$. Now, by Lemma \ref{lem:contract}, $\wom_\bullet(\Y)$ is a contractible Kan complex, which means precisely that $Y^{i_n}$ is surjective for all $n$. As for $p^{\Delta[n]}$, it is just
\[
p\odot\id:\Y\odot\wom_n\To\X\odot\wom_n;
\]
since surjective maps are the regular epimorphisms, they are preserved by any change of base; hence, $p\odot\id$ is surjective if $p$ is.

Now, suppose $p$ is also a homology isomorphism. It suffices to show that $$p^{\Delta[n]}=p\odot\id_{\wom_n}$$ is a homology isomorphism; the conclusion will then follow by Mayer-Vietoris. Notice that we have a decomposition
\[
p=p_0\circ(p\odot\id_{\wom_n})\circ i,
\]
where $i:\Y\to\wom_n(\Y)$ is the canonical inclusion and $p_0:\wom_n(\X)\to\X$ is the evaluation at $0$. Since $p$ is a homology isomorphism by assumption, and $j$ and $p_0$ are by the Poincar\'e Lemma \ref{cor:poincare}, it follows that $p\odot\id_{\wom_n}$ is also a homology isomorphism, by the ``2 out of 3'' property.

Lastly, suppose that $j:K\to L$ is an acyclic cofibration; it suffices to take $j=j_{n,k}:\Lambda^{n,k}\to\Delta[n]$ to be a horn inclusion. Now it suffices to show that
\[
\Y^{j_{n,k}}:\wom_n(\Y)\to\Y^{\Lambda^{n,k}}
\]
is a homology isomorphism (and then conclude by Mayer-Vietoris). Applying the de Rham theorem \ref{thm:derham}, we reduce it to showing that
\[
C_n(\Y)=C(\Delta[n])\otimes\Y\To C(\Lambda^{n,k})\otimes\Y
\]
is a homology isomorphism; however, the normalized simplicial cochain complexes $C(\Delta[n])$ and $C(\Lambda^{n,k})$ are contractible, so the spectral sequences collapse right away, implying the result.
\end{proof}

\begin{rem}
Suppose it was known that the ground ring $\mathbb{F}$ of $\bS$ is a field, and coproducts in $\bS$ commute with finite limits. Then we could apply the argument of \cite{pldr}, Proposition $5.3$ verbatim. Moreover, we could apply it to $\bS_\KK$ for any $\KK\in\bS\Alg$, by ``cancelation of base rings'':
\[
\wom_n(\KK)\odot_\KK\B=(\wom_n(\mathbb{F})\odot_\mathbb{F}\KK)\odot_\KK\B=\wom_n(\mathbb{F})\odot_\mathbb{F}\B.
\]
This trick was used in \cite{haha} to get a simplicial structure on the category of algebras over a linear operad over any ground ring $\KK$ containing $\mathbb{F}=\QQ$.
\end{rem}

\begin{cor}
If $\B$ is cofibrant and $\Y$ is arbitrary (since all objects are fibrant), the simplicial set $\IHom(\B,\Y)$ is a Kan complex.
\end{cor} 
\appendix
\section{Modules over Superalgebras}\label{app:modules}
%(In this subsection: define modules a la Quillen, derivations as sections. Also derivations as infinitesimal automorphisms, show they are the same.)

In this appendix, we review module theory for supercommutative algebras. To this end, it is convenient to first introduce the symmetric monoidal $\KK$-linear abelian category of super $\KK$-modules, $\sv$, with $\KK$ a fixed commutative ring. Its objects are functors $$\ZZ_2 \to \ve,$$ where the former should be regarded as a discrete \emph{$\ve$-enriched} category: the endomorphisms form the unit object (i.e. $\KK$, with composition given by multiplication), and the other morphism objects are null.
\begin{rem}
$\ZZ_2$ viewed in this way becomes a symmetric monoidal $\KK$-linear category, with tensor product given on objects by \emph{addition} in $\ZZ_2$, on morphisms by \emph{multiplication} in $\KK$, and the braiding given by the Koszul sign.
\end{rem}

The category $\sv$ is canonically enriched in $\ve.$ The enrichment is determined by $$\Hom\left(V,W\right):=\IHom\left(V_\even,W_\even\right) \times \IHom\left(V_\odd,W_\odd\right),$$ where $\IHom\left(S,T\right)$ denotes the internal maps from $S$ to $T$, i.e. the $\KK$-linear space of such $\KK$-linear maps. There are two particularly special such super $\KK$-modules, namely $$\underline 0:=\left(\KK,0\right)$$ and $$\underline 1:=\left(0,\KK\right).$$ As functors, these are precisely the representable functors, and are called the even and odd line respectively.

\begin{rem}
There is an enriched version of the Yoneda lemma, since for any $\ii=0,1$ and any super $\KK$-module $V,$ one has
$$\Hom\left(\underline \ii,V\right) \cong V_i.$$
\end{rem}

\begin{rem}
Any super commutative algebra $A$ in $\comsalg_\KK$ has an underlying super $\KK$-module $\left(A_\even,A_\odd\right).$
\end{rem}

%\smallskip
%\textcolor[rgb]{0.00,0.00,1.00}{It is even better to regard $\ZZ_2$ as a discrete \emph{$\Vect$-enriched} category: the endomorphisms form the unit object (i.e. $\ZZ$, with composition given by multiplication), other morphism objects null. Then the lines really are the representable functors. Besides, viewed in this way, $\ZZ_2$ itself becomes a symmetric monoidal category, with tensor product given on objects by \emph{addition} in $\ZZ_2$, on morphisms by \emph{multiplication} in $\ZZ$, and the braiding given by the Koszul sign. The tensor structure on $\Vect^{\ZZ_2}$ is then uniquely characterized by the requirements that the Yoneda embedding be a symmetric monoidal functor, and that the tensor product commute with colimits.}
%\smallskip

The abelian group structure on $\ZZ_2$ comes into play in introducing the symmetric monoidal category structure on $\sv$. For $V$ and $W$ two super $\KK$-modules, define their tensor product by $$\left(V \otimes W\right)_i:=\bigoplus \limits_{j+k=i} V_j \otimes W_k.$$

\begin{rem}
 This makes $\underline 0$ the tensor unit.
\end{rem}

This category also comes equipped with a canonical braiding given by
\begin{eqnarray*}
\tau_{V,W}:V \otimes W &\stackrel{\sim}{\longrightarrow}& W \otimes V\\
v \otimes w &\mapsto& \left(-1\right)^{\deg\left(v\right)\cdot\deg\left(w\right)}\cdot w \otimes v,\\
\end{eqnarray*}
where $\deg\left(v\right)=i$ for all $v \in V_i$.

\begin{rem}
The tensor product and braiding on $\sv$ are uniquely characterized by the requirements that the Yoneda embedding $$\ZZ_2 \hookrightarrow \sv$$ be a symmetric monoidal functor, and that the tensor product commute with colimits.
\end{rem}

The $\KK$-linear category $\sv$ is furthermore symmetric-monoidal closed, under the internal $\IHom$ defined by:
$$\IHom_{\sv}\left(V,W\right)_i:=\prod\limits_{i'\in \ZZ_2} \IHom_{\ve}\left(V_{i'},W_{i+i'}\right).$$
Finally, the category $\sv$ comes equipped with a canonical action of the group $\ZZ_2$. To describe it, it suffices to describe the action of the non-identity element $1$ on a super $\KK$-module $V,$ and this action is simply given by parity reversal: $$\left(V \cdot 1\right)_i:=V_{i+1}=V_{-i}.$$ For $\ii=0,1$ we shall adopt the notation $$\left(V\right)\left[\ii\right]:=V\cdot \ii,$$ and in the particular case $\ii=1$ we will often use the notation $$\Pi V:=\left(V\right)\left[1\right],$$ as this is standard in the literature.

\begin{rem}
For all $\ii$,
$$\left(V\right)\left[\ii\right]\cong\IHom\left(\underline \ii,V\right)\cong V \otimes \underline \ii \cong \underline \ii \otimes V.$$
\end{rem}

\begin{defn}

A \emph{superalgebra} is an internal monoid in $\sv,$ that is an object $\A$ together with a multiplication map $$m:\A \otimes \A \to \A$$ and a unit map $$\eta:\even \to  \A$$ satisfying the usual axioms. A superalgebra is \emph{supercommutative} if in addition to the usual axioms, the following diagram commutes:
$$\xymatrix@C=1.5cm{\A \otimes \A \ar[dd]_-{\tau_{\A,\A}} \ar[rd]^-{m} &\\
& \A\\
\A \otimes \A, \ar[ru]_-{m} & }$$
where $\tau$ denotes the  braiding.
\end{defn}

\begin{rem}
This definition of supercommutative superalgebra readily agrees with the standard one.
\end{rem}

\begin{defn}
Given any superalgebra $\A$, its monoid structure induces a monoid structure on the endofunctor $$\A \otimes \left(\mspace{3mu} \cdot \mspace{3mu}\right):\sv \to \sv,$$ in other words, it makes $$V \mapsto \A \otimes V$$ into a monad. Algebras for this monad are called \emph{left $\A$-modules}. Similarly one can define \emph{right $\A$-modules} as algebras for the monad $\left(\mspace{3mu} \cdot \mspace{3mu}\right)\otimes \A$.
\end{defn}

\begin{rem}
By definition, a left $\A$-module is a super $\KK$-module $V$ together with a map $\rho:\A \otimes V \to V$ satisfying certain axioms. Let us introduce the notation $\rho\left(a,v\right)=:a \cdot v,$ for homogeneous elements $v,$ to make this more akin to modules for ordinary algebras. With this notation, the axioms of an $\A$-module are exactly the same as the usual axioms for a module of a commutative ring.
\end{rem}

\begin{defn}
Suppose that $\A$ is a supercommutative superalgebra, and $V$ is a right $\A$-module. Then $V$ has the canonical structure of a left $\A$-module by the equation $$a \cdot v:=\left(-1\right)^{\deg\left(a\right)\deg\left(v\right)}v \cdot a.$$ Moreover, this left module structure is compatible with the right module structure in the obvious way, making $V$ with these left and right actions an $\A\mbox{-}\A$ bimodule. Denote $V$ with this left module structure by ${}^{L}\!\left(V\right).$ Similarly if one starts with a left $\A$-module $W$, denote the analogously defined right $\A$-module by $\left(W\right)^R.$ Consequently, if $\A$ is supercommutative, the category of left $\A$-modules, the category of right $\A$-modules, and the category of $\A\mbox{-}\A$ bimodules are canonically isomorphic. We hence shall identify all three and denote the abelian category of $\A$-modules by the single category $\A\mbox{-}\Mod.$
\end{defn}

\begin{rem}
We will often abuse notation and denote by the same letter $V$ a left-module, or its associated right-module, or its associated bimodule, when there is no risk of confusion.
\end{rem}

\begin{defn}
Given a left $\A$-module $M,$ the underlying super $\KK$-module $\Pi M$ has the canonical structure of a left $\A$-module as well, by the same formulas defining it for $M.$ More formally, if $$\rho:\A \otimes M \to M$$ is the map exhibiting $M$ as a left module, then $$\A \otimes \Pi M = \A \otimes M \otimes \odd \stackrel{\rho \otimes \odd}{\longlongrightarrow} M \otimes \odd=\Pi M$$ is a map exhibiting $\Pi M$ as a left $\A$-module. Denote the super $\KK$-module $\Pi M$ with this left module structure by $M \Pi.$ Similarly for a right $\A$-module $N$, denote the analogously defined right module structure on the underlying super $\KK$-module $\Pi N$  by $\Pi N$.
\end{defn}

\begin{defn}
If $\A$ is supercommutative, given a left $\A$-module $M,$ there is another way to turn the underlying super $\KK$-module $\Pi M$ into a left $\A$-module, namely by first regarding $M$ as a right $\A$-module $\left(M\right)^R$ and giving $\Pi \left(M\right)^R$ the canonical structure of a \emph{right} $\A$-module, and then regarding this resulting right $\A$-module as a left $\A$-module:
$${}^{L}\!\left(\Pi \left(M\right)^R\right)=:\Pi M.$$
 WARNING: This left $\A$-module structure on $\Pi M$ is \emph{different} than $M \Pi$.
Similarly, denote the analogously constructed right module structure $\Pi N,$ for an $\A$-module $N,$ by $N \Pi.$
\end{defn}

Suppose that $\A$ is supercommutative. Given $$\rho:\A \otimes V \to V$$ and $$\lambda:\A \otimes W \to W$$ two (left) $\A$-modules. Define $$V \mathop{\otimes}_{\A} W:= \varinjlim \left(\mbox{$$\xymatrix@C=3cm{\A \otimes V \otimes W \ar@<+0.65ex>^{\rho \otimes id_W} [r] \ar@<-0.65ex>_{\left(id_V \otimes \lambda\right) \circ \left(\tau_{A,V} \otimes id_W\right)} [r] & V \otimes W}$$}\right),$$
where the colimit is taken in $\sv$ (which may be computed pointwise in $\ve$). Concretely, one may describe $V \mathop{\otimes}_{\A} W$ as having elements which are $\A$-linear combinations of homogeneous elements $v \otimes w,$ such that $$\left(a \cdot v\right) \otimes w=\left(-1\right)^{\deg\left(a\right)\cdot \deg\left(v\right)} v\otimes \left(a \cdot w\right),$$ where each of these homogenous elements $v\otimes w$ has degree $\deg(v)+\deg(w).$ This gives $V \mathop{\otimes}_{\A} W$ the structure of a left $\A$-module via $$a \cdot \left(v\otimes w\right):=\left(a\cdot v\right) \otimes w.$$ The category $\A\mbox{-}\Mod$ of $\A$-modules, together with $\bullet \mathop{\otimes}_{\A} \bullet$ has the structure of a symmetric monoidal category, with the obvious braiding. The unit is given by $\A,$ regarded as module over itself. Moreover, it is symmetric monoidal closed with
$\IHom_{\A\mbox{-}\Mod} \left(V,W\right)_i$ given by the sub $\KK$-module of $\IHom_{\KK\mbox{-}\Mod}\left(V,W[\ii]\right)$ spanned by $\A$-module maps. This super $\KK$-module has the structure of an $\A$-module in the obvious way.

\begin{rem}
Consider the tensor unit $\A$ of $\A\mbox{-}\Mod.$ Notice that $\Pi \A$ is also an $\A$ module. It has the following property:

For any $\A$-module $V,$ $$\Pi \A \mathop{\otimes}_{\A} V \cong \Pi V \cong V \mathop{\otimes}_{\A} \Pi \A.$$
\end{rem}

\begin{rem}
If $\A$ is a commutative $\KK$-algebra, $\iota_!\left(\A\right)=\{\A,0\}$ is a supercommutative $\KK$-algebra and there is a canonical equivalence of categories
$$\iota_!\left(\A\right)\mbox{-}\Mod \simeq \left(\A\mbox{-}\Mod\right)^{\ZZ_2}.$$
\end{rem}

\subsection{Modules as square-zero extensions}\label{sec:sq0}

Recall (Definition \ref{dfn:beckmodule}) for $X$ an object of a category $\bC,$ a \textbf{module} over $X$ is an abelian group object $M$ in the slice category $\bC/X.$ We will show in this subsection that this definition reproduces the notion of a module for a supercommutative algebra as defined in the previous subsection. The argument closely follows \cite{acyclic}.

\begin{defn}
Let $\A$ be a $\bscom_\KK$-algebra, and let $M$ be an $\A$-module. We define the following structure of a supercommutative $\KK$-algebra on the $\KK$-module $\A \oplus M$:
$$\left(a_1,m_1\right) \cdot \left(a_2,m_2\right)=\left(a_1\cdot a_2, a_1\cdot m_2 +m_1\cdot a_2\right).$$ It comes with a canonical projection $$\pi_M:\A \oplus M \to \A$$ which is map of supercommutative algebras. We will write $$\A \oplus M=:M\left[\epsilon\right],$$ and call it the \emph{square zero extension} of $A$ by $M$. Notice that the kernel of $\pi_M$ consists of those elements of the form $\left(0,m\right)$ and is hence canonically isomorphic to $M$ as an $\A$-module. Also $\Ker\left(\pi_M\right)^2=0,$ justifying the terminology ``square zero.'' Moreover, there is a canonical splitting of $\pi_M$ given by $$z_M:a \mapsto \left(a,0\right),$$ which is also an algebra map. Consider the algebra $$M\left[\epsilon\right] \times_\A M\left[\epsilon\right].$$ There is a canonical map
\begin{eqnarray*}
\mu_M:M\left[\epsilon\right] \times_\A M\left[\epsilon\right] &\to& M\left[\epsilon\right]\\
\left(\left(a,m_1\right),\left(a,m_2\right)\right) &\mapsto& \left(a,m_1+m_2\right).
\end{eqnarray*}
It is easy to verify that it is a homomorphism. This gives $\pi_M$ the structure of an abelian group object in the slice category $\bscom_\KK/\A$ with multiplication $\mu_M,$ and unit $z_M.$
\end{defn}

\begin{rem}
Any square zero extension $M\e$ of $\A$ is a split nilpotent extension of $\A$.
\end{rem}

Notice that any element of the form $\left(0,m\right) \in M\e,$ automatically squares to zero. Therefore, it is customary to write $$a+\epsilon m$$ for an element $$\left(a,m\right) \in \M\e \cong \A \oplus M,$$ where $\epsilon$ is to be thought of a formal parameter whose square is zero.

Let $\varphi:M \to N$ be a map of $\A$-modules. Then there is a canonical homomorphism
\begin{eqnarray*}
\varphi[\epsilon]:M\e &\to& N\e\\
a+\epsilon m &\mapsto& a+\epsilon\varphi\left(m\right)
\end{eqnarray*}
which commutes over $\A,$ and $$\varphi \e \circ z_M=z_N.$$ It moreover, by virtue of the additivity of $\varphi,$ it follows that $\varphi[\epsilon]$ is a map of abelian group objects. This defines a functor $$\e:\A\mbox{-}\Mod \to \Ab\left(\bscom_\KK/\A\right)=\Mod\left(\A\right).$$

\begin{prop}\label{prop:modsame}
Let $\KK$ be any super commutative ring, and let $\A$ be a $\KK$-algebra. Then the functor $$\e:\A\mbox{-}\Mod \to \Mod\left(\A\right),$$ is an equivalence of categories.
\end{prop}

\begin{proof}
First, let us show that $\e$ is full and faithful. The fact that $\e$ is faithful is clear. Let $$\theta:M\e \to N\e$$ be a map of abelian group objects. Then since $\theta$ is a morphism over $\A,$ it follows that there is a unique $$\widehat{\theta}:M \to N$$ such that $$\theta\left(a+\epsilon m\right)=a+\epsilon\widehat{\theta}\left(m\right).$$ Since $\theta$ must respect group multiplication, it follows that $\widehat{\theta}$ is a homomorphism of underlying graded abelian groups. Since $\theta$ must be a morphism of supercommutative algebras, the following equality must hold for all $a$ in $\A$ and $m$ in $M$: $$\theta\left(\left(a,0\right)\cdot \left(0,m\right)\right)=\theta\left(\left(a,0\right)\right)\cdot \theta\left(\left(0,m\right)\right).$$ This is equivalent to the condition that $$\widehat{\theta}\left(a\cdot m\right)=a\cdot \widehat{\theta}\left(m\right).$$ This implies $\widehat{\theta}$ is $\A$-linear and hence a map of $\A$-modules. It follows that $\theta=\widehat{\theta}\e,$ so that $\e$ is full and faithful. It remains to show that $\e$ is essentially surjective. Let $\pi:\B \to \A$ be an abelian group object in $\bscom_{\KK}/\A.$ Then it comes with a zero map $$z:\A \to \B$$ which is necessarily a section of $\pi.$ In particular, this means that \emph{as a $\KK$-module}, $\B$ splits as $\B \cong \A \oplus M,$ for some $\KK$-module $M$, and under this identification, $\pi$ is the projection, and $z$ is the map $$a\mapsto \left(a,0\right).$$ Since $\pi$ is a homomorphism, it follows that for $\left(a,m\right), \left(a,m\right)' \in \B,$ multiplication takes the form
\begin{equation}\label{eq:mee}
\left(a,m\right) \cdot \left(a',m'\right) = \left(aa',\varphi\left(a,m,a',m'\right)\right),
\end{equation}
for some function $$\varphi:\A \times M \times \A \times M \to M.$$ Introduce the notation $$a \cdot m:=\varphi\left(a,0,0,m\right)$$ and $$m \cdot a:=\varphi\left(0,a,m,0\right).$$ It is easy to check that this gives $M$ the structure of an $\A\mbox{-}\A$-bimodule. Introduce the notation $$m\cdot m':=\varphi\left(0,m,m',0\right).$$ (This is $\left(0,m\right) \cdot \left(0,m'\right).$) From the distributive law, we can deduce that $$\left(a,m\right) \cdot \left(a',m'\right)=\left(aa',am'+ma'+mm'\right).$$ We claim that $mm'=0$ for all $m$ and $m'$ in $M.$ For this, we will need to use the group multiplication. This is a map $\mu:\B \times_\A \B \to \B.$ It is encoded by a map $$t:\A \times M \times M \to \A \times M$$ such that $$\mu\left(\left(a,m\right),\left(a,m'\right)\right)=t\left(a,m,m'\right).$$ Since $z$ is the group identity, it follows that $$t\left(a,0,m\right)=t\left(a,m,0\right)=\left(a,m\right)$$ for all $a \in \A$ and $m \in M.$ Setting $m=0,$ tells us that for all $a \in \A,$ $$t\left(a,0,0\right)=\left(a,0\right).$$ Notice that $t$ is additive, hence for all $a$ and $m$ we have $$\left(a,m\right)=\left(a,0\right)+t\left(0,0,m\right)=\left(a,0\right)+t\left(0,m,0\right),$$ from which it follows that $$t\left(0,0,m\right)=t\left(0,m,0\right)=\left(0,m\right).$$ Hence, for all $a \in \A,$ and $m,m' \in M,$
\begin{equation}\label{eq:tee}
t\left(a,m,m'\right)=\left(a,m+m'\right).
\end{equation}
Now, since $\mu$ is multiplicative, for all $m$ and $m'$ in $M,$ we have $$\mu\left(\left(\left(1,m\right),\left(1,0\right)\right) \cdot \left(\left(1,0\right),\left(1,m'\right)\right)\right)=\mu\left(\left(\left(1,m\right),\left(1,0\right)\right) \cdot \mu\left(\left(1,0\right),\left(1,m'\right)\right)\right).$$ Using equations (\ref{eq:mee}) and (\ref{eq:tee}), this becomes $$\left(1,m+m'\right)=\left(1+m+m'+mm'\right).$$ It follows that $mm'=0.$ Hence, $\pi:\B \to \A$ is isomorphic to the square zero extension $M\e$ associated to the module $M$.
\end{proof}

\section{Comonoid objects and their coactions.}\label{sec:coactions}
Recall that a \emph{comonoid} (resp. \emph{cogroup}) object in a category $\mathbf{C}$ with finite coproducts is a monoid (resp. group) object in $\mathbf{C}^{\mathrm{op}}$, i.e. an object $\H$ of $\mathbf{C}$ together with morphisms
\[
\Delta:\H\To\H\amalg\H,\quad\epsilon:\H\To\varnothing,
\]
called \emph{comultiplication} and \emph{counit}, respectively (and, in the cogroup case, also $$S:\H\To\H,$$ called the \emph{antipode}), satisfying the standard equations for coassociativity, counitality and the antipode, dual to those for a group:
\[
(\Delta\amalg\mathrm{id})\circ\Delta=(\mathrm{id}\amalg\Delta)\circ\Delta,\quad(\epsilon\amalg\mathrm{id})\circ\Delta
=(\mathrm{id}\amalg\epsilon)\circ\Delta,
\]
\[
\nabla\circ(S\amalg\mathrm{id})\circ\Delta=e\circ\epsilon=\nabla\circ(\mathrm{id}\amalg S)\circ\Delta,
\]
where $\nabla:\H\amalg\H\To\H$ is the codiagonal and $e:\varnothing\To\H$ is the initial object inclusion. In other words, $\H$ is a commutative bialgebra (resp. Hopf algebra) object in the symmetric monoidal category $(\mathbf{C},\amalg,\varnothing)$, with $\nabla$ playing the role of multiplication, and $e$ that of the unit. One says further that $\H$ is \emph{cocommutative} if
\[
\tau\circ\Delta=\Delta,
\]
where $\tau:\H\amalg\H\To\H\amalg\H$ is the transposition.

A (left) \emph{coaction} of a comonoid (or cogroup) object $\H$ on an object $\C$ of $\mathbf{C}$ is a morphism
\[
\Phi:\C\To\H\amalg\C
\]
such that
\[
(\Delta\amalg\mathrm{id})\circ\Phi=(\mathrm{id}\amalg\Phi)\circ\Phi,\quad(\epsilon\amalg\mathrm{id})\circ\Phi=\mathrm{id}
\]
(and similarly for right coactions). Denote the category of objects of $\bC$ with an $\H$-coaction and $\H$-equivariant morphisms between them by $\mathbf{C}^\H$.

The (covariant) Yoneda embedding
\[
Y_{\bC^{\op}}:\mathbf{C}^\mathrm{op}\To\Set^\mathbf{C}
\]
takes every comonoid (resp. cogroup) object $\H\in\mathbf{C}$ to a monoid (resp. group) valued functor $Y_{\bC^{\op}}(\H)=\mathbf{C}(\H,-)$, i.e. a monoid (resp. group) object in $\Set^\mathbf{C}$; a coaction of $\H$ on $\C$ gives rise to an action of $Y_{\bC^{\op}}(\H)$ on $Y_{\bC^{\op}}(\C)$.

\begin{rem}\label{rem:cofreecoact}
The forgetful functor $\bC^\H\to\bC$ (forgetting the coaction) has a \emph{right} adjoint assigning to a $\C\in\bC$ the object $\H\amalg\C$ with the coaction given by
\[
\Delta\amalg\id_\C:\H\amalg\C\To\H\amalg\H\amalg\C.
\]
Colimits in $\bC^\H$ are created in $\bC$; explicitly,
\[
(\underrightarrow{\lim}\,\C_\alpha,\underrightarrow{\lim}\,\Phi_\alpha:\underrightarrow{\lim}\,\C_\alpha\to
\underrightarrow{\lim}\,(\H\amalg\C_\alpha)\cong
\H\amalg(\underrightarrow{\lim}\,\C_\alpha))
\]
is a colimit in $\C^\H$. The same holds for those \emph{limits} in $\bC$ that are preserved by $\H\amalg(\quad)$.
\end{rem}

\begin{rem}\label{rem:freecoact}
Furthermore, if $\H$ is \emph{co-exponentiable} (exponentiable in $\bC^\op$), the forgetful functor $\bC^\H\to\bC$ also has a \emph{left} adjoint. Indeed, let $\C\in\bC$ be any object and $\cop{\C}{\H}$ the co-exponential by $\H$. The co-action
\[
\cop{\C}{\H}\To\H\amalg\cop{\C}{\H}
\]
is the image of the map
\[
\Delta^*:\cop{\C}{(\H\amalg\H)}\To\cop{\C}{\H}
\]
under the composite of the natural isomorphisms
\[
\bC(\cop{\C}{(\H\amalg\H)},\cop{\C}{\H})\To\bC(\cop{({\cop{\C}{\H}})}{\H},\cop{\C}{\H})\To\bC(\cop{\C}{\H},\H\amalg\cop{\C}{\H}),
\]
where $\Delta:\H\to\H\amalg\H$ is the comultiplication. We leave it to the reader to verify that this defines a coaction, and the requisite adjointness.
\end{rem}

Now suppose that we are given categories $\bC$ and $\bD$ with finite limits and colimits and an adjunction
\begin{equation}\label{eq:algmor}
\Adj{R}{\bC}{\bD}{L},
\end{equation}
with $L$ left adjoint to $R.$ Let $\otimes$ denote the coproduct in $\bC$ and $\odot$ the one in $\bD$. Let $\H$ be a comonoid object in $\bC$. Then, since $L$ respects coproducts, $\hat\H=L(\H)$ is a comonoid object in $\bD$. Moreover, if $\Phi:\C\to\H\otimes\C$ is a coaction of $\H$ on $\C$, then $L(\Phi):L(\C)\to\hat\H\odot L(\C)$ is a coaction of $\hat\H$ on $L(\C)$. Observe that, if $H=Y_{\bC^\op}(\H)$, then $\hat\H=Y_{\bD^\op}(\hat\H)=H\circ R$. Denote the functor sending a coaction $(\C,\Phi)$ in $\bC^\H$ to $(L(\C),L(\Phi))\in\bD^{\hat\H}$ by $L^\H$.

\begin{prop}\label{prop:indadj} The functor $L^\H$ has a right adjoint
\[
R_\H:\bD^{\hat{\H}}\To\bC^{\H}.
\]
\end{prop}

\begin{proof}
Given an object $(\D,\Phi:\D\to\hat\H\odot\D)\in\bD^{\hat\H}$ of $\bD^{\hat{\H}}$, define $\D_\mathrm{alg}\in\bC$
to be the following pullback:
$$\xymatrix@C=1.5cm{\D_\mathrm{alg} \ar[r] \ar[d] & \H\otimes R\left(\D\right) \ar[d]\\
R\left(\D\right) \ar[r]_-{R\left(\Phi\right)} & R\left(\hat\H\odot\D\right).}$$
Suppose we are given $(\C,\Phi_\C)\in\bC^\H$, $(\D,\Phi_\D)\in\bD^{\hat{\H}}$ and a map $\psi:L(\C)\To\D$ such that
$$\xymatrix@C=1.5cm{L\left(\C\right) \ar[r]^-{\psi} \ar[d]_-{L\left(\Phi_\C\right)} & \D \ar[d]^-{\Phi_\D}\\
\hat\H\odot L\left(\C\right) \ar[r]_-{\mathrm{id}\odot\psi} & \hat\H\odot\D}$$
commutes. By adjunction, this is equivalent to the commutativity of
$$\xymatrix@C=1.5cm{\C \ar[r]^-{\tilde\psi} \ar[d]_-{\Phi_\C } & R\left(\D\right) \ar[d]^-{\Phi_\C } \\
\H\otimes \C \ar[r]_-{\widetilde{\mathrm{id}\odot\psi}} & R\left(\hat\H\odot\D\right).}$$
Using the universal properties, one sees that $\widetilde{\mathrm{id}\odot\psi}$ factors through
\[
\mathrm{id}\otimes\tilde\psi:\H\otimes\C\To\H\otimes R(\D),
\]
hence, by the commutativity of the diagram, $\tilde\psi$ factors through $\D_\mathrm{alg}$ and furthermore, the restriction of $R(\Phi_\D)$ to $\D_\mathrm{alg}$ factors through $\H\otimes\D_\mathrm{alg}$, thereby defining the coaction map
\[
\Phi_{\D_\mathrm{alg}}=(\Phi_\D)_\mathrm{alg}:\D_\mathrm{alg}\To\H\otimes\D_\mathrm{alg}
\]
making the diagram
$$\xymatrix@C=1.5cm{\C \ar[r]^-{\tilde\psi} \ar[d]_-{\Phi_\C} & \D \ar[d]^-{\Phi_{\D_\mathrm{alg}}}\\
\H\otimes\C \ar[r]_-{\mathrm{id}\otimes\tilde\psi} & \H\otimes\D_\mathrm{alg}}$$
commute. Thus, the assignment $(\D,\Phi_\D)\mapsto(\D_\alg,\Phi_{\D_\alg})$ defines a functor
\[
R_\H:\bD^{\hat\H}\To\bC^\H
\]
which is right adjoint to $L^\H$.
\end{proof}

\begin{defn}
Suppose $\bC$ is a category with finite limits and colimits. Call a coaction $\Phi:\C\To\H\amalg\C$ \emph{trivial} if $\Phi=j_2$, the canonical inclusion of $\C$ into the coproduct; more generally, for any coaction $(\C,\Phi)$ the \emph{invariant subobject} is defined as the equalizer
\[
\C_\H\To\C\rightrightarrows\H\amalg\C
\]
of $\Phi$ and $j_2$.
\end{defn}

The assignment $\C\mapsto(\C,j_2:\C\to\H\amalg\C)$ defines a fully faithful functor
\[
(\quad)_\triv:\bC\To\bC^\H,
\]
which has a right adjoint, namely the functor
\[
(\quad)_\H:\bC^\H\To\bC
\]
associating to a coaction $(\C,\Phi)$ its invariant subobject $\C_\H$. In particular, the map $\C_\H\to\C$ is $\H$-equivariant with respect to the trivial coaction on $\C_\H$.

The functor $(\quad)_\triv$ also has a left adjoint, assigning to an $\H$-coaction its coinvariant quotient $\A^\H$.

\begin{rem}\label{rem:triv} For $\bC$ and $\bD$ as above and a comonoid $\H\in\bC$, the diagram
$$\xymatrix@C=1.5cm{\bD^{\hat\H} \ar[r]^-{(\quad)_{\hat\H}} \ar[d]_-{R_\H} & \bD \ar[d]^-{R}\\
\bC^\H \ar[r]_-{(\quad)_\H} & \bC}$$
commutes up to natural isomorphism; in particular, for a trivial coaction on $\D\in\bD$, $\D_\mathrm{alg}$ is isomorphic to $R(\D)$.
\end{rem}
In contrast with this, the left adjoint $L^\H$ generally fails to preserve invariant subobjects. On the other hand, it maps the trivial $\H$-coaction on $\C\in\bC$ to the trivial $\hat\H$-coaction on $L(\C)$. In case the adjunction \eqref{eq:algmor} is an algebraic morphism between algebraic theories (cf. \cite{dg1}, Appendix), there is a modification of the above construction for which the left adjoint leaves the invariant subobjects alone. The idea is to do everything relative to the invariant subobjects. Observe first that, for each $\D\in\bD$, we have an induced adjunction
\[
\Adj{R_\D}{R(\D)/\bC}{\D/\bD}{L_\D},
\]
where the right adjoint $R_\D$ sends each $f:\D\to\D'$ in $\D/\bD$ to $R(f):R(\D)\to R(\D')$ in $R(\D)/\bC$, while the left adjoint $L_\D$ has the property that, for each $\C\in\bC$, it sends the canonical map $R(\D)\to\C\otimes R(\D)$ to the canonical map $\D\to L(\C)\odot\D$ (see \cite{dg1}, Appendix).

Now, by Remark \ref{rem:triv}, we have a functor
\[
R_\H^\dagger:\bD^{\hat\H}\To\bC^\H\times_{\bC}\bD.
\]
Here, $\bC^\H\times_{\bC}\bD$ denotes the homotopy pullback in the $(2,1)$-category of categories, functors and natural isomorphisms. Its objects are triples $((\C,\Phi),\D^0,\phi)$, where $(\C,\Phi)\in\bC^\H$, $\D^0\in\bD$ and $\phi:R(\D^0)\to\C_\H$ is an isomorphism; the morphisms are the appropriate commutative diagrams. The functor $R_\H^\dagger$ assigns to a $$(\D,\Psi)\in\bD^{\hat\H}$$ the triple $$(R_\H(\D,\Psi),\D_{\hat\H},\phi_\D),$$ where $\phi_\D:R(\D_{\hat\H})\to(\D_\alg)_\H$ is the natural isomorphism from Remark \ref{rem:triv}.

\begin{prop}\label{prop:indadjrel}
If the adjunction \eqref{eq:algmor} is an algebraic morphism of algebraic theories, the functor $R_\H^\dagger$ has a left adjoint $L^\H_\dagger$.
\end{prop}
\begin{proof}
Suppose we are given a triple $((\C,\Phi),\D^0,\phi)\in\bC^\H\times_{\bC}\bD$. Observe that, since $\C_\H$ equalizes $\Phi$ and $j_2$, $\Phi:\C\to\H\otimes\C$ extends to a morphism in $\C_\H/\bC$ from the inclusion $i:\C_\H\to\C$ to $j_2\circ i:\C_\H\to\H\otimes\C$. Writing
\[
\H\otimes\C\cong(\H\otimes\C_\H)\otimes_{\C_\H}\C
\]
and using $\phi$ to identify $\C_\H$ with $R(\D^0)$, we get a map
\[
\Phi':\C\To(\H\otimes R(\D^0))\otimes_{R(\D^0)}\C
\]
in $R(\D^0)/\bC.$ (In fact, $\Phi'$ is a coaction of the comonoid object $R(\D^0)\to\H\otimes R(\D^0)$ on $R(\D^0)\to\C$ in $R(\D^0)/\bC.$) Applying $L_{\D^0}$ to this map we get a map
\[
L_{\D^0}(\Phi'):L_{\D^0}(\C)\To(\hat\H\odot\D^0)\odot_{\D^0}L_{\D^0}(\C)\cong\hat\H\odot L_{\D^0}(\C)
\]
(over $\D^0$). Since $L_{\D^0}$ is left adjoint, this defines a coaction of $\hat\H$ on $L_{\D^0}(\C)$ in $\D^0/\bD$. Finally, define $L^\H_\dagger((\C,\Phi),\D^0,\phi)$ to be $(L_{D^0}(\C),L_{\D^0}(\Phi'))\in\bD^{\hat\H}$ (forgetting the map out of $\D^0$). One can verify that this defines the desired left adjoint.
\end{proof}

Finally, there is one further adjunction
\begin{equation}\label{eq:adjcirc}
\Adj{R_\H^\circ}{\bC^\H}{\bC^\H\times_\bC\bD}{L^\H_\circ}
\end{equation}
where
\[
R_\H^\circ((\C,\Phi),\D^0,\phi)=(\C,\Phi),
\]
while
\[
L^\H_\circ(\C,\Phi)=(RL(\C_\H)\otimes_{\C_\H}\C,L\C_\H,\phi:RL\C_\H\to(RL(\C_\H)\otimes_{\C_H}\C)_\H).
\]
Here, the pushout is taken in $\bC^\H$ (cf. Remark \ref{rem:cofreecoact}), where $\C_\H$ and $RL(\C_\H)$ are endowed with trivial coactions, along the maps
\[
u_\triv:(\C_\H)_\triv\To(RL(\C_\H))_\triv
\]
(where $u:\C_\H\to RL(\C_\H))$ is the unit of \eqref{eq:algmor}) and the counit
\[
(\C_\H)_\triv\To\C.
\]
The composition of \eqref{eq:adjcirc} with
\begin{equation}\label{eq:adjdagger}
\Adj{R_\H^\dagger}{\bC^\H\times_\bC\bD}{\bD^{\hat\H}}{L^\H_\dagger}
\end{equation}
is
\begin{equation}\label{eq:adjH}
\Adj{R_\H}{\bC^\H}{\bD^{\hat\H}}{L^\H}.
\end{equation}
Let us conclude by describing coproducts in $\bC^\H\times_\bC\bD$. First, given a coaction $(\C,\Phi)\in\bC^\H$, and a map $\C_\H\to\Q$ in $\bC$, we can form the pushout $\P$ of
\[
\Q_\triv\oT(\C_\H)_\triv\To\C
\]
in $\bC$ (``change of base''), with $\P_\H\simeq\Q$. Now suppose we are given an
\[
((\A_i,\Phi_i), \A_i^0,\phi_i:(\A_i^0)_\sharp\to(\A_i)_\H),
\]
$i=1,2$, in $\bC^\H\times_\bC\bD$. Define first $\Q$ to be the pushout of the natural maps
\[
(\A_1^0\odot\A_2^0)_\sharp\oT(\A_1)_\H\otimes(\A_2)_\H\To(\A_1\otimes\A_2)_\H.
\]
Then define the coaction $(\A_1\circledast\A_2,\Phi)$ to be the pushout of
\[
((\hat\Q)_\sharp)_\triv\oT\Q_\triv\oT((\A_1\otimes\A_2)_\H)_\triv\To\A_1\otimes\A_2
\]
in $\bC^\H$ (where $\hat\Q=\tau^{\A_1^0\odot\A_2^0}_!(\Q)$, the completion relative to $\A_1^0\odot\A_2^0$). Finally,
\[
((\A_1\circledast\A_2,\Phi),\hat\Q,\phi:(\hat\Q)_\sharp\to(\A_1\circledast\A_2)_\H)
\]
is the desired coproduct.

% ----------------------------------------------------------------
\bibliographystyle{hplain}
\bibliography{derived}
\end{document}